\documentclass[11pt]{article}

\usepackage{a4wide}
\usepackage{amsmath}
\usepackage{amssymb}
\usepackage{dsfont}
\usepackage{verbatim}
\usepackage{color}
\usepackage{graphicx}
\usepackage{caption}
\usepackage{float}
\usepackage{subfigure}
\usepackage{ifpdf}

\newcommand\indicator{\mathds{1}}
\newcommand\prob{\mathds{P}}
\newcommand{\md}{\mathrm{d}}
\newcommand{\me}{\mathrm{e}}
\newcommand{\hf}{\hat f}
\newcommand{\CM}{\text{CM}}
\newcommand{\sign}{\text{sign}}
\newcommand{\R}{\mathbb{R}}
\newcommand{\E}{\mathbb{E}}
\def\argmax{\mathop{\mbox{\sl\em argmax}}}
\newcommand{\indep}{(A0)}
\newcommand{\mono}{(A1)}
\newcommand{\monostar}{($\text{A}^\star$)}
\newcommand{\mart}{(A2)}
\newcommand{\embed}{(A3)}
\newcommand{\regul}{(A4)}
\newcommand{\tqed}{\hfill{$\Box$}\\}
\renewcommand{\hat}{\widehat}
\newcommand{\holder}{\theta}
\newcommand{\hn}{\gamma}
\newcommand{\IK}{\mathbb{K}}
\newcommand{\bcoefA}{\phi}
\newcommand{\bcoefB}{\psi}
\newcommand{\bK}{K_{B,t}}
\newcommand{\bIK}{\mathbb{K}_{B,t}}
\newtheorem{lemma}{Lemma}
\newtheorem{theorem}{Theorem}
\newtheorem{remark}{Theorem}

\renewcommand{\P}{\mathbb{P}}

\begin{document}


\markboth{C\'ecile Durot, Piet Groeneboom, and Hendrik P.~Lopuha\"a}{Journal of Nonparametric Statistics}


\title{Testing equality of functions under monotonicity constraints}

\author{C\'ecile Durot, Piet Groeneboom and Hendrik P.~Lopuha\"a}

\maketitle

\begin{abstract}
We consider the problem of testing equality of functions $f_j:[a,b]\to \mathbb{R}$
for $j=1,2,\ldots,J$ on the basis of $J$ independent samples from possibly
different distributions under the assumption that the functions are monotone.
We provide a uniform approach that covers testing equality of
monotone regression curves, equality of monotone densities and equality of
monotone hazards in the random censorship model. Two test statistics
are proposed based on $L_1$-distances. We show that both statistics
are asymptotically normal and we provide bootstrap implementations, which are
shown to have critical regions with asymptotic level $\alpha$.
\end{abstract}

\section{Introduction}
\label{sec:introduction}
A classical statistical problem is the $k$-sample problem, where one has to decide
whether different samples can be regarded as coming from the same population.
In the non-parametric setting,~\cite{smirnov1939} and~\cite{wilcoxon1947} introduced the first two-sample tests,
one based on the distance between the empirical distribution functions of the two samples
and one based on ranks.
Generalizations of these methods to the $k$-sample problem, with $k\geq 2$,
are given in~\cite{kiefer1959} and~\cite{kruskal-wallis1952}, respectively.
Although Wilcoxon-type tests remain very popular, they are known to be able to detect only a limited range of alternatives.
To overcome this difficulty, several solutions have been proposed.
For an overview and references on this topic, see~\cite{janic-ledwina2000}, who developed a two-sample test inspired by the smooth Neymann test,
in which the problem is reparametrized and modeled via some multivariate exponential family
with an unknown parameter in such a way that the two-sample problem amounts to testing that this parameter is zero.
Another recent approach consist of comparing empirical characteristic functions,
e.g., see~\cite{huskova-meintanis2008}, who extended the univariate two-sample tests from~\cite{epps-singleton1986} and~\cite{meintanis2005}
to the multivariate $k$-sample setting, $k\geq 2$.
For procedures based on comparing kernel density estimators, see Anderson, Hall and Titterington~\cite{anderson-hall-titterington1994},
who consider the $L_2$-distance in the multivariate two sample setting, and~\cite{louani2000}, who uses the $L_1$- and
the $L_{\infty}$-distance, and see~\cite{cao-vankeilegom2006} and references therein, for methods based on empirical likelihood.
A generalization of the Smirnov test to the case of multivariate observations is considered in~\cite{bickel1968,praestgaard1995}, among others.
See also~\cite{baringhaus2004,aslan-zech2005,jureckova-kalina2012} for more references on multivariate $k$-sample tests.
The recent paper by~\cite{jureckova-kalina2012} points out that most of the rank based tests proposed in the literature,
unless one-sided, could be biased against alternatives of interest.
In that paper, multivariate distribution free two-sample tests,
based on the ranks of suitable distances of multivariate observations, are developed.
Unbiasedness and power of those tests are studied against Lehmann alternatives.

The $k$-sample problem arises naturally in survival analysis, where the observations are typically right censored.
Various two-sample tests inspired by the Wilcoxon test have been proposed for right censored data,
see~\cite{harrington-fleming1982,jones-crowley1989, shen-le2000} and references therein.
Other approaches are based for instance on comparison of quantile estimators, see Li, Tiwari and Welles~\cite{li-tiwari-wells1996}.

The $k$-sample problem also arises in the regression setting.
For instance, in medical studies
one wishes to compare the mean response of a treatment group of subjects to that of a control group,
taking into account a covariate such as the dose of drug.
In such cases, one wishes to compare two or more regression curves.
In this setting, \cite{hall-hart1990} developed a two-sample test calibrated by the bootstrap.
The test statistic, which is essentially a scaled version of an integrated $L_2$-distance between smooth estimators of
the regression curves (see also King, Hart and Wehrly~\cite{king-hart-wehrly1991}),
is based on differences between the response variables at given values of the covariate.
Another two-sample test based on these differences,
which in spirit resembles the Kolmogorov-Smirnov statistic, is proposed in~\cite{delgado1993}.
In the case where the covariate values as well as the sample size possibly differ in the two considered samples,
\cite{kulasekara1995} developed a test based on quasi-residuals.
In this setting with possibly heteroscedastic errors,
\cite{munk-dette1998} consider a test based on the estimation of the $L_{2}$-distance between the two regression curves
and generalized the method to the case of $k$-samples, $k\geq 2$,
whereas \cite{neumeyer-dette2003} developed a test calibrated by the bootstrap,
which is based on the difference of two marked empirical processes.
Most of the aforementioned procedures involve  the choice of a tuning parameter.

Nonparametric methods for $k$-sample problems under shape constraints are quite limited.
\cite{groeneboom2012} considers likelihood ratio type two-sample tests
in the current status model, which is closely related to other shape constrained nonparametric models.
The test statistics are shown to be asymptotically normal under the null hypothesis of equal distributions
and the test is calibrated using the bootstrap.
The problem of comparing two monotone fractile regression functions,
which is similar in spirit to the problem of comparing two monotone regression functions, is considered by~\cite{sen2010}.
A test based on the $L_{2}$-distance between two monotone estimators is discussed.
The test is calibrated by the bootstrap, but no limit distribution is provided.

The above mentioned testing problems have one common feature, i.e.,
they essentially test the equality of two or more functions in various sorts of statistical models,
e.g., distribution functions, densities, characteristic functions, hazard rates, or regression curves.
In this paper we consider the following testing problem
\[
H_0:f_1=f_2=\cdots=f_J
\quad\text{ against }\quad
H_1: f_i\ne f_j\text{ for some }i\ne j
\]
subject to the constraint that $f_j:[a,b]\mapsto\mathbb{R}$ is decreasing for all $j=1,2,\ldots,J$,
where $a,b\in \mathbb{R}$ are known.
This general framework includes a $k$-sample test for samples from a decreasing density and a test for equality
of decreasing regression curves or decreasing hazard rates.
We propose two test statistics based on $L_1$-distances between non-parametric Grenander-type estimators $\widehat{f}_{nj}$.
The first one compares mutual distances between the different individual estimators
\begin{equation}
\label{eq:defS1}
S_{n1}=\sum_{i<j}\int_a^b|\widehat f_{ni}(t)-\widehat f_{nj}(t)|\,\md t.
\end{equation}
The second one compares the distances between each individual $\widehat{f}_{nj}$ and a pooled estimator~$\widehat{f}_{n0}$
for the mutual $f_0$ under the null hypothesis:
\begin{equation}
\label{eq:defS2}
S_{n2}=\sum_{j=1}^J\int_a^b|\widehat f_{nj}(t)-\widehat f_{n0}(t)|\,\md t.
\end{equation}
We show that both test statistics are asymptotically normal
and propose a bootstrap procedure to calibrate the test.
Finally, we discuss the special cases of testing equality of monotone regression functions,
monotone densities and monotone hazard rates under random censorship, and show in each of these cases that the bootstrap works.

In Section~\ref{sec:main} we describe the general setup, state our main results,
and discuss the differences between the approach in this paper and the ones
used in Groeneboom, Hooghiemstra and Lopuha\"a~\cite{groeneboom-hooghiemstra-lopuhaa1999},~\cite{kulikov-lopuhaa2005} and~\cite{durot2007}
to prove similar results.
In Section~\ref{sec:calibration} we discuss the calibration of the test and show that the bootstrap works
in the previous mentioned statistical models.
All main proofs are postponed to an appendix at the end of the paper, and details are given in the supplement.

\section{Main results}
\label{sec:main}
For estimating the functions $f_j:[a,b]\to \mathbb{R}$, we suppose that
for each $j=1,2,\ldots,J$, and $t\in[a,b]$,
\begin{equation}
\label{eq:def Fj}
F_j(t)=\int_a^{t} f_j(x)\,\md x,
\end{equation}
is well defined and that we have an estimator $F_{nj}$ at hand based on $n_j$ observations.
We denote by $n=\sum_{j}n_{j}$ the total number of observations and for notational convenience (and possibly arguing along subsequences), we assume that $n_j=c_j\cdot n$ where  $c_j>0$ does not depend on $n$.
Thus,  $\sum_jc_j=1$.
Denote $f_0$ and $F_0$ for the corresponding quantities under the null hypothesis, where we estimate $F_0$ by
\begin{equation}
\label{eq:def Fn0}
F_{n0}=\sum_{j=1}^J c_j F_{nj}.
\end{equation}
Then, for all $j=0,1\ldots,J$, we define $\widehat f_{nj}$ as the left-hand slope of the least concave majorant of $F_{nj}$,
with $\widehat f_{nj}(a)=\lim_{t\downarrow a}\widehat f_{nj}(t)$.
We will frequently make use of results from~\cite{groeneboom-hooghiemstra-lopuhaa1999} and~\cite{durot2007}.
Similar to these papers, we will work under the following assumptions:
\begin{enumerate}
\item[\indep]
The estimators $F_{n1},F_{n2},\ldots,F_{nJ}$ are independent and for every $j=1,2,\ldots,J$, the estimator
$F_{nj}:[a,b]\to\R$ is a cadlag step process.
\item[\mono]
For each $j=1,2,\ldots,J$, the function $f_j:[a,b]\mapsto \mathbb{R}$ is decreasing and continuously differentiable,
such that $0<\inf_{t\in[a,b]}|f_j'(t)|\leq \sup_{t\in[a,b]}|f_j'(t)|<\infty$.
\item[\mart]
For each $j=1,2,\ldots,J$, there exists a constant $C_j>0$, such that for all $x\geq 0$ and $t=a,b$,
the process $M_{nj}=F_{nj}-F_j$ satisfies
\begin{equation*}
\mathds{E}
\left[
\sup_{u\in[a,b],\,x/2\leq |t-u|\leq x}
\left(
M_{nj}(t)-M_{nj}(u)
\right)^2
\right]
\leq
\frac{C_jx}{n_j}.
\end{equation*}
\end{enumerate}
Furthermore, we assume that there exists an embedding either into Brownian motion or into Brownian bridge.
\begin{enumerate}
\item[\embed]
For each $j=1,2,\ldots,J$, there exists a Brownian motion or Brownian bridge $B_{nj}$,
an increasing function $L_{j}:[a,b]\mapsto\mathbb{R}$ with $\inf_{t\in[a,b]}L_j'(t)>0$, and constants $q>6$ and $C>0$,
such that for all $x\in(0,n_j]$:
\[
\prob\left\{
n_j^{1-1/q}
\sup_{t\in[a,b]}
\left|
M_{nj}(t)-n_j^{-1/2}B_{nj}\circ L_{j}(t)
\right|
>x
\right\}
\leq
Cx^{-q}.
\]
\end{enumerate}
It should be noticed that, since the $F_{nj}$ are assumed to be independent,
we can assume without loss of generality that the $B_{nj}$ are independent.
Note that, for $j=1,2,\ldots,J$, we can write
\begin{equation}
\label{eq:Bnj and Wnj}
B_{nj}(t)=W_{nj}(t)-\xi_{nj}t,
\qquad
\text{for }t\in[a,b],
\end{equation}
where the $W_{nj}$ are independent Brownian motions and
$\xi_{nj}\equiv0$, if $B_{nj}$ is Brownian motion,
and $\xi_{nj}\sim N(0,1)$ independent of $B_{nj}$,
if $B_{nj}$ is Brownian bridge.
Finally, we require the following smoothness assumption.
\begin{enumerate}
\item[\regul]
There exist a $\holder\in(3/4,1]$ and $C>0$, such that for all $x,y\in[a,b]$ and $j=1,2,\ldots,J$,
\[
|f_j'(x)-f_j'(y)|\leq C|x-y|^\holder\text{ and }|L_j''(x)-L_j''(y)|\leq C|x-y|^\holder.
\]
\end{enumerate}
These are the usual assumptions when studying the $L_p$-error of isotonic estimators.
It is explained in~\cite{durot2007} that several classical models are covered by the above general framework.
As an example we mention the model where one observes $J$ independent samples where the random variables
in the $j$th sample have a decreasing smooth density function $f_j$.
In this example, $F_{nj}$ is the empirical distribution function based on the $j$th sample and $L_j=F_j$.
This example will be studied in detail in Section~\ref{subsec:density}.
Two other examples, where $f_j$ is either a regression function or a failure rate, are studied in
Sections~\ref{subsec:regression} and~\ref{subsec:censoring}.

Our main result is that under the above assumptions both test statistics defined in~\eqref{eq:defS1}
and~\eqref{eq:defS2} are asymptotically Gaussian under the null hypothesis.
In order to formulate these results more precisely, we introduce the random variables
\begin{equation}
\label{eq: def zetaj}
\zeta_j(c)
=
\argmax_{u\in\R}
\left\{
W_j(u+c)-u^2
\right\},
\quad
\text{ for }
c\in\R
\text{ and }
j=1,2,\ldots,J,
\end{equation}
where the argmax function is the supremum of the times at which the maximum is attained,
$W_1,W_2,\ldots,W_J$ are independent standard two-sided Brownian motions.
We are now in the position to establish asymptotic normality for test statistic $S_{n1}$.
\begin{theorem}
\label{th:main1}
Assume  \indep, \mono, \mart, \embed, \regul\
and let $S_{n1}$ be defined by~\eqref{eq:defS1}.
Let $\zeta_j$ be defined in~\eqref{eq: def zetaj}, for $j=1,2,\ldots,J$, with independent standard Brownian motions $W_1,W_2,\ldots,W_J$.
If $f_0=f_1=\cdots=f_J$, then $n^{1/6}(n^{1/3}S_{n1}-m_1)$
converges in distribution, as $n\to\infty$, to the Gaussian law with mean zero and variance $\sigma_1^2$,
where
\[
\begin{split}
\sigma^2_1
&=
8
\sum_{i<j}
\sum_{l<m}
\int_a^b\int_0^{\infty}
\mathrm{cov}\left(|Y_{si}(t)-Y_{sj}(t)|, |Y_{sl}(0)-Y_{sm}(0)|\right)\,\md t\,\md s,\\
m_1
&=
\sum_{i<j}
\int_{a}^{b}
|4f_{0}'(s)|^{1/3}
\E\left|
Y_{si}(0)-Y_{sj}(0)
\right|
\,\md s,
\end{split}
\]
with
\begin{equation}
\label{eq:def Ysj}
Y_{sj}(t)
=
\frac{L_{j}'(s)^{1/3}}{c_j^{1/3}}
\zeta_{j}\left(\frac{c_j^{1/3}t}{L_{j}'(s)^{1/3}}\right),
\quad
\text{for $j=1,2,\ldots,J$.}
\end{equation}
\end{theorem}
In the special case where all $L_j=L$ are the same, after change of variables $u=t/L(s)^{1/3}$, we find
\[
\begin{split}
\sigma^2_1
&=
8\int_a^b L'(s)\,\md s
\sum_{i<j}
\sum_{l<m}
\int_0^\infty
\mathrm{cov}
\left(
\left|\frac{\zeta_i(c_i^{1/3}u)}{c_i^{1/3}}-\frac{\zeta_j(c_j^{1/3}u)}{c_j^{1/3}}\right|,
\left|\frac{\zeta_l(0)}{c_l^{1/3}}-\frac{\zeta_m(0)}{c_m^{1/3}}\right|
\right)
\,\md u,\\
m_1
&=
\int_{a}^{b}
|4f_0'(s)L'(s)|^{1/3}
\,\md s\sum_{i<j}
\E\left|\frac{\zeta_i(0)}{c_i^{1/3}}-\frac{\zeta_j(0)}{c_j^{1/3}}\right|.
\end{split}
\]
This resembles the normalizing constants in Theorem 2 in~\cite{durot2007} for the case
$p=1$.
An example is the monotone density model, where under the null hypothesis $L'_j=f_j=f_0$,
in which case~$m_1$ and $\sigma^2_1$ coincide with the normalizing constants in Theorem~1.1 in~\cite{groeneboom-hooghiemstra-lopuhaa1999}.
In particular, the limiting variance $\sigma_1^2$ does not depend on the underlying density $f_0$.
Another example, where the limiting variance does not depend on $f_{0}$, is the monotone regression model, where $L_j(t)=(t-a)(b-a)\tau^2_j$, with $\tau_j^2$ being the variance of the measurement error.

To establish a similar result for $S_{n2}$ is more complex, due to the fact that we have to deal with
differences between a marginal estimator $\widehat{f}_{ni}$ and the pooled estimator $\widehat{f}_{n0}$,
which are both (partly) build from the same $i$th sample.
First of all, we need that,
under the null hypothesis $f_0=f_1=\cdots=f_J$, the above assumptions also hold to some extent for~$j=0$.
Clearly, assumption~\indep\ then becomes irrelevant and~\mono\ is immediate, as well as the first part of~\regul.
Because, under the null hypothesis,
\begin{equation}
\label{eq:relation Mn0}
M_{n0}(t)=F_{n0}(t)-F_0(t)=\sum_{j=1}^J c_j(F_{nj}(t)-F_j(t))=\sum_{j=1}^J c_jM_{nj}(t),
\end{equation}
the inequality (Jensen)
\begin{equation}
\label{eq:jensen}
(a_1+\cdots+a_k)^2\leq k(a_1^2+\cdots+a_k^2),
\end{equation}
yields that \mart\, also holds for $j=0$.
The remaining assumptions require the definitions of~$L_0$ and $B_{n0}$.
To have~\embed\ for $j=0$, we need to show that $M_{n0}$ can be approximated by~$n^{-1/2}B_{n0}\circ L_0$,
for a given increasing function $L_0:[a,b]\to\mathbb{R}$ with $\inf_{t\in[a,b]}L_0'(t)>0$
and a Gaussian process $B_{n0}$.
From relation~\eqref{eq:relation Mn0} and assumption~\embed, for $j=1,2,\ldots,J$,
it follows that we then must have
\begin{equation}
\label{eq:defBn0}
B_{n0}(t)
=
\sum_{j=1}^Jc_j^{1/2}B_{nj}\circ L_j\circ L_0^{-1}(t).
\end{equation}
Furthermore, when $B_{nj}=W_{nj}$, for $j=1,2,\ldots,J$ in~\embed, then
\[
B_{n0}\circ L_0(t)
=
\sum_{j=1}^Jc_j^{1/2}W_{nj}\circ L_j(t)
\stackrel{d}{=}
W\left(\sum_{j=1}^Jc_jL_j(t)\right),
\]
where $W$ denotes Brownian motion.
Hence, using the monotonicity of $L_j$, for $j=0,1,\ldots,J$, from comparing the covariance functions,
it follows that we must have
\begin{equation}
\label{eq:defL0}
L_0(t)=\sum_{j=1}^J c_j L_j(t).
\end{equation}
This $L_0$ is increasing such that $\inf_{t\in[a,b]}L_0'(t)>0$ and
the second part of~\regul\ for $j=0$ follows immediately from the one for~$j=1,2,\ldots,n$.
Note that, in contrast to~$B_{nj}$, for $j=1,2,\ldots,J$, the process~$B_{n0}$ is not necessarily a
Brownian motion or a Brownian bridge.
However, we do have the following version of condition~\embed.
\begin{lemma}
\label{lem:embed0}
Assume \embed\ and suppose $f_0=f_1=\cdots=f_J$.
Let $B_{n0}$ and $L_0$ be defined by~\eqref{eq:defBn0} and~\eqref{eq:defL0}, and let
$M_{n0}=F_{n0}-F_0$, where $F_{n0}$ is defined in~\eqref{eq:def Fn0}.
Then, there exists $C>0$, such that for all $x\in(0,n]$:
\begin{equation}
\label{eq:embed0}
\prob\left\{
n^{1-1/q}
\sup_{t\in[a,b]}
\left|
M_{n0}(t)-n^{-1/2}B_{n0}\circ L_0(t)
\right|
>x
\right\}
\leq
Cx^{-q}.
\end{equation}
\end{lemma}
Now that we have established assumptions (A0)-(A4) for the pooled estimator,
we proceed by introducing a suitable variable, such as the one defined in~\eqref{eq: def zetaj}, for the case $j=0$.
However, this case is more complex and we have to distinguish between two of them.
First, for each fixed $t\in[a,b]$, define
\begin{equation}
\label{eq: def zetat0}
\begin{split}
\widetilde{\zeta}_{t0}(c)
&=
\argmax_{u\in\R}
\left\{
\widetilde{W}_{t0}(u+c)
-
u^2
\right\},\\
\widehat{\zeta}_{t0}(c)
&=
\argmax_{u\in\R}
\left\{
\widehat{W}_{t0}(u+c)
-
u^2
\right\},
\end{split}
\end{equation}
where
\begin{equation}
\label{eq:def Wt0}
\begin{split}
\widetilde{W}_{t0}(u)
&=
\sum_{j=1}^J\left(\frac{c_jL_j'(t)}{L_0'(t)}\right)^{1/2}W_j(u)\\
\widehat{W}_{t0}(u)
&=
\sum_{j=1}^J
c_j^{1/2}
W_{j}
\left(
n^{1/3}
\left\{
L_j\circ L_0^{-1}(L_0(t)+n^{-1/3}u)
-
L_j(t)
\right\}
\right),
\end{split}
\end{equation}
with $W_1,W_2,\ldots,W_J$ being the independent standard Brownian motions used to define~\eqref{eq: def zetaj}
and~$L_0$ defined in~\eqref{eq:defL0}.
Note that for $t\in[a,b]$ fixed, due to~\eqref{eq:defL0}, the processes~$\widetilde{W}_{t0}$
and~$\widehat{W}_{t0}$ are distributed as standard Brownian motion, which means
that $\widetilde{\zeta}_{t0}(c)$ and $\widehat{\zeta}_{t0}(c)$ have the same distribution as $\zeta_j(c)$.
We are now in the position to formulate our second main theorem.
\begin{theorem}
\label{th:main2}
Assume  \indep, \mono, \mart, \embed, \regul\
and let $S_{n2}$ be defined by~\eqref{eq:defS2}.
Let $\zeta_j$, $\widetilde{\zeta}_{t0}$ and $\widehat{\zeta}_{t0}$
be defined in~\eqref{eq: def zetaj} and~\eqref{eq: def zetat0}, respectively,
with independent standard Brownian motions $W_1,W_2,\ldots,W_J$.
If $f_0=f_1=\cdots=f_J$, then $n^{1/6}(n^{1/3}S_{n2}-m_{2})$
converges in distribution, as $n\to\infty$, to the Gaussian law with mean zero and variance~$\sigma_2^2$,
where
\[
\sigma^2_2
=
8
\sum_{i=1}^J
\sum_{j=1}^J
\int_a^b\int_0^{\infty}
\mathrm{cov}\left(|Y_{si}(t)-Y_{s0}(t)|, |Y_{sj}(0)-Y_{s0}(0)|\right)\,\md t\,\md s,
\]
with $Y_{sj}$ defined in~\eqref{eq:def Ysj} and
\[
Y_{s0}(t)
=
L_{0}'(s)^{1/3}
\widetilde{\zeta}_{t0}\left(\frac{t}{L_{0}'(s)^{1/3}}\right).
\]
Furthermore, $m_{2}$ may depend on $n$ and is defined by
\[
m_{2}
=
\sum_{j=1}^J
\int_{a}^{b}
|4f_{0}'(t)|^{1/3}
\E\left|
L_0'(t)^{1/3}\widehat{\zeta}_{t0}(0)
-
\frac{L_j'(t)^{1/3}}{c_j^{1/3}}\zeta_j(0)
\right|
\,\md t.
\]
If in addition, $L_j=a_j L$, for all $j=1,2,\ldots,J$,
for a given function $L:[a,b]\to\R$ and given real numbers~$a_j$,
then $\widehat{\zeta}_{t0}=\widetilde{\zeta}_{t0}$ and $m_{2}$ no longer depends on $n$.
\end{theorem}
The difference between the limiting bias $m_{2}$ and $\mathbb{E}(n^{1/3}S_{n2})$ will be shown to be of the order~$o(n^{-1/6})$.
Although $\widehat{\zeta}_{t0}(0)$ can be approximated further by $\widetilde{\zeta}_{t0}(0)$, this approximation is not sufficiently strong to
cancel the factor $n^{1/6}$.
This difficulty does not play a role for $\mathrm{var}(n^{1/3}S_{n2})$, for which we only need a consistent approximation.
For this reason the limiting bias $m_{2}$ may still depend on $n$, whereas the limiting variance $\sigma_2^2$ is independent of $n$.
Only in specific cases, such as $L_j=a_j L$, the limiting bias will also not depend on $n$.
Examples are the monotone density model, where under the null hypothesis~$L_j=F_j=F_{0}$, and the monotone regression model,
where $L_j(t)=(t-a)(b-a)\tau_j^2$.
Similar to Theorem~\ref{th:main1}, in these two cases the limiting variance is again independent of the underlying distribution.

The explicit expressions for the normalizing constants in Theorems~\ref{th:main1} and~\ref{th:main2}
are intractable for the purpose of building a statistical test because they depend on the $f_j$'s and the $L_j$'s in a complicated manner.
Therefore, in order to implement our statistical test,
we prefer to approximate the limit distribution of our test statistics using bootstrap methods,
as described in the following section.

\bigskip

Before doing so, we believe it is useful to give the main line of reasoning used to prove Theorems~\ref{th:main1} and~\ref{th:main2}
and explain the main differences with the type of argument used to prove similar results
in~\cite{groeneboom-hooghiemstra-lopuhaa1999},~\cite{kulikov-lopuhaa2005} and~\cite{durot2007}.
First note that it suffices to prove the results for the case $[a,b]=[0,1]$.
This is explained in more detail in the following remark.
\begin{remark}\label{remark}
Suppose that for $t\in[a,b]$ and $j=1,2,\ldots,J$, $f_j(t)$ satisfies conditions~\indep, \mono, \mart, \embed\ and~\regul\ with
corresponding $F_j$, $L_j$ and $F_{nj}$ on $[a,b]$.
Then this case can be transformed to the case~$[0,1]$ by
considering $(b-a)f_{j}(a+x(b-a))$ for $x\in[0,1]$.
It is straightforward to see that, for $j=1,2,\ldots,J$, these are functions on $[0,1]$ that
satisfy~\indep, \mono, \mart, \embed\ and~\regul\ with
corresponding functions $F_j(a+x(b-a))$, $L_j(a+x(b-a))$ and $F_{nj}(a+x(b-a))$ for $x\in[0,1]$.
Moreover, note that the transformed estimator $(b-a)\hat f_{nj}(a+x(b-a))$ is the left-hand slope of the least concave majorant of the process
$\{F_{nj}(a+u(b-a)), u\in[0,1]\}$ at the point $u=x$.
When Theorems~\ref{th:main1} and~\ref{th:main2} have been established for the case $[0,1]$,
then the results for the general case $[a,b]$ follow immediately and the expressions of the limiting constants
can by found by plugging in the transformed expressions.
\end{remark}
Hence, in the rest of the section we assume $[a,b]=[0,1]$.
When we define
\begin{eqnarray}
\label{eq:def Wn0}
W_{n0}(t)&=&\sum_{j=1}^J c_j^{1/2}W_{nj}\circ L_j\circ L_0^{-1}(t),\\
\label{eq:def xi0}
\xi_{n0}(t)&=&\sum_{j=1}^J c_j^{1/2}\xi_{nj}L_j\circ L_0^{-1}(t),
\end{eqnarray}
with $B_{n0}$ and $L_0$ defined by~\eqref{eq:defBn0} and~\eqref{eq:defL0}, respectively,
then similar to~\eqref{eq:Bnj and Wnj}, also $B_{n0}(t)=W_{n0}(t)-\xi_{n0}(t)$,
where the process $W_{n0}(t)$ is a standard Brownian motion and~$\xi_{n0}$ is independent of $B_{n0}$.
For every $j=0,1,\ldots,J$, define
\begin{equation}
\label{eq:def F_nj^S}
\begin{split}
F_{nj}^E(t)&=F_{nj}(t),\\
F_{nj}^B(t)&=n_j^{-1/2}B_{nj}(L_j(t))+F_j(t),\\
F_{nj}^W(t)&=n_j^{-1/2}W_{nj}(L_j(t))+F_j(t).
\end{split}
\end{equation}
and for $S=E,B,W$
\begin{equation}
\label{eq:def F_n0^S}
F_{n0}^S(t)=\sum_{j=1}^J c_jF_{nj}^S(t).
\end{equation}
Obviously, $F_{nj}^E$ is the estimator for~\eqref{eq:def Fj} or $F_{0}$ and, although
$F_{nj}^B$ and $F_{nj}^W$ are not estimators in the sense that they are built from observations,
we can define the corresponding slope processes
\begin{equation*}
\hf_{nj}^S(t)=\text{left slope of the least concave majorant (LCM) of $F_{nj}^S(u)$ at $u=t$,}
\end{equation*}
for $j=0,1,\ldots,J$ and $S=E,B,W$.
When investigating the asymptotic behavior of $\widehat{f}_{ni}^E$, one typically
exploits the fact that $F_{ni}^E$ can be approximated by $F_{ni}^W$, using~\embed\
and~\eqref{eq:def F_nj^S}.
However, if two processes are uniformly close, then the slopes of the concave majorants of both processes are not necessarily uniformly close.
For this reason, we introduce the (generalized) inverse of $\widehat{f}_{nj}^S$,
defined for $a\in\R$ by
$\widehat{U}_{nj}^S(a)=\sup\{t\in[0,1]: \widehat{f}_{nj}^S(t)\geq a\}$,
with the convention that the supremum of an empty set is zero.
It is fairly easy to see that
\begin{equation}
\label{eq:def Unj}
\widehat{U}_{nj}^S(a)=\argmax_{t\in[0,1]}\left\{F_{nj}^S(t)-at\right\},
\end{equation}
for $a\geq 0$, and that
\begin{equation}
\label{eq:switch}
\widehat{f}_{nj}^S(t)\geq a
\quad\text{if and only if}\quad
\widehat{U}_{nj}^S(a)\geq t.
\end{equation}
This means $\widehat{U}_{nj}^S$ is closely connected to $\widehat{f}_{nj}^S$, but its asymptotic behavior is more tractable
because if two processes are close, then also the locations of their maxima are close.
For this reason, the approach used in~\cite{groeneboom-hooghiemstra-lopuhaa1999}, \cite{kulikov-lopuhaa2005} and~\cite{durot2007},
which is originally due to~\cite{groeneboom1985},
is to switch from $L_p$-errors in terms of $\widehat{f}_{ni}^E$ to
$L_p$-errors in terms of $\widehat{U}_{ni}^E$.
The next lemma provides such an approximation suitable for our purposes.
For $S=E$, this result is similar to Corollary~2.1 in~\cite{groeneboom-hooghiemstra-lopuhaa1999},
Lemma~2.1 in~\cite{kulikov-lopuhaa2005} and equality~(21) in~\cite{durot2007}.
For later purposes, e.g, see~\eqref{eq:transition complete}, we also establish the approximation for the cases $S=B,W$.
\begin{lemma}
\label{lem:switch}
Assume \mono, \mart, \embed,~and suppose $f_0=f_1=\dots=f_J$.
Then for $i,j=0,1,\ldots,J$ and $S=E,B,W$,
\begin{equation}
\label{switch int}
n^{1/3}\int_0^1
|\widehat{f}_{ni}^S(t)-\widehat{f}_{nj}^S(t)|\,\md t
=
n^{1/3}\int_{f(1)}^{f(0)}
|\widehat{U}_{ni}^S(a)-\widehat{U}_{nj}^S(a)|\,\md a
+
o_p(n^{-1/6}).
\end{equation}
\end{lemma}
Proceeding in the spirit of~\cite{groeneboom-hooghiemstra-lopuhaa1999}, \cite{kulikov-lopuhaa2005} and~\cite{durot2007},
the next step would be to replace $\widehat{U}_{nj}^E$ on the right hand side of~\eqref{switch int}
by $\widehat{U}_{nj}^W$.
In statistical models where \embed\ holds with $B_{nj}$ being Brownian motion,
e.g., the regression model and the random censoring model,
this can be done by means of Lemma~5 in~\cite{durot2007}.
However, in statistical models, where $B_{nj}$ is Brownian bridge rather than Brownian motion,
e.g., the density model, this is no longer possible.
In such models, the approximation of $\widehat{U}_{nj}^E$ by $\widehat{U}_{nj}^B$ is relatively easy,
due to assumption~\embed.
This assumption will ensure that $F_{nj}^E-F_{nj}^B$ is of order smaller
than $n^{-5/6}$, which in turn guarantees that $\widehat{U}_{nj}^E-\widehat{U}_{nj}^B$ will be sufficiently small.
\begin{lemma}
\label{lem:E to B}
Assume \mono, \embed, \regul\ and suppose $f_0=f_1=\dots=f_J$.
Then, for each $j=0,1,2,\ldots,J$,
\[
n^{1/3}
\int_{f(1)}^{f(0)}
|\widehat{U}_{nj}^E(a)-\widehat{U}_{nj}^B(a)|\,\md a
=
o_p(n^{-1/6}).
\]
\end{lemma}
This result is similar to Corollary~3.1 in~\cite{groeneboom-hooghiemstra-lopuhaa1999}
for the density model, but is now extended to our general setup.
However, it is not possible to establish a similar result for~$\widehat{U}_{ni}^B$ and~$\widehat{U}_{ni}^W$,
because $F_{ni}^B-F_{ni}^W$ is of order $n^{-1/2}$, which is too large.
This difficulty is solved in~\cite{groeneboom-hooghiemstra-lopuhaa1999}, \cite{kulikov-lopuhaa2005}
and~\cite{durot2007}, by subtle use of relationship~\eqref{eq:Bnj and Wnj}.
These approaches apply to the $L_p$-error in terms of a single inverse $\widehat{U}_{ni}^B$,
but they do not extend to our current situation,
because the right hand side of~\eqref{switch int} involves the $L_1$-error
between two different inverses~$\widehat{U}_{ni}^B$ and $\widehat{U}_{nj}^B$.

For our current setup, we solve this problem by returning to the slopes themselves.
Whereas closeness of $F_{nj}^E$ and~$F_{nj}^B$ is not sufficient to obtain
suitable bounds on $\widehat{f}_{nj}^E-\widehat{f}_{nj}^B$, the situation is different
for~$F_{nj}^B$ and~$F_{nj}^W$.
The reason is that the difference between Brownian bridge~$B_{nj}$ and Brownian motion~$W_{nj}$ is only a straight line.
We can then obtain the following slope equivalent of Corollary~3.3 in~\cite{groeneboom-hooghiemstra-lopuhaa1999},
which is concerned with a similar approximation for the location processes defined in~\eqref{eq:def Unj}.
\begin{lemma}
\label{lem:transition}
Assume \mart, \embed.
Then, for $i,j=0,1,\ldots,J$,
\begin{equation*}
n^{1/3}\int_0^1 |\widehat{f}^B_{ni}(t)-\widehat{f}^B_{nj}(t)|\,\md t
=
n^{1/3}\int_0^1 |\widehat{f}^W_{ni}(t)-\widehat{f}^W_{nj}(t)|\,\md t
+o_p(n^{-1/6}).
\end{equation*}
\end{lemma}
After having established Lemmas~\ref{lem:switch}, \ref{lem:E to B} and~\ref{lem:transition},
in order to prove Theorems~\ref{th:main1} and~\ref{th:main2}, we will use the following line of reasoning:
\begin{equation}
\label{eq:transition complete}
\begin{split}
n^{1/3}\int_0^1
|\widehat{f}_{ni}^E(t)-\widehat{f}_{nj}^E(t)|\,\md t
&=
n^{1/3}\int_{f(1)}^{f(0)}
|\widehat{U}_{ni}^E(a)-\widehat{U}_{nj}^E(a)|\,\md a
+
o_p(n^{-1/6})\\
&=
n^{1/3}\int_{f(1)}^{f(0)}
|\widehat{U}_{ni}^B(a)-\widehat{U}_{nj}^B(a)|\,\md a
+
o_p(n^{-1/6})\\
&=
n^{1/3}\int_0^1
|\widehat{f}_{ni}^B(t)-\widehat{f}_{nj}^B(t)|\,\md t
+
o_p(n^{-1/6})\\
&=
n^{1/3}\int_0^1
|\widehat{f}_{ni}^W(t)-\widehat{f}_{nj}^W(t)|\,\md t
+
o_p(n^{-1/6})\\
&=
n^{1/3}\int_{f(1)}^{f(0)}
|\widehat{U}_{ni}^W(a)-\widehat{U}_{nj}^W(a)|\,\md a
+
o_p(n^{-1/6}).
\end{split}
\end{equation}
Of course, in models where $B_{nj}$ is Brownian motion, then $\widehat{U}_{nj}^B=\widehat{U}_{nj}^W$,
so that Lemma~\ref{lem:transition} becomes irrelevant, and we can obtain
\[
\begin{split}
n^{1/3}\int_0^1
|\widehat{f}_{ni}^E(t)-\widehat{f}_{nj}^E(t)|\,\md t
&=
n^{1/3}\int_{f(1)}^{f(0)}
|\widehat{U}_{ni}^E(a)-\widehat{U}_{nj}^E(a)|\,\md a
+
o_p(n^{-1/6})\\
&=
n^{1/3}\int_{f(1)}^{f(0)}
|\widehat{U}_{ni}^W(a)-\widehat{U}_{nj}^W(a)|\,\md a
+
o_p(n^{-1/6}),
\end{split}
\]
immediately, either using Lemma~5 in~\cite{durot2007} or as a special case of Lemma~\ref{lem:E to B}.
Once the test statistics~\eqref{eq:defS1} and~\eqref{eq:defS2} can be expressed in terms of integrals
\[
n^{1/3}\int_{f(1)}^{f(0)}|\widehat{U}_{ni}^W(a)-\widehat{U}_{nj}^W(a)|\,\md a,
\]
the proof of asymptotic normality follows the same line of reasoning as used in~\cite{groeneboom-hooghiemstra-lopuhaa1999},
\cite{kulikov-lopuhaa2005} and~\cite{durot2007}, using the independent increments property of Brownian motion.

\section{Calibration of the test}
\label{sec:calibration}
This section is devoted to the calibration of the test in
three different models that are covered by the general setup
for testing equality of monotone functions on $[a,b]$.
Because the limit distribution of the test statistic under the null hypothesis is intractable,
we use a bootstrap procedure to calibrate the test, see Subsection~\ref{subsec:bootstrap}.
The bootstrap procedure involves an estimator that is defined in Subsection~\ref{subsec:tildef}.
The three different models that we investigate are described in
Subsections~\ref{subsec:regression}, \ref{subsec:density} and~\ref{subsec:censoring} below.

\subsection{The bootstrap procedure}
\label{subsec:bootstrap}
It is known that the standard bootstrap typically does not work for Grenander-type estimators,
e.g., see~\cite{kos2008,banerjee-sen-woodroofe2010}.
These authors propose a smooth bootstrap based on generating from a kernel smoothed Grenander-type estimator.
\cite{caroll-delaigle-hall2011} and~\cite{du-parmeter-racine2013} discuss a smooth bootstrap based
on a monotonized kernel estimator, which consists of replacing the equal weights $1/n$ of the kernel estimator
by general weights $p_j$, $j=1,2,\ldots,n$ in such a way that the resulting estimator is monotone.

Here we also consider a smoothed bootstrap.
This will require the use of a smooth estimator $\widetilde f_n$ which, under the null hypothesis
$f_1=\dots=f_J=f_0$, satisfies bootstrap versions of assumptions \indep-\regul.
The following general property will be sufficient for our purposes.
\begin{enumerate}
\item[\monostar]
The estimator $\widetilde f_n$ is continuously differentiable on $[a,b]$.
Furthermore, there exists an event~$A_n$ and real numbers $\holder\in(3/4,1]$
and $\varepsilon_n>0$,
such that $\prob(A_n)\to 1$  and $n^{\gamma}\varepsilon_n\to\infty$ for any $\gamma>0$, as $n\to\infty$,
and such that the following three properties hold on $A_n$:
\begin{eqnarray}
\label{eq: approxtildef}
\sup_{t\in[a,b]}|\widetilde f_n(t)-f_0(t)|&=&o(n^{-1/3}),\\
\label{eq: approxtildef'}
\sup_{t\in[a,b]}|\widetilde f_n'(t)-f_0'(t)|&=&o(n^{-1/6}),
\end{eqnarray}
and for all $x,y\in[a,b]$,
\begin{equation}
\label{eq: tildef'}
|\widetilde f_n'(x)-\widetilde f_n'(y)|\leq |x-y|^{\holder}/\varepsilon_n.
\end{equation}
\end{enumerate}
Condition~\eqref{eq: approxtildef} comes naturally from minimax rates considerations for kernel density estimators,
in situations where the underlying density satisfies~\regul.
Condition~\eqref{eq: approxtildef'} ensures  the bootstrap version of assumption~\mono: if both~{\mono} and~\eqref{eq: approxtildef'} hold,
then  there exist positive numbers $C_0,C_1$ such that on $A_n,$ the function $\widetilde f_n$ is decreasing with
\begin{equation}
\label{eq: boundtildef'}
C_0<\inf_{t\in[a,b]}|\widetilde f_n'(t)|\leq \sup_{t\in[a,b]}|\widetilde f_n'(t)|<C_1.
\end{equation}
Condition~\eqref{eq: tildef'} ensures part of the bootstrap version of assumption~\regul.
It will have the same effect as assumption~\regul,
because $n^\gamma\varepsilon_n\to\infty$ for all $\gamma>0$ (typically one should think of~$\varepsilon_n=1/\log n$).
Bootstrap versions of \indep, \mart, \embed, and the second part of~\regul\ require the definitions of estimators
for $F_{j}$ and $L_j$.
This will be taken care off in Sections~\ref{subsec:regression}, \ref{subsec:density} and~\ref{subsec:censoring} for
the three different models that are covered by the general setup.

By means of the estimator $\widetilde f_n$, we aim to build bootstrap versions $S_{nk}^\star$
of test statistics~$S_{nk}$, for $k=1,2$,
in such a way that under the null hypothesis and conditionally on the original observations,
$n^{1/6}(S_{nk}^\star-m_k)$ converges in distribution to
the Gaussian law with mean zero and variance $\sigma_k^2$, in probability,
i.e.,
\begin{equation}\label{eq:cvboot}
\sup_{t\in\R}
\left|\prob^\star\left\{\frac{n^{1/6}(S_{nk}^\star-m_k)}{\sigma_k}\leq t\right\}-\Phi(t)\right|
\to0,
\qquad
\text{in probability, as }n\to\infty,
\end{equation}
where $m_k$ and $\sigma_k^2$ are the limit bias and variance given in Theorems~\ref{th:main1} and~\ref{th:main2},
$\Phi$ denotes the distribution function of the standard Gaussian law, and $\prob^\star$
is the conditional probability given the original observations.
In this case, for a fixed level $\alpha\in(0,1)$, one can compute (or merely approximate thanks to Monte-Carlo simulations)
the $\alpha$-upper percentile point $q_{nk}^\star(\alpha)$ of the conditional distribution of $S_{nk}^\star$
and consider the critical region
\begin{equation}
\label{eq:critical}
\left\{S_{nk}>q_{nk}^\star(\alpha)\right\}.
\end{equation}
If assumptions \indep --\regul~ are fulfilled, then Theorems~\ref{th:main1} and~\ref{th:main2}
together with (\ref{eq:cvboot}) ensures that the test with critical region (\ref{eq:critical})
has asymptotic level $\alpha$.

Below, we will provide an estimator $\widetilde f_n$ satisfying \monostar~ under the null hypothesis in the general framework of
Section~\ref{sec:main}.
Subsequently, we provide, in the three different models that are covered by this general framework,
a construction of $S_{nk}^\star$ that ensures that the test with critical region (\ref{eq:critical})
has asymptotic level $\alpha$.

\subsection{Estimators for the bootstrap procedure}
\label{subsec:tildef}
In this subsection, we discuss possible estimators to be used in the bootstrap procedure.
For simplicity, we assume here that under the null hypothesis,
the function $f_0=f_1=\dots=f_J$ is twice continuously differentiable.
We consider a sequence of positive real numbers $h_n$ and a kernel function $K:\R\to\R$
supported on $[-1,1]$, which is symmetric around zero and three times continuously differentiable on $\R$,
such that $\int K(t)\,\md t=1$.

Based on $h_n$ and~$K$, we
consider a kernel-type estimator $\widetilde f_n$, corrected at the boundaries in such a way that $\tilde f_{n}$ and
$\widetilde f_{n}'$ converge to $f_{0}$ and $f_{0}'$, respectively,
with a fast rate over the whole interval $[a,b]$ (whereas we recall that the non-corrected kernel estimator may show difficulties at the boundaries).
For every $t\in[a+h_{n},b-h_{n}]$ we define
\begin{equation}
\label{eq:tildefmain}
\widetilde f_n(t)
=
\frac{1}{h_n}\int_\R K\left(\frac{t-x}{h_n}\right)\,\md F_{n0}(x).
\end{equation}
At the boundaries $[a,a+h_{n})$ and $(b-h_{n},b]$, we discuss two possible bias corrections.

The first one is local
linear fitting (see e.g.~\cite{wand-jones1995}) that was used by~\cite{kos2008} in a similar context.
It is defined as follows: for every $t\in[a,a+h_n]\cup[b-h_n,b]$,
\begin{equation}
\label{eq:tildefKosorok}
\widetilde f_n(t)=\widetilde f_n(u_n)+\widetilde f_n'(u_n)(t-u_n),
\end{equation}
where $u_n=a+h_n$ for $t\in[a,a+h_n]$ and $u_n=b-h_n$ for $t\in[b-h_n,b]$.
Note that~\eqref{eq: tildef'} holds with $\holder=1$ provided that the second derivative of $\widetilde f_{n}$ is bounded from above by~$1/\varepsilon_{n}$.
It can be proved that with the boundary correction \eqref{eq:tildefKosorok},
the supremum norm of $\widetilde f_{n}''$ is of order $1+ h_{n}^{-5/2}n^{-1/2}\sqrt{\log (1/h_{n})}$
if $h_{n}$ is of order at least $n^{-2/3}$ and $f_{0}$ is twice continuously differentiable, so the optimal choice $h_{n}\sim n^{-1/5}$
is allowed thanks to the presence of $\varepsilon _{n}$.
In~\cite{kos2008}, the author requires the second derivative to be bounded independently of $n$,
which rules out the choice $h_{n}\sim n^{-1/5}$.
However, it can be checked that his result still holds under the assumption that the supremum norm
of~$\widetilde f_{n}''$ is bounded by some $1/\varepsilon_{n}$ satisfying our assumptions,
which means that the choice $h_{n}\sim n^{-1/5}$ is actually allowed in his case.

Another method to correct the bias is the use of boundary kernels
(see e.g.~\cite{gasser-muller1979}, \cite{gasser-muller-mammitzsch1985}).
One possibility is to construct linear combinations of $K(u)$ and $uK(u)$
with coefficients depending on the value near the boundary (see e.g.~\cite{zhang-karunamuni1998}).
For $t\in[a,a+h_n]\cup[b-h_n,b]$, define
\begin{equation}
\label{def:boundary kernel estimate}
\widetilde f_{n}(t)
=
\int
\frac1{h_n}\bK\left(\frac{t-x}{h_n}\right)
\,\md F_{n0}(x),
\end{equation}
with
\begin{equation*}
\bK(u)
=
\begin{cases}
\bcoefA\left(\frac{t-a}{h_n}\right)K(u)+\bcoefB\left(\frac{t-a}{h_n}\right)uK(u) & t\in [a,a+h_n],\\
\bcoefA\left(\frac{b-t}{h_n}\right)K(u)-\bcoefB\left(\frac{b-t}{h_n}\right)uK(u) & t\in [b-h_n,b],
\end{cases}
\end{equation*}
for $u\in\R$,
where for $s\in[-1,1]$, the coefficients~$\bcoefA(s)$ and $\bcoefB(s)$ are determined by
\begin{equation}
\label{eq:def coef boundary kernel}
\begin{split}
\bcoefA(s)\int_{-1}^{s} K(u)\,\md u&+\bcoefB(s)\int_{-1}^{s} uK(u)\,\md u=1,\\
\bcoefA(s)\int_{-1}^{s} uK(u)\,\md u&+\bcoefB(s)\int_{-1}^{s} u^2K(u)\,\md u=0.
\end{split}
\end{equation}
The following lemma guarantees that $\widetilde f_n$ with one of the above two boundary corrections satisfies condition~\monostar.
The proof is somewhat technical and has been put in the supplement.
\begin{lemma}
\label{lem: tildefn}
Let $\widetilde f_n(t)$ be defined by~\eqref{eq:tildefmain} for all $t\in[a+h_{n},b-h_{n}]$
and either by~\eqref{eq:tildefKosorok} or by~\eqref{def:boundary kernel estimate} on the boundaries $[a,a+h_{n})$ and $(b-h_{n},b_n]$.
Assume $h_{n}=R_{n}n^{-\hn}$, where $0<R_n+R_{n}^{-1}=O_{P}(1)$ and $\hn\in(1/6,1/5]$.
If $f_0=f_1=\dots=f_J$ is twice continuously differentiable on $[a,b]$ and \embed~ holds with
$\sup_{t\in[0,1]} L_j'(t)<\infty$, for each $j=1,2,\ldots,J$,
then $\widetilde f_n$ satisfies \monostar.
\end{lemma}
It may be more natural to use the least concave majorant~$\widehat{F}_{n0}$ of $F_{n0}$ instead of $F_{n0}$
in the definition of $\widetilde f_{n}$.
In this case the estimator is a smoothed Grenander type estimator, corrected at the boundaries.
Whether this estimator satisfies~\monostar\ will depend on how close $\widehat{F}_{n0}$ and $F_{n0}$ are.
For the density model, \cite{kiefer-wolfowitz1976} showed that the difference between $\widehat{F}_{n0}-F_{n0}$ is of the order $n^{-2/3}\log n$,
and a similar result has been obtained by~\cite{durottocquet2003} for the regression model;
see also~\cite{balabdaoui-wellner2007}.
This type of result for our general setting is proved in \cite{durot-lopuhaa2013}, whereby we prove the following lemma in the Appendix.
\begin{lemma}
\label{lem: brevefn}
Let $\widetilde f_n(t)$ be defined by~\eqref{eq:tildefmain} for all $t\in[a+h_{n},b-h_{n}]$
and either by~\eqref{eq:tildefKosorok} or by~\eqref{def:boundary kernel estimate} on the boundaries $[a,a+h_{n})$ and $(b-h_{n},b_n]$,  with $F_{n0}$ replaced by its least concave majorant~$\widehat{F}_{n0}$.
Assume $h_{n}=R_{n}n^{-\hn}$, where $0<R_n+R_{n}^{-1}=O_{P}(1)$ and $\hn\in(1/6,1/5]$.
If $f_0=f_1=\dots=f_J$ is twice continuously differentiable on $[a,b]$ and \embed~ holds with
$\sup_{t\in[0,1]} L_j'(t)<\infty$, for each $j=1,2,\ldots,J$,
then $\widetilde f_n$ satisfies \monostar.
\end{lemma}
It is tempting to consider the pooled Grenander type estimator $\widetilde f_n=\hat f_{n0}$ itself in the bootstrap procedure
since $\hat f_{n0}$ does not depend on any tuning parameter. Such a bootstrap procedure was used in~\cite{sen2010} in a similar context as our (a two-sample test
for monotone fractile regression functions), but no theoretical result was provided for their  procedure. Also in our context, we were not able to prove that bootstrapping from $\hat f_{n0}$ works, since $\hat f_{n0}$ does not satisfy~\monostar.  From the results in~\cite{banerjee-sen-woodroofe2010} and~\cite{kos2008} it appears that
bootstrapping from the Grenander does not work when the statistic of interest is $\widehat{f}_{n0}$ at a fixed point.
However, in our situation we are bootstrapping statistics that are integrals of the difference of two Grenander estimators, so it is not clear whether the results by~\cite{banerjee-sen-woodroofe2010} and~\cite{kos2008} apply. We investigated bootstrapping from $\hat f_{n0}$ in a simulation study reported in Section~\ref{sec:simulation}.

We end this section, by discussing possible estimators $F_{nj}$ in the three models that are covered by our setup.
Furthermore, for these models, we propose bootstrap versions of our test statistic for which we show that
the test with critical region (\ref{eq:critical}) has asymptotic level $\alpha$.
It suffices to specify bootstrap versions $F_{nj}^\star$, for $j=1,2,\ldots,J$.
Consequently, $F_{n0}^\star$ is defined similar to~\eqref{eq:def Fn0}, and
for $k=1,2$, bootstrap versions $S_{nk}^\star$ are defined similar to~\eqref{eq:defS1} and~\eqref{eq:defS2},
with~$\widehat{f}_{nj}^\star$ being the left-hand slope of the least concave majorant of~$F^\star_{nj}$.\\

\subsection{Monotone regression function}
\label{subsec:regression}
For each $j=1,2,\ldots,J$, we have observations $Y_{ij}$, for $i=1,2,\ldots,n_j$, satisfying
$Y_{ij}={f}_j(t_{ij})+\epsilon_{ij}$,
where $\E(\epsilon_{ij})=0$ and $t_{ij}=a+(b-a)i/n_j$, which means that the observation points are uniformly spread on $[a,b]$.
We assume that the $\epsilon_{ij}$'s are independent and that for each $j=1,2,\ldots,J$,
the variables $\epsilon_{ij}$, $i=1,2,\ldots,n_j$, have the same distribution with a finite variance $\tau_j^2>0$.
In this case,  the estimator for~$F_j$ is
\begin{equation}
\label{eq:def Fnj regression}
F_{nj}(t)
=
\frac{1}{n_j}\sum_{i=1}^{n_j} Y_{ij}\indicator\{t_{ij}\leq t\}.
\end{equation}
To define the bootstrap version of the test statistic, we first define $\widehat\epsilon_{ij}=Y_{ij}-\widetilde f_n(t_{ij})$,
where $\tilde f_{n}$ satisfies \monostar~under $H_{0}$. One can consider for instance one of the estimators from Lemmas~\ref{lem: tildefn} and~\ref{lem: brevefn} with $F_{n0}$ defined by~\eqref{eq:def Fn0} and~\eqref{eq:def Fnj regression}.
Then we define
$$
\widetilde\epsilon_{ij}=\widehat\epsilon_{ij}-\bar\epsilon_{j},
\quad
\text{where }
\bar\epsilon_{j}=\frac1{n_j}\sum_{i=1}^{n_j}\widehat\epsilon_{ij},
$$
for $j=1,2,\ldots,J$ and $i=1,2,\ldots,n_j$.
Then, conditionally on the original observations~$Y_{ij}$, we define independent random variables  $\epsilon_{ij}^\star$ as follows.
For $j=1,2,\ldots,J$ fixed, each random variable $\epsilon_{ij}^\star$ is uniformly distributed on
$\{\widetilde\epsilon_{mj},\ m=1,2,\dots,n_j\}$.
Finally we set
$$
Y_{ij}^\star=\widetilde f_n(t_{ij})+\epsilon_{ij}^\star,
$$
for $j=1,2,\ldots,J$ and $i=1,2,\ldots,n_j$.
For $j=0,1,\ldots,J$, we define bootstrap versions~$F_{nj}^\star$ in the same manner as $F_{nj}$
in~\eqref{eq:def Fnj regression} and~\eqref{eq:def Fn0}, just by replacing $Y_{ij}$ by $Y_{ij}^\star$.
The following theorem states that the bootstrap calibration (\ref{eq:critical}) is consistent under appropriate assumptions.
\begin{theorem}
\label{th:regression}
Suppose $\max_{1\leq j\leq J}\E|\epsilon_{ij}|^q<\infty$, for some $q>6$,
and that \mono\ and~\regul\ hold with $L_j(t)=(t-a)(b-a) \tau_j^2$.
Let $\widetilde f_n$ be an estimator that satisfies \monostar~under $H_{0}$ with $f_0=f_1=\dots=f_J$,
then for $k=1,2$, the test with critical region~(\ref{eq:critical}) has asymptotic level $\alpha$.
\end{theorem}

\subsection{Monotone density}
\label{subsec:density}
For each $j=1,2,\ldots,J$, we have independent observations ${X}_{ij}$, for $i=1,2,\ldots,n_j$,
with density~${f}_j:[a,b]\to\R$, where $a$ and $b$ are known real numbers.
The $J$ samples are assumed to be independent and we aim at testing that all observations are from the same density, that is ${f}_{1}=\dots={f}_{J}$.
The estimator for $F_j$ in this case is the empirical distribution function
\begin{equation}
\label{eq:def Fnj density}
F_{nj}(t)
=
\frac{1}{n_j}\sum_{i=1}^{n_j} \indicator\{X_{ij}\leq t\}.
\end{equation}
Let $\widetilde f_n$ be a genuine estimator for $f_0=f_1=\cdots=f_J$ in the sense that $\widetilde f_{n}$ is a density function.
One possibility is to use one of the estimators $\widetilde f_{n}$ from Lemmas~\ref{lem: tildefn} and~\ref{lem: brevefn} as a starting point.
Note that this function $\widetilde f_{n}$ need not integrate to one and may even be negative.
However, the function can be shifted upwards by
$-\min\{\inf \widetilde f_n,0\}$, so that it is positive on $[a,b]$, and then normalized so that it integrates to one.
It can be shown that if the original function satisfies condition~\monostar, so does the shifted and normalized version.
Moreover, the normalizing constant need not being computed
when generating from this function by means of rejection sampling.

To define the bootstrap versions of $F_{nj}$,
conditionally on the original observations~$X_{ij}$,
we define independent random variables  $X_{ij}^\star$, $j=1,2,\dots,J$, $i=1,2,\dots,n_j$,
with the same density~$\widetilde f_n$.
Then we define $F_{nj}^\star$ in the same manner as $F_{nj}$,
just by replacing~$X_{ij}$ by $X_{ij}^\star$ in~\eqref{eq:def Fnj density}.
\begin{theorem}\label{th:density}
Suppose that \mono~and~\regul~hold with $L_j'=f_j$ and $\inf_{t\in[a,b]}f_j(t)>0$, for each $j=1,2,\dots,J$.
Let $\widetilde f_n$ be a genuine
estimator that satisfies \monostar~under $H_{0}$ with $f_0=f_1=\dots=f_J$, then for $k=1,2$,
the test with critical region (\ref{eq:critical}) has asymptotic level $\alpha$.
\end{theorem}

\subsection{Random censorship with monotone hazard}
\label{subsec:censoring}
For each $j=1,2,\ldots,J$,
we have right-censored observations $(X_{ij},\Delta_{ij})$, for $i=1,2,\ldots,n_j$,
where $X_{ij}=\min(T_{ij},Y_{ij})$ and $\Delta_{ij}=\indicator\{T_{ij}\leq Y_{ij}\}$.
For each $j=1,2,\ldots,J$, the failure times~$T_{ij}$ are assumed to be nonnegative independent with density $g_j$
and to be independent of the i.i.d.~censoring times~$Y_{ij}$ that have distribution function $H_j$.
The $J$ samples are assumed to be independent.
The parameters of interest are the failure rates $f_j=g_j/(1-G_j)$ on $[0,b]$, where
$G_j=1-\exp(-F_j)$ is the distribution function corresponding to $g_j$.
Note that in this setting, we only consider the case $a=0$, since this is more natural.

The estimator for the cumulative hazard $F_j$ is defined via the Nelson-Aalen estimator~$N_{nj}$ as follows:
let $t_{1j}<\cdots<t_{mj}$ denote the ordered distinct uncensored failure times in the $j$th sample
and $n_{kj}$ the number of $i\in\{1,2,\dots,n_j\}$ with $X_{ij}\geq t_{kj}$,
then~$N_{nj}$ is constant on~$[t_{ij},t_{i+1,j})$ with
\[
N_{nj}(t_{ij})
=
\sum_{k\leq i}\frac{1}{n_{kj}},
\]
and $N_{nj}(t)=0$ for all $t<t_{1j}$ and $N_{nj}(t)=N_{nj}(t_{mj})$ for all $t\geq t_{mj}$.
The estimator $F_{nj}$ is the restriction of $N_{nj}$ to $[0,b]$.
Finally, as an estimator for the distribution function~$H_j$, we take
the Kaplan-Meier estimator $H_{nj}$ based on the $j$th sample.

Let $\widetilde f_n$ be a genuine estimator for $f_0=f_1=\cdots=f_J$ in the sense that
$\widetilde f_{n}$ is a non-negative failure rate.
One possibility is to use one of the estimators $\widetilde f_{n}$ from Lemmas~\ref{lem: tildefn} and~\ref{lem: brevefn} as a starting point, and to shift it upwards by $-\min\{-\inf \widetilde f_n,0\}$, so that it is positive on $[0,b]$.
It can be shown that if the original function satisfies condition~\monostar, so does the shifted version.
To define a bootstrap version of~$F_{nj}$,
conditionally on the original observations, we first define independent random variables
$T_{ij}^\star$ and $Y_{ij}^\star$, for $j=1,2,\dots,J$ and $i=1,2,\dots,n_j$,
where $T_{ij}^\star$ has failure rate $\widetilde f_n$ and $Y_{ij}^\star$ has distribution function $H_{nj}$.
Then we set
$X_{ij}^\star=\min(T_{ij}^\star,Y_{ij}^\star)$
and
$\Delta_{ij}^\star=\indicator\{T_{ij}^\star\leq Y_{ij}^\star\}$.
Finally, we define $F_{nj}^\star$ in the same manner as $F_{nj}$,
just replacing the $(X_{ij},\Delta_{ij})$'s by the $(X_{ij}^\star,\Delta_{ij}^\star)$'s in the definition.
\begin{theorem}
\label{th:censoring}
Suppose that \mono~and \regul~ hold with $L_j'=f_j/((1-G_j)(1-H_j))$ and that
for each $j=1,2,\dots,J$, $\inf_{t\in[0,b]}f_j(t)>0$,
$G_j(b)<1$, $\lim_{t\uparrow b}H_j(t)<1$,
and~$H_j$ has a bounded continuous first derivative on $[0,b]$.
Let $\widetilde f_n$ be a non-negative estimator that satisfies~\monostar~under $H_{0}$ with $f_0=f_1=\dots=f_J$,
then for $k=1$, the test with critical region~(\ref{eq:critical}) has asymptotic level $\alpha$.
\end{theorem}
If, in addition, we assume that the censoring variables all have the same distribution function~$H$,
then, instead of generating $Y_{ij}^\star$ from distribution function $H_{nj}$ as above, one should merely generate,
the bootstrap censoring times $Y_{ij}^\star$ as an $n$-sample of
independent random variables with common distribution function $H_{n}$,
the Kaplan-Meier estimator of $H$ based on all $n$ observations.
With this construction of the bootstrap censoring variables we obtain a similar result.
\begin{theorem}\label{th:censoring2}
Under the assumptions of Theorem~\ref{th:censoring}
with $H=H_1=\cdots=H_J$, for $k=1,2$,
the test with critical region~(\ref{eq:critical}) has asymptotic level $\alpha$.
\end{theorem}

\section{Simulation Study}
\label{sec:simulation}
To investigate the performance of bootstrapping the test statistics we have performed a simulation study.
To alleviate notation, in this section we sometimes omit subscript $n$, so a bandwidth is denoted by $h$ rather than $h_{n}$.
Moreover, we define $K_h(x)=h^{-1}K(x/h)$.

\subsection{Setup}
\label{subsec:setup}
We consider a 3-sample test in the monotone density model.
The three densities $f_1,f_2$ and $f_3$ are chosen from
the family of exponential densities truncated to the interval~$[0,3]$:
\begin{equation}
\label{def:truncated exponential}
f(x,\lambda)
=
\begin{cases}
\lambda \me^{-\lambda x}(1-\me^{-3\lambda})^{-1} & ,\lambda>0;\\
1/3 & ,\lambda=0,
\end{cases}
\end{equation}
for $x\in[0,3]$ and $f(x,\lambda)=0$ otherwise.
Under the null hypothesis $f_1=f_2=f_3=f_0$ the bootstrap samples are generated from a pooled estimate
for $f_0$ based on the pooled sample of size $n=n_1+n_2+n_3$.
We have several options to construct the smooth estimator~$\widetilde{f}_{n,h}$.
One can either smooth the empirical distribution function or smooth the Grenander estimator.
Furthermore we can correct the estimator at the boundaries either by~\eqref{eq:tildefKosorok} or by~\eqref{def:boundary kernel estimate}.
According to Lemmas~\ref{lem: tildefn} and~\ref{lem: brevefn} the bootstrap works for each of these combinations.
For the different possibilities, we first investigated their performance
when determining the bandwidth of the kernel estimate in a data-adaptive way.
\begin{figure}
\centering
\subfigure[$\mathrm{LSCV}(h)$ for the smoothed Grenander (solid) with boundary correction~\eqref{eq:tildefKosorok}.\label{fig:LSCV Grenander}]{\includegraphics[width=0.45\textwidth]{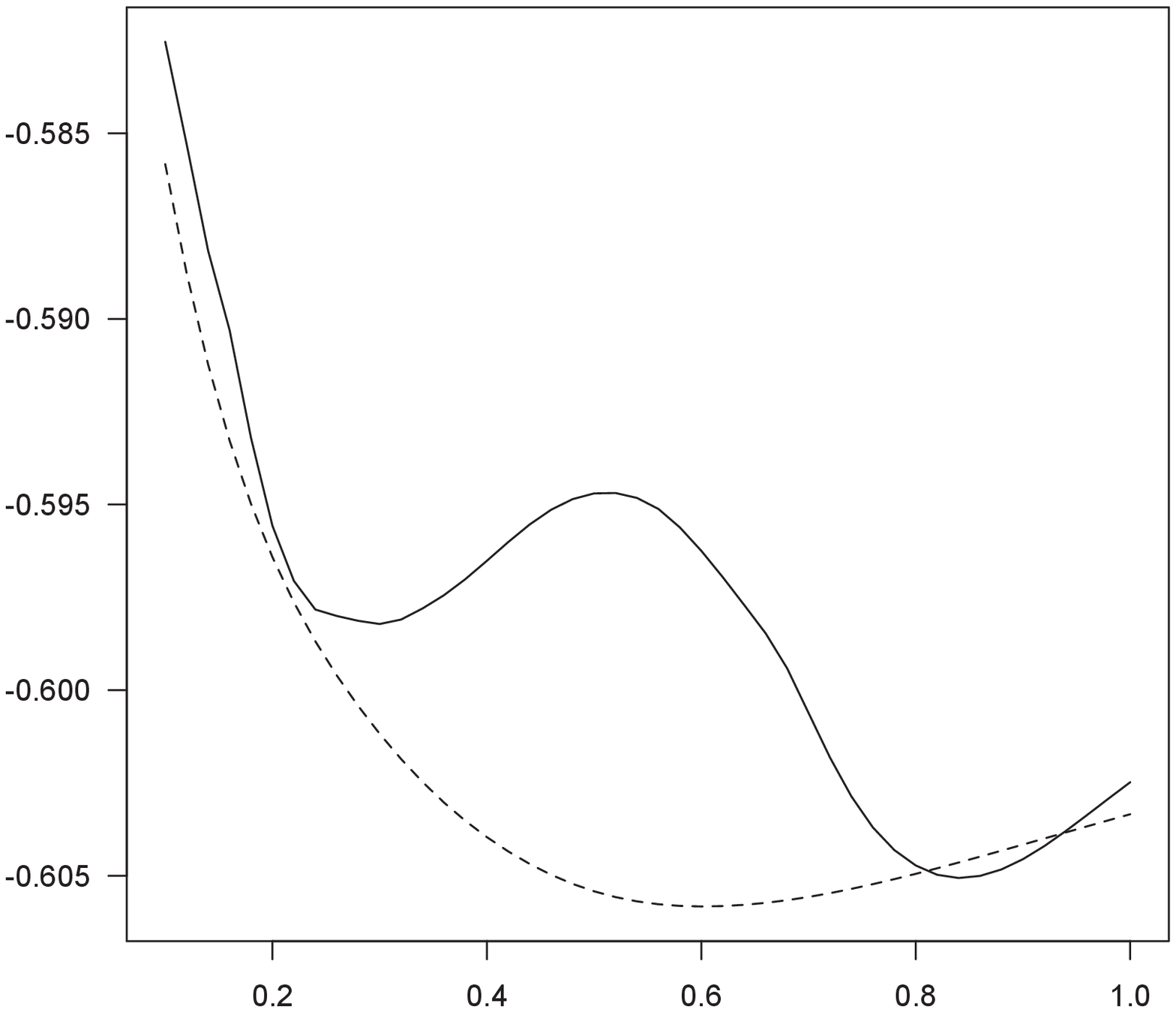}}
\hspace{1cm}
\subfigure[$\mathrm{LSCV}(h)$ for the ordinary kernel estimate (solid) with boundary correction~\eqref{eq:tildefKosorok}.\label{fig:LSCV Kernel}]{\includegraphics[width=0.45\textwidth]{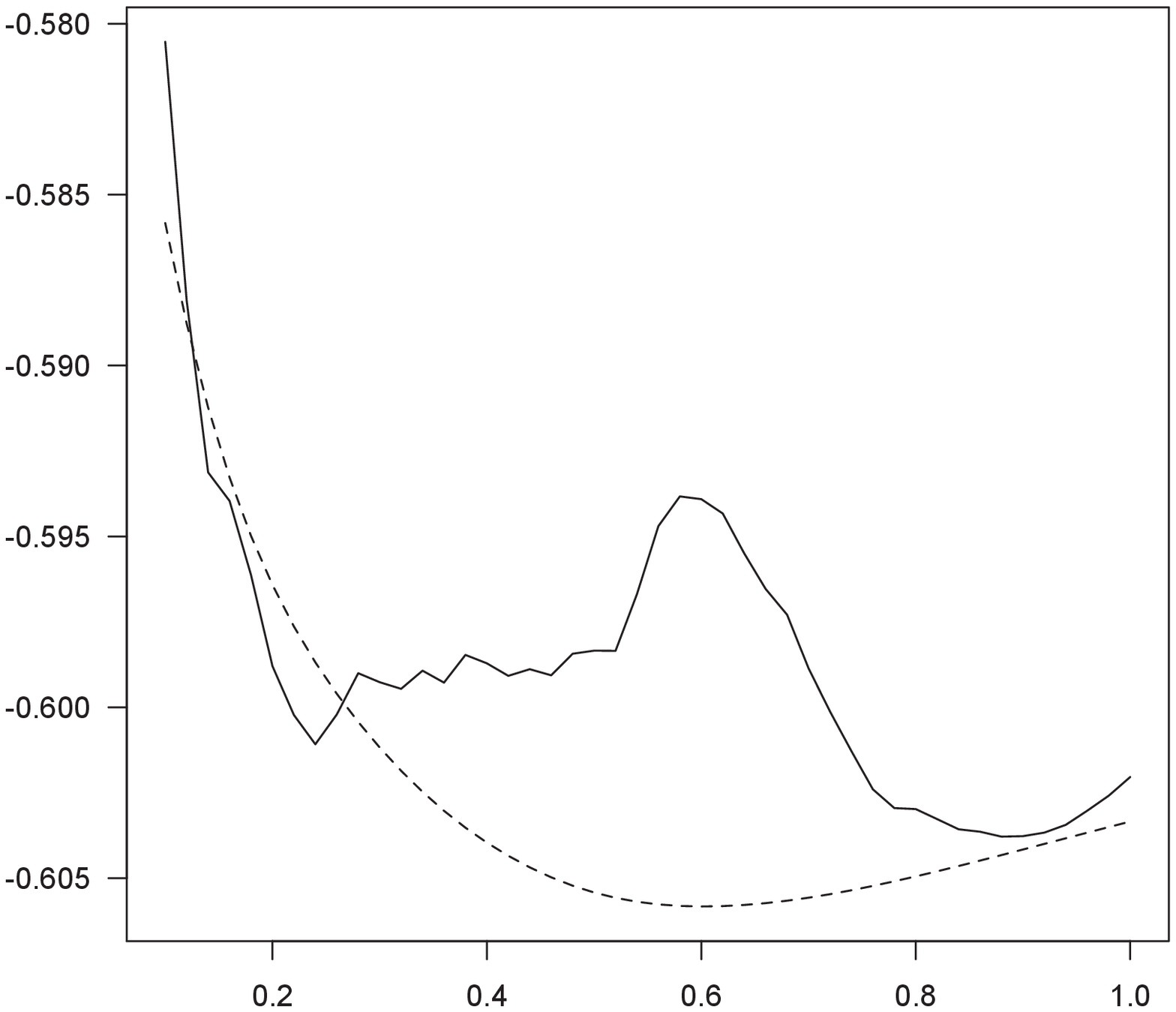}}
\caption{Cross-validation functions for boundary correction~\eqref{eq:tildefKosorok} (solid) and for the smoothed Grenander with boundary correction~\eqref{def:boundary kernel estimate} (dashed).}
\label{fig:LSCV}
\end{figure}
\paragraph*{Choice of bandwidth.}
We choose to use a so-called ``first generation method" instead of a ``second generation" adaptive plug-in method,
where we would have to assume the existence of third derivatives.
These methods also would be more complicated if we wish to take the Grenander estimator as the starting point of our smoothing method.
The first generation method of our preference is least squares cross-validation, adapted to the present situation,
where we possibly want to smooth the Grenander estimator instead of the empirical distribution function, as in ordinary density estimation.
For illustrative purposes, we consider the family of truncated exponentials in~\eqref{def:truncated exponential}.
In our experiments, the least squares cross validation function, as a function of the bandwidth $h$, is given by
\begin{equation}
\label{crossval_function}
\mathrm{LSCV}(h)=\int \widetilde{f}_{n,h}(t)^2\,dt-\frac{2n}{n-1}\int \widetilde{f}_{n,h}(x)\,\md F_{n0}(x)+\frac{2K(0)}{(n-1)h},
\end{equation}
where $\widetilde{f}_{n,h}$ is the (smooth) estimate of the density, based on the pooled samples, with bandwidth $h$
and $F_{n0}$ is the empirical distribution function of the pooled samples.
Note that if $\widetilde{f}_{n,h}$ is the ordinary kernel estimator
determined with the empirical distribution function $F_{n0}$,
then $\mathrm{LSCV}(h)+\int f^2(t)\,\md t$ is an unbiased estimator of the mean integrated squared error.

It turns out that least squares cross-validation does not work very well for the boundary correction method~\eqref{eq:tildefKosorok}.
One typically gets a very non-convex function, as illustrated in Figure~\ref{fig:LSCV}.
Figure~\ref{fig:LSCV Grenander} displays the cross-validation curve~\eqref{crossval_function}, for $h\in[0,1]$, for the smoothed Grenander
and Figure~\ref{fig:LSCV Kernel} displays the same curve for the ordinary kernel estimator,
both with boundary correction~\eqref{eq:tildefKosorok} (solid curves),
for a pooled sample of size $n=300$ from a truncated exponential with parameter $\lambda=1$.
It seems clear that the kernel estimate based on the Grenander gives a somewhat smoother cross-validation function
than the kernel estimate based on the empirical distribution function, but both curves are very non-convex.
For comparison, the cross-validation function for the smoothed Grenander with boundary correction~\eqref{def:boundary kernel estimate} has been added (dashed curves).

\begin{figure}
\centering
\subfigure[$\mathrm{LSCV}(h)$ for the smoothed Grenander (solid) and ordinary kernel estimator (dashed).\label{fig:estimates LSCV}]{\includegraphics[width=0.45\textwidth,clip=]{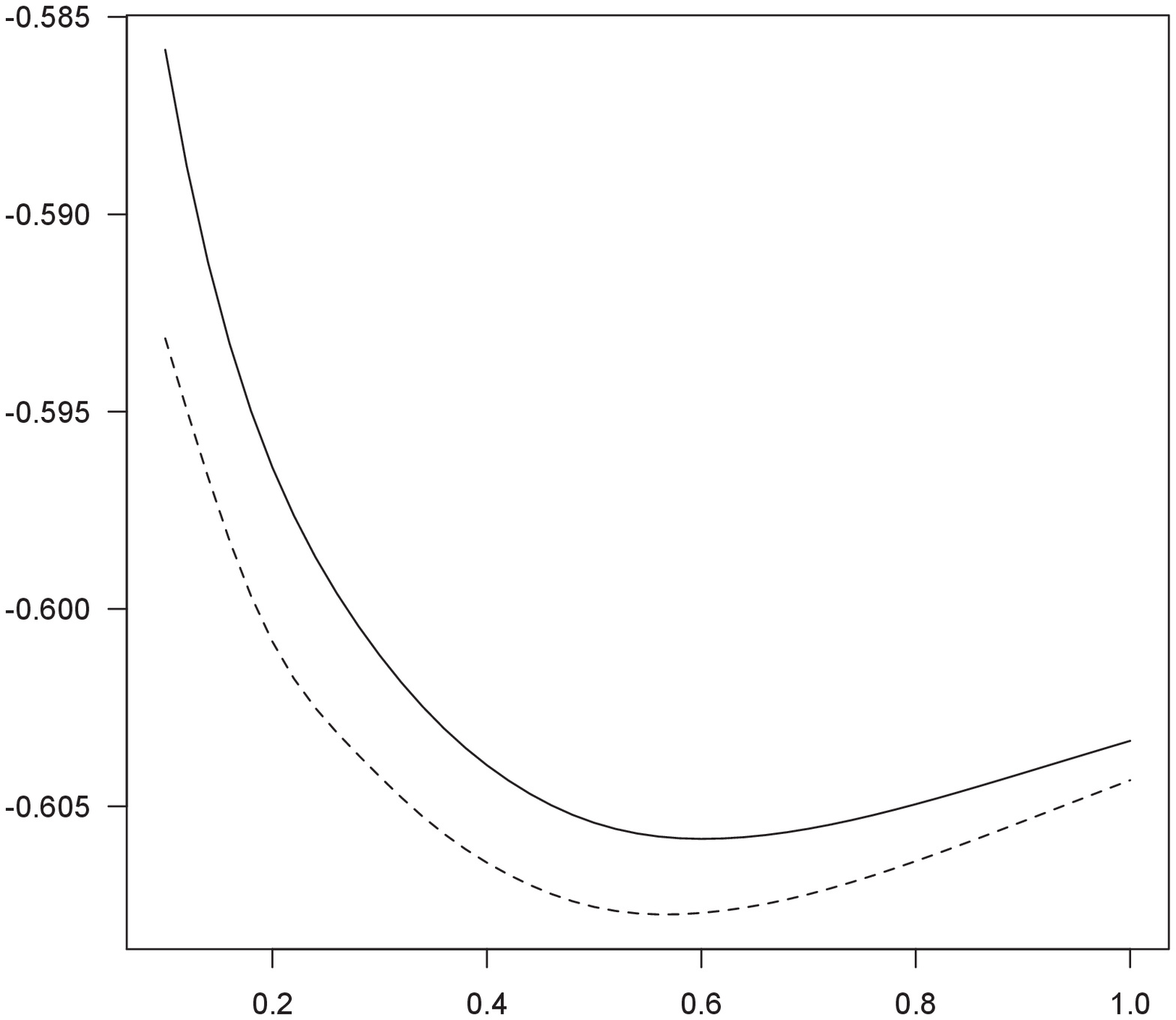}}
\hspace{1cm}
\subfigure[Smoothed Grenander (solid), ordinary kernel estimate (dotted), and the true density (dashed).\label{fig:estimates density}]{\includegraphics[width=0.45\textwidth,clip=]{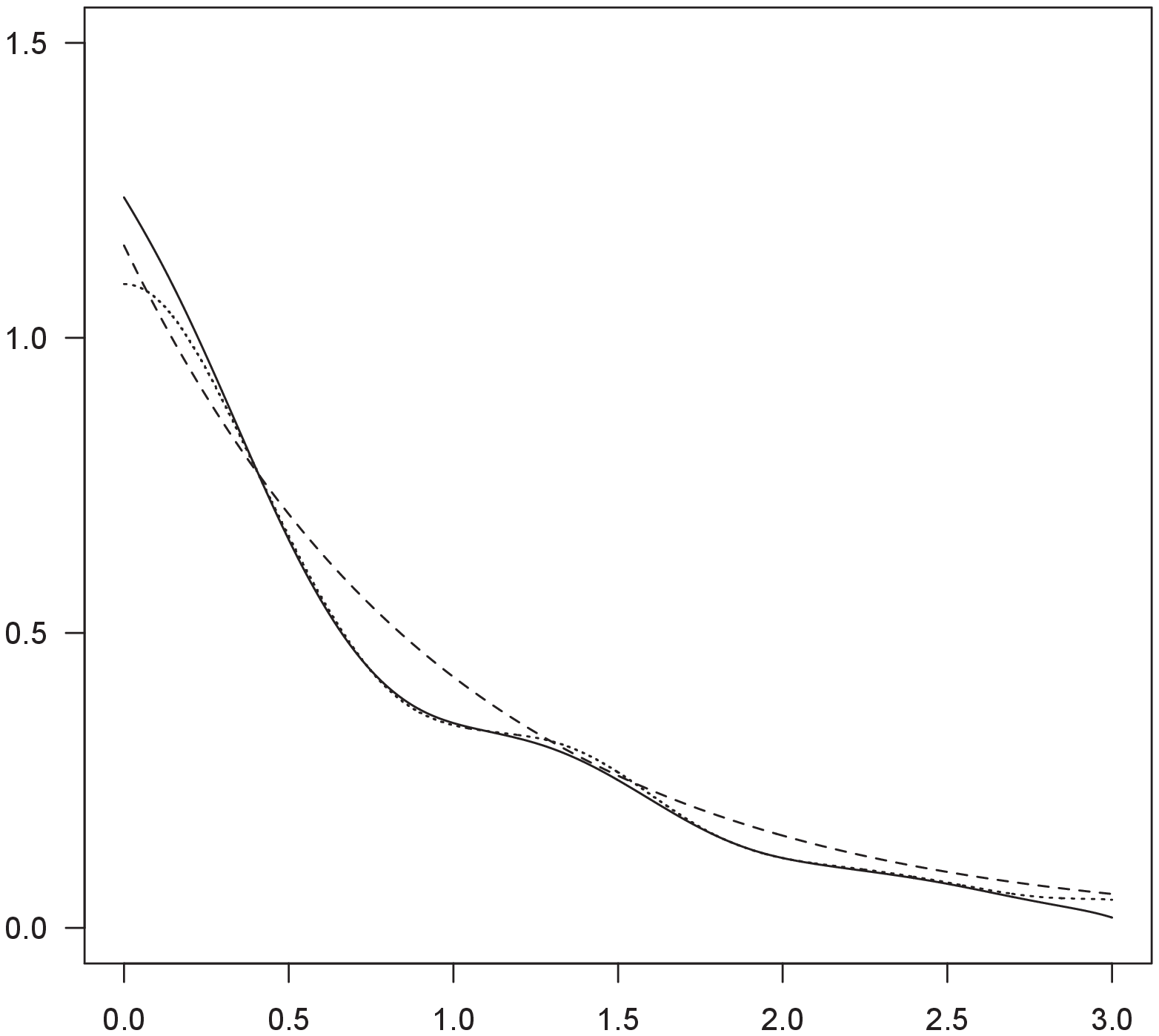}}
\caption{Cross-validation functions for the smoothed Grenander and the ordinary kernel estimator with boundary correction~\eqref{def:boundary kernel estimate} and the resulting density estimates.}
\label{fig:estimates}
\end{figure}

Boundary correction method~\eqref{def:boundary kernel estimate} generally leads to a convex cross-validation curve with a clear minimum,
as shown in Figure~\ref{fig:estimates LSCV}.
The cross-validation curve for the smoothed Grenander (solid) with the boundary correction~\eqref{def:boundary kernel estimate}
attains its minimum for a value of $h$ close to $0.6$, and the resulting kernel estimate (solid) is shown in Figure~\ref{fig:estimates density}.
The cross-validation curve in Figure~\ref{fig:estimates LSCV} of the ordinary kernel estimate (dashed) with the boundary correction~\eqref{def:boundary kernel estimate}
lies completely below the cross-validation curve of the kernel estimate based on the Grenander
and Figure~\ref{fig:estimates density} indeed shows that this kernel estimate is closer to the real density (dashed)
than the smoothed Grenander for $h=0.6$.
On the other hand, this kernel estimate will not necessarily be decreasing and we actually prefer a decreasing density like the smoothed Grenander,
which belongs to the allowed class of densities, for generating the bootstrap samples.
After a thorough investigation of the different possibilities, the overall performance of the smoothed Grenander with
boundary correction~\eqref{def:boundary kernel estimate} seems to be the best.
Therefore, for $t\in[h,3-h]$ our smooth estimate is defined as
\begin{equation}
\label{def:kernel estimate}
\widetilde f_{n,h}(t)
=
\int K_h(t-x)\,\md \widehat{F}_{n0}(x),
\end{equation}
where $K$ is a symmetric kernel with support $[-1,1]$ and $\widehat{F}_{n0}$ is the least concave majorant
of the empirical distribution function $F_{n0}$.
We correct the kernel density estimate at the boundaries of $[0,3]$ by means of~\eqref{def:boundary kernel estimate}
with $\widehat{F}_{n0}$ instead of $F_{n0}$.
\paragraph*{Simulating the level and power under alternatives.}
To investigate the finite sample power at a given combination $(\lambda_1,\lambda_2,\lambda_3)$,
we generate three samples of sizes $n_j$ from $f_j=f(\cdot,\lambda_j)$, for $j=1,2,3$, and compute the value of the test statistics $S_{n1}$
and~$S_{n2}$, as defined in~\eqref{eq:defS1} and~\eqref{eq:defS2}.
We then generate 1000 times three bootstrap samples of sizes~$n_1$, $n_2$ and~$n_3$ from the pooled estimate $\widetilde{f}_{n,h}$,
compute the values $S_{n1}^\star$ and~$S_{n2}^\star$ of both test statistics and determine their
5th upper-percentiles $q_{nk}^\star(0.05)$, for $k=1,2$.
This whole procedure is repeated $B$ times and
we count the number of times the values of the test statistics~$S_{n1}$ and~$S_{n2}$ exceed the corresponding~5th upper-percentiles~$q_{n1}^\star(0.05)$
and~$q_{n2}^\star(0.05)$, respectively.
By dividing this number by $B$, this provides an approximation of the finite sample power of both test statistics
at underlying truncated exponentials with parameters~$\lambda_1$, $\lambda_2$ and $\lambda_3$.

In view of the comments made right after~Lemma~\ref{lem: brevefn},
we compare the behavior of the smooth bootstrap procedure with bootstrapping from the pooled Grenander estimate itself.
To this end, we also run the same procedure as described above, but then generate the bootstrap samples from $\widehat{f}_{n0}$
instead of the smooth estimate~$\widetilde f_{n,h}$.

To investigate the performance under the null hypothesis, we take $\lambda_1=\lambda_2=\lambda_3$
equal to the values $0.1, 0.5, 1,2,\ldots,6$ and equal sample sizes $n_j=100$ and $n_j=250$, for $j=1,2,3$.
The simulated levels are determined by means of $B=10\,000$ repetitions.
The simulations to investigate the finite sample power at alternatives are done with sample sizes $n_1=n_2=n_3=100$
and alternatives for which $\lambda_1=\lambda_2=1$ and $\lambda_3$ varies between~0 and 3.5 by steps of 0.1.
To save computer time, we determined the simulated power at each $\lambda_3$ by means of $B=1000$ repetitions.
\paragraph*{Benchmark with true power.}
Finally, in order to calibrate the finite sample power obtained from bootstrapping, we also approximate the true finite sample power for
a given choice $(\lambda_1,\lambda_2,\lambda_3)$.
To this end, we generate $10\,000$ samples of size $n=n_1+n_2+n_3$ from the mixture density
\[
f_0(x)
=
c_1
\frac{\lambda_1 \me^{-\lambda_1 x}}{1-\me^{-\lambda_1}}
+
c_2
\frac{\lambda_2 \me^{-\lambda_2 x}}{1-\me^{-\lambda_2}}
+
c_3\frac{\lambda_3 \me^{-\lambda_3 x}}{1-\me^{-\lambda_3}},
\]
where $c_j=n_j/n$, for $j=1,2,3$.
We consider this as  the least favorable density among all densities under the null hypothesis,
in case of three truncated exponentials with parameters $\lambda_1$, $\lambda_2$ and $\lambda_3$.
For each of the samples we compute the value of the test statistics $S_{n1}$ and $S_{n2}$,
and use this to determine the 5th upper-percentiles $q_{nk}(0.05)$, $k=1,2$, for both test statistics.
Next, we generate another 10\,000 times three samples of sizes $n_j$ from $f_j=f(\cdot,\lambda_j)$,
compute both test statistics and count the number of times it exceeds the corresponding 5th percentile $q_{nk}(0.05)$.
Dividing these numbers by 10\,000 provides an approximation of the true finite sample power
for a given choice $(\lambda_1,\lambda_2,\lambda_3)$. Note that such a calibration is not implementable in practice since it requires knowledge of $f_{1},f_{2}$ and $f_{3}$, but it may serve as a benchmark for the power obtained from bootstrapping, in the simulations.
\subsection{Implementation}
\label{subsec:implementation}
We believe it useful to spend some words on how bootstrapping from the smoothed Grenander has been implemented.
First consider the estimate defined in~\eqref{def:kernel estimate} for $t\in[h,3-h]$.
One possibility to implement this estimate would be to use numerical integration of $\hat f_{n0}$.
However, one can also avoid this by using summation by parts.
Let $p_1,\dots,p_m$ be the jump sizes of the Grenander estimator at the points of jump $\tau_1<\dots<\tau_m\in(0,3)$,
where~$\tau_m$ is the largest order statistic. Note that $\hat f_{n0}$ is left-continuous and that $\hat f_{n0}$ always has a jump down to zero at the last order statistic.
We now define
\begin{equation}
\label{int_kernel}
\mathbb{K}_h(x)=\int_{x/h}^{\infty}K(u)\,du,\qquad x\in\R.
\end{equation}
In our simulations we
took $K(u)=(35/32)\left(1-u^2\right)^31_{[-1,1]}(u)$.
Then, when defining $\tau_0=0$, for $t\in[h,3-h]$, we can write,
\[
\begin{split}
\widetilde f_{n,h}(t)
&=\sum_{i=1}^m\left\{\sum_{j=i}^m p_j\right\}\int_{\tau_{i-1}}^{\tau_{i}}K_h(t-x)\,\md x
=\sum_{j=1}^m p_j\int_0^{\tau_j}K_h(t-x)\,\md x\\
&=\sum_{j=1}^m p_j\int_{(t-\tau_j)/h}^{t/h}K(u)\,\md u =\sum_{j=1}^m p_j\int_{(t-\tau_j)/h}^{1\wedge(t/h)}K(u)\,\md u
=\sum_{j=1}^m p_j\mathbb{K}_h(t-\tau_j),
\end{split}
\]
so that for $t\in[h,3-h]$,
the estimate $\widetilde f_{n,h}(t)$ can now be computed  as a finite sum over the jumps $p_i$ of the Grenander estimator $\hat f_{n0}$.
We then still have to define $\widetilde f_{n,h}(t)$ for $t\in[0,h)\cup(3-h,3]$.
To this end, for $j=0,1,2$, let
\[
\IK_h^{(j)}(x)=\int_{-\infty}^{x/h} u^jK(u)\,\md u.
\]
Note that $\IK_h^{(0)}(t)=1-\IK_h(t)$, where $\IK_h$ is defined in~\eqref{int_kernel}.
As before, we get for $t<h$,
\[
\begin{split}
\widetilde{f}_{n,h}(t)
&=
\int
\left\{\bcoefA\left(\frac{t}{h}\right)K_h(t-x)+\bcoefB\left(\frac{t}{h}\right)\frac{t-x}{h}K_h(t-x)\right\}\widehat f_{n0}(x)\,dx\\
&=\bcoefA\left(\frac{t}{h}\right)\sum_{j=1}^m p_j\int_{(t-\tau_j)/h}^{t/h}K(u)\,du
+\bcoefB\left(\frac{t}{h}\right)\sum_{j=1}^m p_j\int_{(t-\tau_j)/h}^{t/h}uK(u)\,du\\
&=\bcoefA\left(\frac{t}{h}\right)\sum_{j=1}^m p_j\left\{\IK_h^{(0)}(t)-\IK_h^{(0)}(t-\tau_j)\right\}
+\bcoefB\left(\frac{t}{h}\right)\sum_{j=1}^m p_j\left\{\IK_h^{(1)}(t)-\IK_h^{(1)}(t-\tau_j)\right\},
\end{split}
\]
where $\phi$ and $\psi$ are defined by~\eqref{eq:def coef boundary kernel},
and similarly for $t>3-h$,
\[
\widetilde{f}_{n,h}(t)
=
\bcoefA\left(\frac{3-t}h\right)\sum_{j=1}^{m} p_j\left\{1-\IK_h^{(0)}(t-\tau_j)\right\}+\bcoefB\left(\frac{3-t}h\right)\sum_{j=1}^{m} p_j\IK_h^{(1)}(t-\tau_j).
\]
This means that also near the boundaries of $[0,3]$, the estimator $\widetilde f_{n,h}(t)$ can be computed in terms of finite sums
over the jumps of the Grenander estimator $\hat f_{n0}$.
\subsection{Results}
\label{subsec:results}
We first investigate the level of the tests under the null hypothesis of all $\lambda$'s equal
to some~$\lambda_0$, where we vary $\lambda_0$ over $0.1,0.5,1,2,\ldots,6$.
We set the significance level $\alpha=0.05$ and perform the bootstrap experiments
with $n_1=n_2=n_3=100$ and $n_1=n_2=n_3=250$ .
The results are listed in Table~\ref{tab:levels}.
\begin{table}
\centering
\begin{tabular}{cccccccccc}
\hline
\hline
\\[-10pt]
& \multicolumn{4}{c}{$n_1=n_2=n_3=100$} & & \multicolumn{4}{c}{$n_1=n_2=n_3=250$} \\
& \multicolumn{2}{c}{Genander} & \multicolumn{2}{c}{Smooth} & & \multicolumn{2}{c}{Grenander} & \multicolumn{2}{c}{Smooth} \\
\\[-10pt]
\hline
\\[-10pt]
$\lambda_0$ & $S_{n1}$ & $S_{n2}$ & $S_{n1}$ & $S_{n2}$ & & $S_{n1}$ & $S_{n2}$ & $S_{n1}$ & $S_{n2}$\\
\\[-10pt]
\hline
\\[-10pt]
0.1 &.0122	 & 	.0215	 & 	.0109	 & 	.0195 &	 & 	.0173	 & 	.0256	 & 	.0166	 & 	.0241\\
0.5 &.0337	 & 	.0352	 & 	.0381	 & 	.0424 &	 & 	.0393	 & 	.0402	 & 	.0407	 & 	.0500\\
1 & .0417	 & 	.0405	 & 	.0474	 & 	.0477 &  & 	.0422	 & 	.0424	 & 	.0479	 & 	.0499\\	
2 & .0462	 & 	.0436	 & 	.0560	 & 	.0592 &	 & 	.0466	 & 	.0449	 & 	.0546	 & 	.0504\\	
3 & .0473	 & 	.0443	 & 	.0582	 & 	.0591 &	 & 	.0496	 & 	.0474	 & 	.0554	 & 	.0545\\	
4 & .0474	 & 	.0444	 & 	.0545	 & 	.0560 &	 & 	.0498	 & 	.0478	 & 	.0559	 & 	.0563\\	
5 & .0455	 & 	.0486	 & 	.0593	 & 	.0577 &	 & 	.0530	 & 	.0541	 & 	.0487	 & 	.0508\\	
6 & .0441	 & 	.0459	 & 	.0600	 & 	.0597 &	 & 	.0553	 & 	.0496	 & 	.0541	 & 	.0526\\	
\\[-10pt]
\hline
\hline
\end{tabular}
\caption{Simulated levels of $S_{n1}$ and $S_{n2}$ under the null hypothesis.\label{tab:levels} }
\end{table}
It can be seen that close to $\lambda_0=0$, which corresponds to the uniform distribution,
the attained level is much too small.
For large $\lambda_0$ the attained levels tend to be somewhat too large.
Note that the simulated levels obtained from bootstrapping from the Grenander itself are comparable
to the ones obtained from the smooth bootstrap.

Next, we investigate the power under alternatives of the form $f_1=f_2=f(\cdot,1)$ and $f_3=f(\cdot,\lambda)$ with $n_1=n_2=n_3=100$.
A picture of the power estimates of the smoothed Grenander, using cross-validation for the bandwidth choice,
is shown in Figure~\ref{fig:powers[0,3]} together with the estimates obtained by bootstrapping from the Grenander estimator.
Figure~\ref{fig:Grenander2[0,3]} displays the powers simulated by generating bootstrap samples from the ordinary Grenander estimator (solid curves)
and the direct estimates of the true power (dashed curves).
The top solid and dashed curves correspond to test statistic $S_{n2}$.
This test statistic seems to be uniformly more powerful than test statistic $S_{n1}$, which corresponds to the bottom solid and dashed curves.
Figure~\ref{fig:powers_smooth[0,3]} displays the powers simulated by generating bootstrap samples from the smoothed Grenander estimator (solid curves)
and the same direct estimates of the true power (dashed curves).
Again the top solid and dashed curves correspond to test statistic $S_{n2}$.

The simulated powers in Figure~\ref{fig:Grenander2[0,3]}, based on bootstrapping from the ordinary Grenander, tend to be conservative.
The simulated powers in Figure~\ref{fig:powers_smooth[0,3]}, based on bootstrapping from the smoothed Grenander tend to be slightly anti-conservative.
Note that, similar to the simulated levels in Table~\ref{tab:levels},
there is hardly any difference between the results when using the smooth bootstrap or when bootstrapping from the ordinary Grenander.
Although, we have no theoretical evidence, up to this point there is no reason to think that bootstrapping from the ordinary
Grenander does not work.
\begin{figure}
\centering
\subfigure[Simulated powers (solid) from bootstrapping the Grenander.\label{fig:Grenander2[0,3]}]{\includegraphics[width=0.45\textwidth]{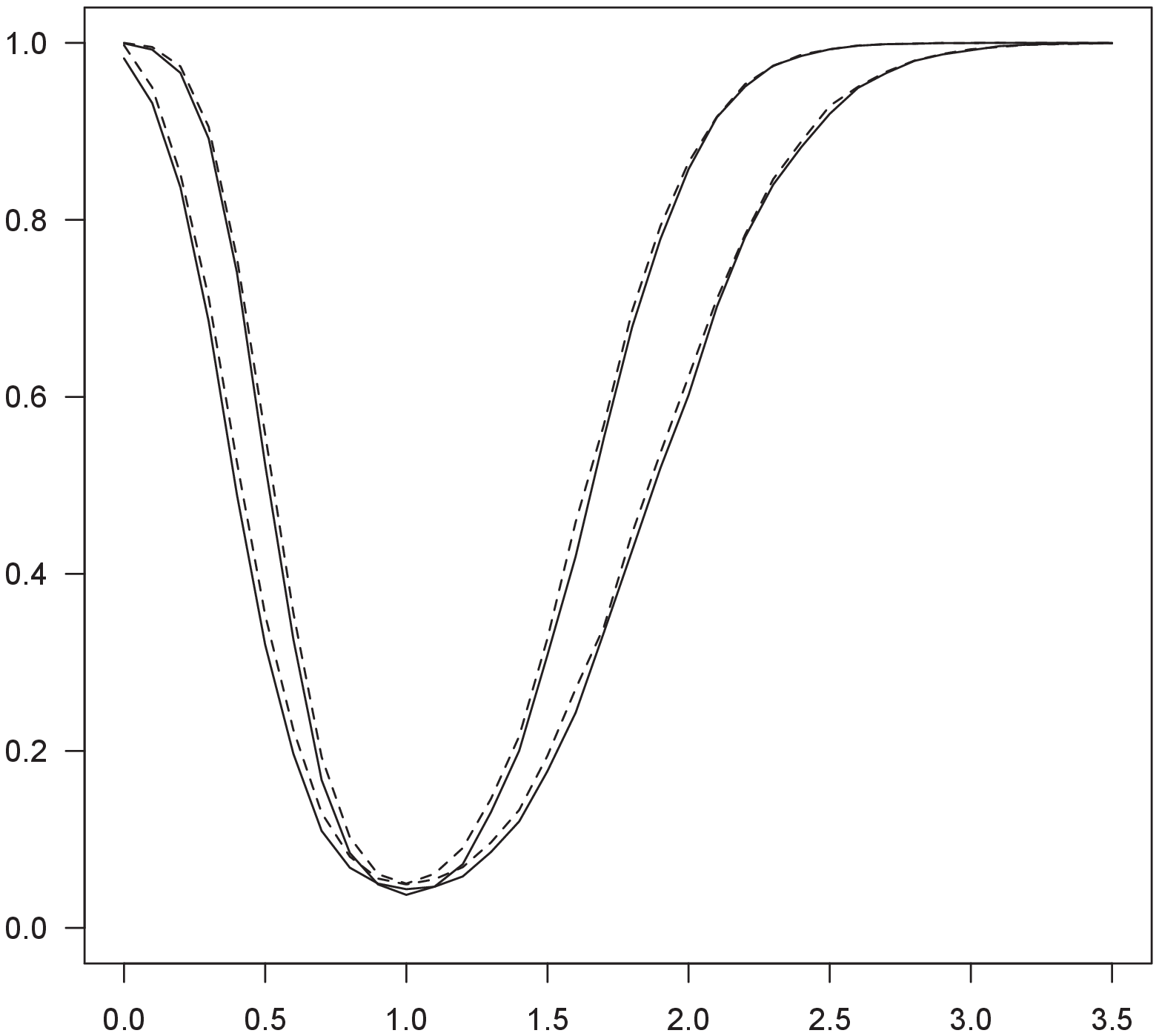}}
\hspace{1cm}
\subfigure[Simulated powers (solid) from smooth bootstrap.\label{fig:powers_smooth[0,3]}]{\includegraphics[width=0.45\textwidth]{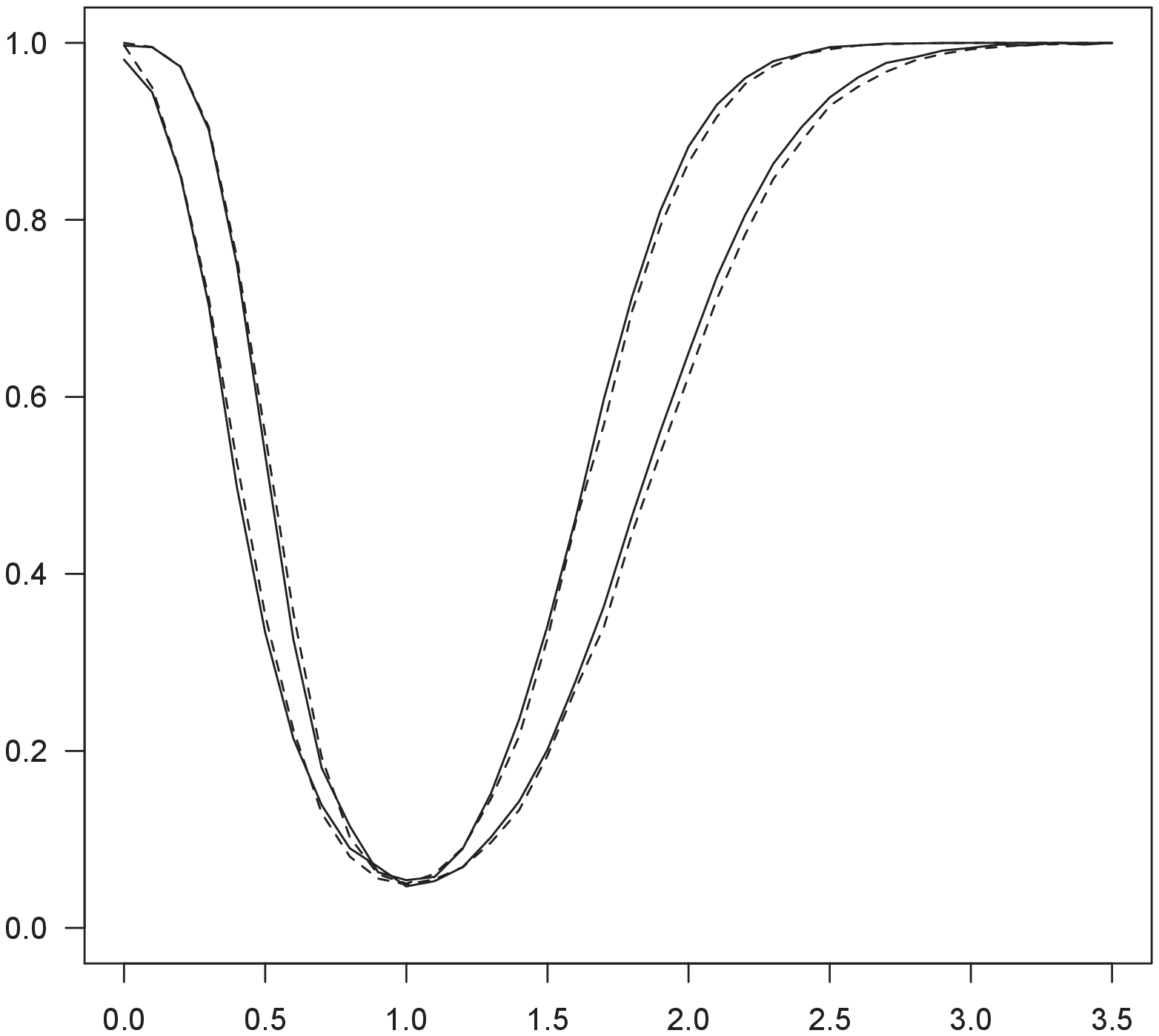}}
\caption{Simulated powers (solid) from bootstrapping and estimated true powers (dashed) of $S_{n1}$ and $S_{n2}$, for $\lambda=0,0.1,0.2,\dots,3.5$. The level of the test is taken to be $0.05$.}
\label{fig:powers[0,3]}
\end{figure}

\section{Appendix}
\label{sec:appendix}
\subsection{Proof of the lemmas in Section~\ref{sec:main}}
\label{subsec:proof lemmas}
The proof of Lemma~\ref{lem:embed0} is straightforward and has been put in the supplement.
The proof of Lemma~\ref{lem:switch} is along the lines of the proof of equality (21)
in~\cite{durot2007} and has also been put in the supplement.
Similarly, the proof of Lemma~\ref{lem:E to B} follows the same reasoning as
the proof of Corollary 3.1 in \cite{groeneboom-hooghiemstra-lopuhaa1999} and has been deferred to the supplement.

\bigskip

We proceed by establishing Lemma~\ref{lem:transition}
to make the transition to Brownian motion.
For this we first prove that under $f_0=f_1=\cdots=f_J$,
standardized slopes converge in distribution
to the slopes of the LCM of the process $W(s)-s^2+2xs$.
For $j=0,1,2,\ldots,J$ and $S=E,B,W$, define
\begin{equation}
\label{eq:def phi_nj^S}
\widehat{\phi}^S_{nj}(t)=n_j^{1/3}\left\{\hf_{nj}^S(t)-f_j(t)\right\}.
\end{equation}
Then $\widehat{\phi}^S_{nj}(t)$ is the slope at $s=0$ of the LCM of the process
\begin{equation}
\label{eq:def Z_nj^S}
\begin{split}
Z_{nj,t}^S(u)
&=
n_j^{2/3}
\left\{
M_{nj}^S\big(t,t+n_j^{-1/3}u\big]
+
F_j\big(t,t+n_j^{-1/3}u\big]-f_j(t)n_j^{-1/3}u
\right\}.
\end{split}
\end{equation}
For $j=0,1,\ldots,J$ and $t\in[0,1]$, define scaling constants
\[
A_j(t)=\frac{2^{1/3}}{|f_j'(t)|^{1/3}L'_j(t)^{1/3}}>0
\quad\text{ and }\quad
B_j(t)=\frac{4^{1/3}L'_j(t)^{1/3}}{|f_j'|^{2/3}}>0,
\]
and let
$I_{nj}(t)
=
\left\{
u:t+n_j^{-1/3}B_j(t)u\in[0,1]
\right\}$.
\begin{lemma}
\label{lem:conv slopes}
Assume \mart, \embed.
Suppose $f_0=f_1=\cdots=f_J$ and for $t\in(0,1)$, let
\[
\Phi_{nj,t}^S(x)=A_j(t)\widehat{\phi}^S_{nj}\big(t+n_j^{-1/3}B_j(t)x\big)
\]
where $\widehat{\phi}^S_{nj}$ is defined by~\eqref{eq:def phi_nj^S}.
Then, for $S=E,B,W$, and $x\in \bigcap_{j=0}^J I_{nj}(t)$ fixed,
the vector
$(\Phi_{n0,t}^S(x),\Phi_{n1,t}^S(x),\ldots,\Phi_{nJ,t}^S(x))$ converges in distribution to
$(\widetilde{\Phi}_{t0}(x),\Phi_1(x),\ldots,\Phi_J(x))$,
where
\[
\Phi_j(x)
=
\text{the slope at $s=x$ of the LCM of the process $W_j(s)-s^2+2xs$},
\]
where $W_1,W_2,\ldots,W_J$ are independent standard Brownian motions
and $\widetilde{\Phi}_{t0}$ is defined similarly with
the standard Brownian motion $\widetilde{W}_{t0}$ defined in~\eqref{eq:def Wt0}.
\end{lemma}
\noindent\textbf{Proof.}
For $j=0,1,\ldots,J$, $t\in(0,1)$ fixed and $a_j\in\mathbb{R}$, consider the event
\[
\Phi_{nj,t}^S(x)\leq a_j
\quad\Leftrightarrow\quad
\widehat{f}_{nj}^S\left(t+n_j^{-1/3}B_j(t)x\right)
\leq
f_j\left(t+n_j^{-1/3}B_j(t)x\right)+n_j^{-1/3}A_j(t)^{-1}a_j,
\]
which, according to~\eqref{eq:switch}, is equivalent to
\begin{equation}
\label{eq:argmax event}
B_j(t)^{-1}n_j^{1/3}
\left\{
\widehat{U}_{nj}^S
\left(
f_j(t+n_j^{-1/3}B_j(t)x)+n_j^{-1/3}A_j(t)^{-1}a_j
\right)-t
\right\}
\leq x,
\end{equation}
where $n_0=n$.
By~\eqref{eq:def Unj}, the left hand side of~\eqref{eq:argmax event} is the argmax over $u\in I_{nj}(t)$ of the process
\begin{equation}
\label{eq:def Znjx}
\begin{split}
n_j^{2/3}
\bigg\{
M_{nj}^S\left(t,t+n_j^{-1/3}B_j(t)u\right]
&+
F_j\left(t,t+n_j^{-1/3}B_j(t)u\right]\\
&-
f_j\left(t+n_j^{-1/3}B_j(t)x\right)n_j^{-1/3}B_j(t)u
\bigg\}\\
&-
A_j(t)^{-1}B_j(t)a_ju.
\end{split}
\end{equation}
To cover all cases $j=0,1,\ldots,J$ simultaneously, write
\begin{equation}
\label{eq:def xi_nj^S}
n_0=n
\quad\text{and define}\quad
\xi_{nj}^S(t)
=
\begin{cases}
\xi_{n0}(t) & ,j=0\\
\xi_{nj}t  & ,j=1,2,\ldots,J.
\end{cases}
\end{equation}
where $\xi_{n0}$ is defined in~\eqref{eq:def xi0}.
For $t+n_j^{-1/3}s\in[0,1]$, write
\[
\begin{split}
n_j^{2/3}M_{nj}^S\big(t,t+n_j^{-1/3}s\big]
&=
n_j^{1/6}
W_{nj}\circ L_j\big(t,t+n_j^{-1/3}s)\big]\\
&\quad+
n_j^{2/3}\left(M_{nj}^S-n_j^{-1/2}B_{nj}\circ L_{j}\right)\big(t,t+n_j^{-1/3}s\big]\\
&\qquad+
n_j^{1/6}\xi_{nj}^S\circ L_j\big(t,t+n_j^{-1/3}s\big],
\end{split}
\]
where $\xi_{nj}^S$ and $n_0$ are defined in~\eqref{eq:def xi_nj^S}.
According to (A3) and Lemma \ref{lem:embed0},
\[
\sup_{t+n_j^{-1/3}s\in[0,1]}
\left|
n_j^{2/3}\left(M_{nj}^S-n_j^{-1/2}B_{nj}\circ L_{j}\right)\big(t,t+n_j^{-1/3}s\big]
\right|
=
n_j^{2/3}O_p\left(n_j^{-1+1/q}\right)
=o_p(1).
\]
For every $j=1,2,\ldots,J$ and $k>0$ fixed, we have
\[
\sup_{|s|\leq k}
\left|n_j^{1/6}\xi_{nj}^S\circ L_j\big(t,t+n_j^{-1/3}s\big]\right|
\leq
kn_j^{-1/6}|\xi_{nj}|\sup_{t\in[0,1]}L_j'(t)=o_p(1).
\]
Furthermore
\[
n_0^{1/6}\xi_{n0}\circ L_0\big(t,t+n_0^{-1/3}s\big]=
n^{1/6}
\sum_{j=1}^J
c_j^{1/2}\xi_{nj}L_j\big(t,t+n_0^{-1/3}s\big].
\]
Hence,  we have
\[
\sup_{|u|\leq k}
\left|n_0^{1/6}\xi_{n0}\circ L_0\big(t,t+n^{-1/3}u\big]\right|
\leq
kn^{-1/6}
\sum_{j=1}^J c_j^{1/2}|\xi_{nj}|\sup_{t\in[0,1]}L_j'(t)=o_p(1).
\]
Finally, if we define
\begin{equation}
\label{eq:def Inj Wtj}
\begin{split}
I_{nj}(a)
&=
\left[
n_j^{1/3}(L_j(0)-L_j(g_j(a))),\, n_j^{1/3}(L_j(1)-L_j(g_j(a)))
\right],\\
W_{tj}(y)
&=
n_j^{1/6}
\left\{
W_{nj}\big(L_j(t)+n_j^{-1/3}y)-W_{nj}(L_j(t)\big)
\right\},
\end{split}
\end{equation}
with $W_{nj}$, for $j=1,2,\ldots,J$, being independent Brownian motions from~\eqref{eq:Bnj and Wnj}
and $W_{n0}$ is the Brownian motion defined by (\ref{eq:def Wn0}),
then for $j=1,2,\dots,J$,
\[
\begin{split}
Z_{nt,j}(s)
&=
n_j^{1/6}
W_{nj}\circ L_j\big(t,t+n_j^{-1/3}s)\big]\\
&=
W_{tj}\left(n_j^{1/3}\big(L_j(t+n_j^{-1/3}s)-L_j(t)\big)\right)
\approx
W_{tj}\left(L_j'(t)s\right).
\end{split}
\]
Because Brownian motion is uniformly continuous on compacta, it follows that
for each $j=1,2,\ldots,J$, the process $Z_{nj}(s)$ converges in the uniform topology on compacta to the process
$Z_{tj}(s)=L_j'(t)^{1/2}W_j(s)$,
where $W_1,W_2,\ldots,W_J$ are independent standard Brownian motions.
Furthermore, for $j=0$, according to~\eqref{eq:def Inj Wtj}, we have that
\begin{equation}
\label{eq:def Zn0}
Z_{nt,0}(s)
=
n^{1/6}
W_{n0}\circ L_0\big(t,t+n^{-1/3}s)\big]
=
\sum_{j=1}^J
c_j^{1/3}
Z_{nt,j}(c_j^{1/3}s),
\end{equation}
which converges in distribution to the process
\begin{equation*}
\widetilde{Z}_{t0}(s)=
\sum_{j=1}^J
c_j^{1/2}L_j'(t)^{1/2}W_j(s)
=
L_0'(t)^{1/2}\widetilde{W}_{t0}(s),
\end{equation*}
where $\widetilde{W}_{t0}$ is defined in~\eqref{eq:def Wt0}.
We then conclude that for each $j=1,2,\ldots,J$ and $S=E,B,W$,
the process in~\eqref{eq:def Znjx}
converges in the uniform topology on compacta to the process
\[
\begin{split}
&
L'_j(t)^{1/2}B_j(t)^{1/2}W_{j}\left(u\right)
-
\frac12|f_j'(t)|B_j(t)^2u^2+|f_j'(t)|B_j(t)^2xu
-
A_j(t)^{-1}B_j(t)a_ju\\
&=
L'_j(t)^{1/2}B_j(t)^{1/2}
\left\{
W_{j}(u)-u^2+2xu-a_ju
\right\},
\end{split}
\]
and to $L'_0(t)^{1/2}B_0(t)^{1/2}\{\widetilde{W}_{t0}(u)-u^2+2xu-a_0u\}$ in the case $j=0$.
According to Lemma~4 in~\cite{durot2007}, together with assumptions~\mono\ and~\mart,
the argmax on the left hand side
of~\eqref{eq:argmax event} is of order~$O_p(1)$.
This means we can apply Theorem~2.7 from~\cite{kim-pollard1990}.
Together with the fact that $\argmax\{H(u)\}=\argmax\{aH(u)+b\}$ for constants $a>0$ and $b\in \mathbb{R}$,
this yields that for each $j=0,1,\ldots,J$, $t\in(0,1)$ and~$S=E,B,W$ fixed, the argmax in~\eqref{eq:argmax event}
converges in distribution to $\nu_j(x-a_j/2)$ and $\widetilde{\nu}_{t0}(x-a_0/2))$, respectively, where
\begin{equation}
\label{eq:def Vj}
\begin{split}
\nu_j(c)
&=
\argmax_{u\in \mathbb{R}}
\left\{
W_j(u)-(u-c)^2
\right\},
\quad j=1,2,\ldots,J,\\
\widetilde{\nu}_{t0}(c)
&=
\argmax_{u\in \mathbb{R}}
\left\{
\widetilde{W}_{t0}(u)-(u-c)^2
\right\}.
\end{split}
\end{equation}
To extend this to joint convergence,
note that since the processes $Z_{nt,j}$, for $j=1,2,\ldots,J$, are independent
and $Z_{nt,0}$ satisfies~\eqref{eq:def Zn0}, they converge in distribution jointly to~$(\widetilde{Z}_{t0},Z_{t1},\ldots,Z_{tJ})$.
This implies joint convergence of the argmax's in~\eqref{eq:argmax event};
see~e.g., Theorem~6.1 in~\cite{huang-wellner1995},
which is only proven for two argmax's but which can trivially be extended to joint convergence
of more than two.
We conclude that
\[
\begin{split}
\prob\left(
\bigcap_{j=0}^J
\left\{
A_j(t)x)\Phi_{nj,t}^S(x)\leq a_j
\right\}
\right)
&\to
\prob\left(
\bigcap_{j=0}^J
\left\{
\nu_j\left(x-\frac{a_j}2\right)
\leq x
\right\}
\right)\\
&=
\prob\left(
\bigcap_{j=0}^J
\left\{
2\nu_j\left(0\right)
\leq a_j
\right\}
\right),
\end{split}
\]
using the fact that $\nu_j(c)-c\stackrel{d}{=}\nu_j(0)$, for all $j=0,1,\ldots,J$ and $c\in\mathbb{R}$.
Now, let
\[
D_{j}(s,x)=\text{left hand slope of $W_j(u)-u^2+2xu$ at $u=s$.}
\]
Then, similar to~\eqref{eq:switch}, one has $\nu_j(c)\leq t$ if and only if $D_j(t,0)\leq -2c$,
and it is straightforward to deduce that $2\nu_j(c)$ has the same distribution as $D_j(c,0)+4c$.
Furthermore, by properties of the LCM, one also has~$D_j(s,x)=D_j(s,0)+2x$.
It follows that
\[
\prob\left(
\bigcap_{j=0}^J
\left\{
2\nu_j\left(0\right)
\leq a_j
\right\}
\right)
=
\prob\left(
\bigcap_{j=0}^J
\left\{
D_j(x,x)
\leq a_j
\right\}
\right),
\]
and since $\Phi_j(x)=D_j(x,x)$, this proves the lemma.
\tqed

\noindent
Another ingredient to establish Lemma~\ref{lem:transition},
is a mixing property of the Brownian motion slope process.
\begin{lemma}
\label{lem:mixing}
Suppose that $f_0=f_1=\cdots,f_J$.
Then the process
\[
\big\{\big(\hf_{n0}^W(t),\hf_{n1}^W(t),\ldots,\hf_{nJ}^W(t)\big):t\in[0,1]\big\}
\]
is strong mixing.
More specifically, for $d>0$,
$$
\sup|\prob(A\cap B)-\prob(A)\prob(B)|\leq \alpha_{n}(d)=C_1e^{-C_2n d^3},
$$
where $C_1,C_2>0$ only depend on $f_0=f_1=\cdots=f_J$ and $c_1,c_2,\ldots,c_J$,
where the supremum is taken over all sets
$A\in\sigma\{\hf_{nj}^W(s):j=0,1,\ldots,J,\,0<s\leq t\}$ and
$B\in\sigma\{\hf_{nj}^W(u):j=0,1,\ldots,J,\,t+d\leq u<1\}$.
\end{lemma}
Its proof is along the lines of the proof of Lemma~4.6 in~\cite{kulikov-lopuhaa2008} and
has been put in the supplement.
Finally, we need the following result on the slopes of dependent Brownian motions with drift.
\begin{lemma}
\label{lem:cator}
For $i,j=0,1,\ldots,J$, $i\ne j$, $c>0$, and $t\in[0,1]$ fixed,
we have
$\prob\left(c\Phi_i(0)=\Phi_j(0)\right)=0$,
where $\Phi_1,\Phi_2,\ldots,\Phi_J$ and $\Phi_0=\widetilde{\Phi}_{t0}$
are defined in Lemma~\ref{lem:conv slopes}.
\end{lemma}
\noindent\textbf{Proof:}
When $i,j\geq 1$, the statement is trivially true,
because $\Phi_j(0)$ has the same distribution as $2\nu_j(0)$, as defined in~\eqref{eq:def Vj},
which has a bounded symmetric density according to Lemma~3.3 in~\cite{groeneboom-hooghiemstra-lopuhaa1999}.
Consider the case $i=0$ and $j\geq 1$.
Because $\widetilde{\Phi}_{t0}(0)$ also has the same distribution as $2\widetilde{\nu}_{t0}(0)$,
it is equivalent to prove that $c^{-1}\nu_j(0)=\widetilde{\nu}_{t0}(0)$ with probability zero.
By Brownian scaling, we have that
\[
c^{-1}\nu_j(0)
=
\argmax_{u\in\mathbb{R}}
\left\{
W_j(u)-c^{3/2}u^2
\right\}
\]
and from~\eqref{eq:def Wt0} it follows that
$\widetilde{\nu}_{t0}(0)=\argmax_{u\in\mathbb{R}}
\left\{
a_jW_j(u)+aW(u)-u^2
\right\}$,
where
\[
a_j=c_j^{1/2}\left(\frac{L_j'(t)}{L_{0}'(t)}\right)^{1/2},
\quad
a=\left(\sum_{j\neq i}c_{j}\frac{L_{j}'(t)}{L_{0}'(t)}\right)^{1/2},
\quad
W=\frac{1}{a}\sum_{m\neq j}c_{m}^{1/2}\left(\frac{L_{m}'(t)}{L_{0}'(t)}\right)^{1/2}W_{m}.
\]
Note that $W_j$ and $W$ are independent standard Brownian motions.
According to~\cite{cator2012}, with probability one there does not exist $u\in\mathbb{R}$ such that $u$ is a local maximum for
the process $X_1(u)=W_j(u)-c^{3/2}u^2$ and the process $X_2(u)=a_jW_j(u)+aW(u)-u^2$.
This proves the lemma.
\tqed

\noindent
\textbf{Proof of Lemma~\ref{lem:transition}.}
We follow the line of reasoning as in Corollary 3.3 in \cite{groeneboom-hooghiemstra-lopuhaa1999}.
To cover all cases $j=0,1,\ldots,J$ simultaneously, first introduce
\[
X_{nj}(s)=
\begin{cases}
\xi_{nj}L_j(s) & ,j=1,2,\ldots,J,\\
\xi_{n0}(L_0(s))& ,j=0,
\end{cases}
\]
where $\xi_{nj}$ is defined in~\eqref{eq:Bnj and Wnj}, for $j=1,2,\ldots,J$, and in~\eqref{eq:def xi0}, for $j=0$.
Note that according to~\eqref{eq:def xi0},
$X_{n0}(s)
=
\sum_{j=1}^J c_j^{1/2}X_{nj}(s)$.
Next, for $j=0,1,\ldots,J$, introduce the process
\[
Z_{nj,t}^\xi(s)=Z_{nj,t}^B(s)+n_j^{-1/6}X_{nj}'(t)s,
\]
where $Z_{nj,t}^B$ is defined in~\eqref{eq:def Z_nj^S} and $n_0=n$,
and denote $\widehat{\phi}_{nj}^\xi(t)$ as the slope of the least concave majorant
of~$Z_{nj,t}^\xi(s)$ at~$s=0$.
Then
\begin{equation}
\label{eq:xi slope}
\widehat{\phi}^\xi_{nj}(t)=\widehat{\phi}^B_{nj}(t)+n_j^{-1/6}X_{nj}'(t).
\end{equation}
Because
$Z_{nj,t}^B(s)=Z_{nj,t}^W(s)-n_j^{1/6}\{X_{nj}(t+n_j^{-1/3}s)-X_{nj}(t)\}$,
it follows that
\[
Z_{nj,t}^\xi(s)=Z_{nj,t}^W(s)-n_j^{1/6}\left(X_{nj}(t+n_j^{-1/3}s)-X_{nj}(t)-n_j^{-1/3}X_{nj}'(t)s\right).
\]
Let $[\tau_1,\tau_2]$ be the segment of the LCM of ${Z}_{nj,t}^\xi$ that contains zero,
and $[\tau_1',\tau_2']$ the segment of the LCM of ${Z}_{nj,t}^W$ that contains zero.
Define $a=\max(\tau_1,\tau_1')\leq 0$ and $b=\min(\tau_2,\tau_2')\geq 0$.
Note that we always have $a<b$, otherwise $\tau_1=0$ or $\tau_1'=0$, which is impossible by definition of the argmax.
Then for any $j=0,1,\ldots,J$,
\[
\widehat{\phi}_{nj}^\xi(t)
=
\frac{Z_{nj,t}^\xi(a)-Z_{nj,t}^\xi(b)}{a-b}
\quad\text{ and }\quad
\widehat{\phi}_{nj}^W(t)
=
\frac{Z_{nj,t}^W(a)-Z_{nj,t}^W(b)}{a-b},
\]
and
\begin{equation}
\label{eq:emb slope}
\left|
\widehat{\phi}_{nj}^\xi(t)
-
\widehat{\phi}_{nj}^W(t)
\right|
\leq
n_j^{-1/2}\sup_{s\in[0,1]}|X_{nj}''(s)|,
\end{equation}
uniformly in $t$.
This means that it remains to show that
\[
\int_0^1 |\widehat{\phi}^B_{ni}(t)-\widehat{\phi}^B_{nj}(t)|\,\md t
=
\int_0^1 |\widehat{\phi}^\xi_{ni}(t)-\widehat{\phi}^\xi_{nj}(t)|\,\md t+o_p(n^{-1/6}).
\]
However, if we define for $S=B,W$ and $i,j=0,1,\ldots,J$,
\[
\begin{split}
\widehat{\psi}_{ij}^S(t)
&=
\widehat{\phi}_{ni}^S(t)-\widehat{\phi}_{nj}^S(t)\\
X_{ij}(t)
&=
c_i^{-1/6}X_{ni}(t)-c_j^{-1/6}X_{nj}(t),
\end{split}
\]
where $c_0=1$, then according to~\eqref{eq:xi slope}, it is equivalent to show
\[
n^{1/6}
\int_0^1
\left\{
|\widehat{\psi}_{ij}^B(t)+n^{-1/6}X_{ij}'(t)|-|\widehat{\psi}_{ij}^B(t)|
\right\}
\,\md t
=
o_p(1).
\]
Let $\epsilon>0$.
Then
\[
\begin{split}
&n^{1/6}
\int_0^1
\left\{
|\widehat{\psi}_{ij}^B(t)+n^{-1/6}X_{ij}'(t)|-|\widehat{\psi}_{ij}^B(t)|
\right\}\,\md t\\
&=
n^{1/6}
\int_0^1
\left\{
|\widehat{\psi}_{ij}^B(t)+n^{-1/6}X_{ij}'(t)|-|\widehat{\psi}_{ij}^B(t)|
\right\}
1_{[0,\epsilon]}(|\widehat{\psi}_{ij}^B(t)|)\,\md t\\
&\quad+
n^{1/6}
\int_0^1
\left\{
|\widehat{\psi}_{ij}^B(t)+n^{-1/6}X_{ij}'(t)|-|\widehat{\psi}_{ij}^B(t)|
\right\}
1_{(\epsilon,\infty)}(|\widehat{\psi}_{ij}^B(t)|)\,\md t.
\end{split}
\]
Because of the independence between $\xi_{nj}$ and $B_{nj}$, and hence between $X_{nj}(t)$ and $\widehat{\phi}^B_{nj}(t)$, the expectation of the first
term on the right hand side is bounded from above by
\[
\sup_{s\in[0,1]}\mathbb{E}|X_{ij}'(s)|
\int_0^1
\mathbb{P}\left(|\widehat{\psi}^B_{ij}(t)|\leq \epsilon\right)\,\md t.
\]
According to Lemma~\ref{lem:conv slopes}, it follows that
for all $i,j=0,1,\ldots,J$,
\[
\mathbb{P}\left(|\widehat{\psi}^B_{ij}(t)|\leq \epsilon\right)
\to
\mathbb{P}
\left(
\left|
\frac{\Phi_i(0)}{c_{i1}(t)}-\frac{\Phi_j(0)}{c_{j1}(t)}
\right|
\leq \epsilon\right),
\]
where $\Phi_{0}$ is short for $\widetilde{\Phi}_{t0}$.
By right continuity and Lemma~\ref{lem:cator},
\[
\lim_{\epsilon\downarrow0}
\prob\left(\left|\frac{\Phi_{i}(0)}{c_{i1}(t)}-\frac{\Phi_{0}(0)}{c_{01}(t)}\right|\leq\epsilon\right)
=
\prob\left(\frac{\Phi_{i}(0)}{c_{i1}(t)}=\frac{\Phi_{0}(0)}{c_{01}(t)}\right)=0
\]
It follows that
\[
\lim_{\epsilon\downarrow0}
\,
\limsup_{n\to\infty}
\mathbb{E}
\left[
n^{1/6}
\int_0^1
\left\{
|\widehat{\psi}_{ij}^B(t)+n^{-1/6}X_{ij}'(t)|-|\widehat{\psi}_{ij}^B(t)|
\right\}
1_{[0,\epsilon]}(|\widehat{\psi}_{ij}^B(t)|)\,\md t
\right]
=0.
\]
For the remaining integral we write
\[
\begin{split}
&n^{1/6}
\int_0^1
\left\{
|\widehat{\psi}^B_{ij}(t)+n^{-1/6}X_{ij}'(t)|-|\widehat{\psi}^B_{ij}(t)|
\right\}
1_{(\epsilon,\infty)}(|\widehat{\psi}^B_{ij}(t)|)\,\md t\\
&=
\int_0^1
1_{(\epsilon,\infty)}(|\widehat{\psi}^B_{ij}(t)|)
\frac{2X_{ij}'(t)\widehat{\psi}^B_{ij}(t)+n^{-1/6}X_{ij}'(t)^2}{|\widehat{\psi}^B_{ij}(t)+n^{-1/6}X_{ij}'(t)|+|\widehat{\psi}^B_{ij}(t)|}
\,\md t\\
&=
\int_0^1
1_{(\epsilon,\infty)}(|\widehat{\psi}^B_{ij}(t)|)
\frac{2X_{ij}'(t)\widehat{\psi}^B_{ij}(t)}{|\widehat{\psi}^B_{ij}(t)+n^{-1/6}X_{ij}'(t)|+|\widehat{\psi}^B_{ij}(t)|}
\,\md t
+O_p(n^{-1/6})\\
&=
-
\int_0^1
X_{ij}'(t)\,
\sign(\widehat{\psi}^B_{ij}(t))1_{(\epsilon,\infty)}(|\widehat{\psi}^B_{ij}(t)|)
\,\md t
+O_p(n^{-1/6}),
\end{split}
\]
using the fact that for $|\widehat{\psi}^B_{ij}(t)|>\epsilon$,
\[
\left|
\frac{2\widehat{\psi}^B_{ij}(t)}{|\widehat{\psi}^B_{ij}(t)+n^{-1/6}X_{ij}'(t)|+|\widehat{\psi}^B_{ij}(t)|}
-
\frac{\widehat{\psi}^B_{ij}(t)}{|\widehat{\psi}^B_{ij}(t)|}
\right|
\leq
\frac{n^{-1/6}|X_{ij}'(t)|}{\epsilon}+O_p(n^{-1/6}).
\]
For $t\in[0,1]$ and $S=B,W$, let $\widehat{Y}_{ij}^S(t)=\sign(\widehat{\psi}^S_{ij}(t))1_{(\epsilon,\infty)}(|\widehat{\psi}^S_{ij}(t)|)$.
Then, again by independence of $\xi_{nj}$ and $B_{nj}$, we get
\[
\mathbb{E}
\left\{
\int_0^1
X_{ij}'(t)\widehat{Y}_{ij}^B(t)\,\md t
\right\}^2
=
2\iint_{0<s<t<1}
\mathbb{E}\left[X_{ij}'(s)X_{ij}'(t)\right]
\mathbb{E}\left[\widehat{Y}_{ij}^B(s)\widehat{Y}_{ij}^B(t)\right]\,\md s\,\md t,
\]
for the cases $j=1,2,\ldots,J$ and according to~\eqref{eq:jensen},
\[
\begin{split}
\mathbb{E}
\left\{
\int_0^1
X_{i0}'(t)\widehat{Y}_{i0}^B(t)\,\md t
\right\}^2
&\leq
J
\sum_{j=1}^J
c_j
\mathbb{E}
\left\{
\int_0^1
X_{ij}'(t)\widehat{Y}_{i0}^B(t)\,\md t
\right\}^2\\
&=
2J
\sum_{j=1}^J
c_j
\iint_{0<s<t<1}
\mathbb{E}\left[X_{ij}'(s)X_{ij}'(t)\right]
\mathbb{E}\left[\widehat{Y}_{i0}^B(s)\widehat{Y}_{i0}^B(t)\right]\,\md s\,\md t,
\end{split}
\]
where for $i,j=1,2,\ldots,J$,
\[
\mathbb{E}\left[X_{ij}'(s)X_{ij}'(t)\right]
=
\left(
\frac{L_i'(s)L_i'(t)}{c_i^{1/3}}+\frac{L_j'(s)L_j'(t)}{c_j^{1/3}}
\right).
\]
Furthermore, for all $i,j=0,1,\ldots,J$,
\[
\left|
\mathbb{E}\widehat{Y}_{ij}^B(s)\widehat{Y}_{ij}^B(t)-\mathbb{E}\widehat{Y}_{ij}^W(s)\widehat{Y}_{ij}^W(t)
\right|
\leq
\mathbb{E}|\widehat{Y}_{ij}^B(s)-\widehat{Y}_{ij}^W(s)|
+
\mathbb{E}|\widehat{Y}_{ij}^B(t)-\widehat{Y}_{ij}^W(t)|.
\]
For every $t\in[0,1]$, we have
\[
\begin{split}
\mathbb{E}|\widehat{Y}_{ij}^B(t)-\widehat{Y}_{ij}^W(t)|
&\leq
2\mathbb{P}(|\widehat{\psi}^B_{ij}(t)-\widehat{\psi}^W_{ij}(t)|>2\epsilon)
+
\mathbb{P}(|\widehat{\psi}^B_{ij}(t)|\leq\epsilon)
+
\mathbb{P}(|\widehat{\psi}^W_{ij}(t)|\leq\epsilon).
\end{split}
\]
Note that for all $j=0,1,\ldots,J$,
\begin{equation}
\label{eq:supp diff}
\sup_{t\in[0,1]}
\mathbb{E}|\widehat{\phi}^B_{nj}(t)-\widehat{\phi}^W_{nj}(t)|
=
O(n_j^{-1/6}).
\end{equation}
This can be shown similar to \eqref{eq:emb slope} using that
\[
Z_{nj,t}^B(s)=Z_{nj,t}^W(s)-n_j^{1/6}\left(X_{nj}(t+n_j^{-1/3}s)-X_{nj}(t)\right)
\]
and
\[
|X_{nj}(t+n_j^{-1/3}b)-X_{nj}(t+n^{-1/3}a)|
\leq
n_j^{-1/3}(b-a)\sup_{t\in[0,1]} |X_{nj}'(t)|.
\]
The Markov inequality together with \eqref{eq:supp diff} yields that
$\mathbb{P}(|\widehat{\psi}^B_{ij}(t)-\widehat{\psi}^W_{ij}(t)|>2\epsilon)$ tends to zero,
uniformly in $t\in[0,1]$.
As before, according to Lemma~\ref{lem:conv slopes},
\[
\lim_{\epsilon\downarrow0}
\,
\limsup_{n\to\infty}
\mathbb{P}(|\widehat{\psi}^B_{ij}(t)|\leq\epsilon)
=0
\]
and likewise for $\mathbb{P}(|\widehat{\psi}^W_{ij}(t)|\leq\epsilon)$.
We conclude that for all $0<s<t<1$, $i,j=0,1,\ldots,J$ and~$\epsilon>0$,
\[
\lim_{\epsilon\downarrow0}
\,
\limsup_{n\to\infty}
\left|
\mathbb{E}\widehat{Y}_{ij}^B(s)\widehat{Y}_{ij}^B(t)-\mathbb{E}\widehat{Y}_{ij}^W(s)\widehat{Y}_{ij}^W(t)
\right|
=0.
\]
Finally, for all $j=0,1,\ldots,J$, write
\[
\mathbb{E}\left[\widehat{Y}_{ij}^W(s)\widehat{Y}_{ij}^W(t)\right]
=
\text{cov}\left(\widehat{Y}_{ij}^W(s),\widehat{Y}_{ij}^W(t)\right)+\mathbb{E}[\widehat{Y}_{ij}^W(s)]\mathbb{E}[\widehat{Y}_{ij}^W(t)].
\]
Because $(\widehat{\phi}^W_{ni}(t),\widehat{\phi}^W_{nj}(t))$ is strong mixing according to
Lemma~\ref{lem:mixing}, also $\widehat{\psi}^W_{ij}(t)$ is strong mixing.
Then according to~\cite{ibragimov-linnik1971} (see also Lemma 3.1 in \cite{groeneboom-hooghiemstra-lopuhaa1999}),
for every $0<s<t<1$ we get that
\[
\left|
\text{cov}\left(\widehat{Y}_{ij}^W(s),\widehat{Y}_{ij}^W(t)\right)
\right|
\leq
Ke^{-C_2n(t-s)}\to0.
\]
Also for every $t\in(0,1)$, according to Lemma~\ref{lem:conv slopes},
writing $V_0$ and $\Phi_0$ for $\widetilde{\nu}_{t0}$ and $\widetilde{\Phi}_{t0}$,
\[
\begin{split}
\mathbb{E}[\widehat{Y}_{ij}^W(t)]
&=
\prob(\widehat{\psi}^W_{ij}(t)>\epsilon)-\prob(\widehat{\psi}^W_{ij}(t)<-\epsilon)\\
&\to
\prob\left(
\frac{\Phi_i(0)}{c_{1i}(t)}-\frac{\Phi_j(0)}{A_j(t)}
> \epsilon\right)
-
\prob\left(
\frac{\Phi_i(0)}{c_{1i}(t)}-\frac{\Phi_j(0)}{A_j(t)}
<-\epsilon
\right)\\
&=
\prob\left(
\frac{2V_i(0)}{c_{1i}(t)}-\frac{2\nu_j(0)}{A_j(t)}
> \epsilon\right)
-
\prob\left(
\frac{2V_i(0)}{c_{i1}(t)}-\frac{2\nu_j(0)}{A_j(t)}
<-\epsilon
\right)=0,
\end{split}
\]
because $(-V_i(0),-\nu_j(0))$ has the same distribution as $(V_i(0),\nu_j(0))$.
It follows that
\[
\lim_{\epsilon\downarrow0}
\,
\limsup_{n\to\infty}
\mathbb{E}
\left\{
\int_0^1
X_{ij}'(t)\widehat{Y}_{ij}^B(t)\,\md t
\right\}^2
=0.
\]
This proves the lemma.
\hfill
\tqed
\subsection{Proof of Theorems~\ref{th:main1} and~\ref{th:main2}}
\label{subsec:proof theorem1}
In this section we assume that assumptions \indep, \mono, \mart, \embed, \regul\ hold.
Now that we have established~\eqref{eq:transition complete} thanks to Lemmas~\ref{lem:switch}, \ref{lem:E to B} and~\ref{lem:transition},
we investigate integrals of the type
\[
n^{1/3}\int_{f(1)}^{f(0)}|\widehat{U}_{ni}^W(a)-\widehat{U}_{nj}^W(a)|\,\md a.
\]
We proceed as in step 2 in~\cite{durot2007} and approximate
$n_j^{1/3}\left(L_j(\widehat U_{nj}^W(a)-L_j(g_j(a))\right)$ by~$\widetilde{V}_{nj}(g_j(a))$,
where for all $j=0,1,\dots,J$,
\begin{equation}
\label{eq:def tilde Vnj}
\widetilde{V}_{nj}(t)=
\argmax_{|u|\leq\log n}
\left\{
W_{tj}(u)
-
\frac{|f_j'(t)|}{2(L_j'(t))^2}u^2
\right\}
\end{equation}
with $W_{tj}$ the Brownian motion defined in~\eqref{eq:def Inj Wtj}.

\begin{lemma}
\label{lem:transitionV}
There exists a $C>0$, such that for each $j=0,1,\dots,J$ and $a\in[f_j(1)+n^{-1/3}(\log n)^2, f_j(0)-n^{-1/3}(\log n)^2]$,
$$\E\left|n_j^{1/3}
\left(
L_j(\widehat{U}_{nj}^W(a))-L_j(g_j(a))
\right)
-\widetilde V_{nj}(g_j(a))\right|\leq C\frac{n^{-1/6}}{\log n}.$$
\end{lemma}
The proof is along the lines of the proof of step 2 in~\cite{durot2007} and has been put in the supplement.
Combining this with~\eqref{eq:transition complete} yields the following lemma, whose proof is straightforward and has been
put in the supplement.
\begin{lemma}
\label{lem:11}
Assume $f_0=f_1=\cdots=f_J$.
Then for every $i,j=0,1,\dots,J$, we have
$$
n^{1/3}\int_0^1
|\widehat{f}_{ni}^E(t)-\widehat{f}_{nj}^E(t)|\,\md t=\int_{0}^{1}
\left|
\frac{\widetilde{V}_{ni}(t)}{L_i'(t)c_i^{1/3}}
-
\frac{\widetilde{V}_{nj}(t)}{L_j'(t)c_j^{1/3}}
\right|
|f_0'(t)|\,\md t
+
o_p(n^{-1/6}).$$
\end{lemma}

\bigskip

\noindent
From Lemma~\ref{lem:11} we conclude that under $f_0=f_1=\cdots=f_J$,
for $k=1,2$, the test statistic~$S_{nk}$ satisfies
$$
n^{1/3}S_{nk}=\int_0^1Y_{nk}(t)\,\md t+o_p(n^{-1/6}),
$$
where
\begin{equation}
\label{eq:def Ynk}
\begin{split}
Y_{n1}(t)
&=
\sum_{i<j} \left|
\frac{\widetilde{V}_{ni}(t)}{L_i'(t)c_i^{1/3}}
-
\frac{\widetilde{V}_{nj}(t)}{L_j'(t)c_j^{1/3}}
\right|
|f_0'(t)|,\\
Y_{n2}(t)
&=
\sum_{j=1}^J \left|
\frac{
\widetilde{V}_{n0}(t)}{L_0'(t)}
-
\frac{\widetilde{V}_{nj}(t)}{L_j'(t)c_j^{1/3}}
\right|
|f_0'(t)|.
\end{split}
\end{equation}
Therefore, in order to prove Theorems~\ref{th:main1} and~\ref{th:main2}, it remains to show that under $f_0=f_1=\cdots=f_J$,
for $k=1,2$,
$$
n^{1/6}\left(\int_0^1Y_{nk}(t)\,\md t-m_k\right)
$$
converges in distribution to a centered Gaussian variable with a finite variance~$\sigma^2_k$.
To determine $m_k$ and $\sigma^2_k$, we have to deal with joint distributions of $\widetilde{V}_{ni}(t)$ and $\widetilde{V}_{nj}(t)$,
for different $i,j=0,1,\ldots,J$, and
with covariances between $\widetilde{V}_{ni}(s)$ and~$\widetilde{V}_{ni}(t)$,
for $s$ and $t$ close to each other.

We will approximate $\widetilde{V}_{nj}(s)$ and $\widetilde{V}_{nj}(t)$ with the variable~$V_{tj}(s)$ defined as follows.
For all  $t\in(0,1)$ and for $j=1,2,\dots,J$, we define
\begin{equation}
\label{eq:def Vtj}
V_{tj}(s)=\argmax_{u\in\R}\left\{W_{tj}(u)-\frac{|f_0'(s)|}{2(L_j'(s))^2}u^2\right\},
\end{equation}
where the process $W_{tj}(u)$ is defined in~\eqref{eq:def Inj Wtj}.
Recall that the processes $W_{tj}$, for $j=1,2,\ldots,J$, are independent Brownian motions,
whose joint distribution of $(W_{t1},W_{t2},\ldots,W_{tJ})$ does not depend on $n$.
Note that from~\eqref{eq:def Inj Wtj} it follows that
\begin{equation}
\label{eq:Wt0 in Wtj}
W_{t0}(u)
=
\sum_{j=1}^J
c_j^{1/3}
W_{tj}
\left(
n_j^{1/3}
\left\{
L_j\circ L_0^{-1}(L_0(t)+n^{-1/3}u)
-
L_j(t)
\right\}
\right).
\end{equation}
Although, due to~\eqref{eq:defL0}, $W_{t0}$ itself is distributed as standard Brownian motion,
the joint distribution of $(W_{t0},W_{t1},\ldots,W_{tJ})$ does depend on $n$.
For this reason, we approximate $W_{t0}$ with
\begin{equation}
\label{eq:def tildeWtj}
\widetilde W_{t0}(u)=\sum_{j=1}^Jc_j^{1/3}W_{tj}\left(c_j^{1/3}\frac{L_j'(t)}{L_0'(t)}u\right),
\end{equation}
and define
\begin{equation}
\label{eq:def Vt0}
V_{t0}(s)=\argmax_{u\in\R}\left\{\widetilde W_{t0}(u)-\frac{|f_0'(s)|}{2(L_0'(s))^2}u^2\right\}.
\end{equation}
Note that from~\eqref{eq:defL0} it follows that also $\widetilde W_{t0}$ is distributed like a standard Brownian motion,
but this time it is a linear combination of the $W_{tj}$'s not depending on $n$,
so that the joint distribution of $(\widetilde{W}_{t0},W_{t1},W_{t2},\ldots,W_{tJ})$
is independent of $n$.
The latter is important to determine an expression for the limiting variance $\sigma^2_2$ not depending on~$n$.
On the other hand, the approximation of~$\widetilde{V}_{n0}(t)$ by $V_{t0}(t)$ is not sufficient to replace
the expectation of $n^{1/3}S_{n2}$ by a constant~$m_2$ not depending on $n$.
For this we will use
\begin{equation}
\label{eq:def Vt0'}
V_{t0}'(s)=\argmax_{u\in\R}\left\{W_{t0}(u)-\frac{|f_0'(s)|}{2(L_0'(s))^2}u^2\right\},
\end{equation}
where $W_{t0}$ is from~\eqref{eq:Wt0 in Wtj}.
The following lemma provides approximations of $\widetilde{V}_{nj}(t)$.
Its proof is somewhat technical, but uses the same kind of reasoning
as in step 3 in~\cite{durot2007} and has been deferred to the supplement.
\begin{lemma}
\label{lem:approx tilde Vn}
For $r>1$,
\begin{equation}
\label{eq:tildeV VV'}
\begin{split}
\E^{1/r}\left|\widetilde{V}_{nj}(t)-V_{tj}(t)\right|^r
&=
o(n^{-1/6}),
\quad
\text{for $j=1,2\ldots,J$,}\\
\E^{1/r}\left|\widetilde{V}_{n0}(t)-V_{t0}'(t)\right|^r
&=
o(n^{-1/6}),
\end{split}
\end{equation}
uniformly in $t\in(0,1)$.
Furthermore, let $A>0$ and $r\in(1,2\holder)$, where $\holder>3/4$ is taken from~\regul.
Then,
\begin{equation}
\label{eq:tildeV V logn}
\E^{1/r}\left|\widetilde{V}_{nj}(t)-V_{tj}(s)\right|^r=o(\log n)^{-1},
\end{equation}
uniformly in $j=0,1,\ldots,J$, and $s,t\in(0,1)$, such that $|s-t|\leq A n^{-1/3}\log n$.
\end{lemma}

In the following lemma,
we prove that the variance of the above variable has a finite limit under $f_1=f_2=\cdots=f_J$.
As the proof follows the line of reasoning used in the proof of step 5 in~\cite{durot2007},
we only present a sketch of the proof.
A detailed proof can be found in the supplement.
\begin{lemma}
\label{lem:limit variance}
For $k=1,2$, let
$$
v_{nk}=\mathrm{var}\left(\int_0^1Y_ {nk}(t)\,\md t\right).
$$
Under $f_1=f_2=\cdots=f_J$, $n^{1/3}v_{nk}$ has a finite limit $\sigma_k^2$,
as $n\to\infty$, where $\sigma_1^2$ is defined in Theorem~\ref{th:main1} and $\sigma_2^2$ is defined in
Theorem~\ref{th:main2}.
\end{lemma}
\textbf{Sketch of proof.}
For $k=1,2$, we have
$$
v_{nk}=2\int_0^1\int_s^1
\mathrm{cov}(Y_ {nk}(t),Y_ {nk}(s))\,\md t\,\md s.
$$
Note that by definition of $\widetilde V_{nj}(t)$ in~\eqref{eq:def tilde Vnj},
the random variable $Y_{nk}(t)$, defined in~\eqref{eq:def Ynk}, depends only on the increments of $W_{nj}$ over a neigbourhood
of $L_j(t)$ with radius of the order~$O(n^{-1/3}\log n)$,
for $j=0,1,\dots,J$.
But for every $s,t\in[0,1]$ and $j=0,1,\dots,J$, we have
$|L_j(t)-L_j(s)|\geq|t-s|\inf_{u\in[0,1]}|L_j'(u)|$,
where the infimum is positive according to~\embed.
Setting $a_n=An^{-1/3}\log n$, for some large enough $A>0$,
we find that $Y_{nk}(t)$ is independent of $Y_{nk}(s)$ for every $|t-s|\geq a_n$.
This means that in for $k=1$, together with Lemma~\ref{lem:approx tilde Vn},
\begin{equation}\notag
\begin{split}
v_{n1}
&=
2\int_0^1\int_s^{1\wedge(s+a_n)}\mathrm{cov}(Y_ {n1}(t),Y_ {n1}(s))\,\md t\,\md s\\
&=
2\sum_{i<j}\sum_{l<m}
\int_0^1\int_s^{1\wedge(s+a_n)}
|f_0'(s)|^2C_{iljm}(s,t)\,\md t\,\md s+o(n^{-1/3}),
\end{split}
\end{equation}
where
$$
C_{iljm}(s,t)
=
\mathrm{cov}\left(\left|
\frac{{V}_{ti}(s)}{L_{i}'(s)c_i^{1/3}}
-
\frac{{V}_{tj}(s)}{L_{j}'(s)c_j^{1/3}}
\right|,\left|
\frac{{V}_{sl}(s)}{L_{l}'(s)c_{l}^{1/3}}
-
\frac{{V}_{sm}(s)}{L_{m}'(s)c_{m}^{1/3}}
\right|\right).
$$
Let $d_j(s)=|f_0'(s)|/(2L_j'(s)^2)$.
From~\eqref{eq:def Vtj} and~\eqref{eq:def Inj Wtj}, we have for $j=1,2,\ldots,J$,
that $d_j(s)^{2/3}V_{tj}(s)$ has the same distribution as
\[
\argmax_{u\in\R}
\left\{
W_{j}\left(
u+n^{1/3}(t-s)|f_0'(s)/2|^{2/3}\frac{c_j^{1/3}}{L_j'(s)^{1/3}}
\right)
-
u^2
+
R_n'(s,t,u)
\right\}
\]
where $W_j$ is a standard Brownian motion and where for every $|t-s|\leq a_n$,
and $R_n'(s,t,u)$ can be shown to be negligible.
For $j=1,2,\ldots,J$, let $\zeta_j$ be defined by~\eqref{eq: def zetaj}.
We then conclude that for $j=1,2,\ldots,J$,
$$
\mathbb{E}
\left|
\frac{V_{tj}(s)}{L_{j}'(s)c_j^{1/3}}
-
\left(\frac{4L_{j}'(s)}{c_j|f_0'(s)|^2}\right)^{1/3}
\zeta_j\left(n^{1/3}(t-s)|f_0'(s)/2|^{2/3}\frac{c_j^{1/3}}{L_j'(s)^{1/3}}\right)
\right|
=
o(1/\log n).
$$
Change of variable $t'=n^{1/3}(t-s)|f_0'(s)/2|^{2/3}$, then gives
\[
\begin{split}
n^{1/3}v_{n1}
=
8\sum_{i<j}\sum_{l<m}
\int_0^1\int_0^{a_n'}
\mathrm{cov}&\left(|Y_{si}(t')-Y_{sj}(t')|,|Y_{sl}(0)-Y_{sm}(0)|\right)\,\md t'\,\md s
+
o(1),
\end{split}
\]
where $a_n'=A\log n|f_0'(s)/2|^{2/3}$ and where for $j=1,2,\ldots,J$,
$Y_{sj}(t)$ is defined in~\eqref{eq:def Ysj}.
We finish the proof for the case $k=1$, by showing that
there exist absolute constants $K$ and~$K'$ such that
\begin{equation}\notag
|\mathrm{cov}\left(|Y_{si}(t)-Y_{sj}(t)|, |Y_{sl}(0)-Y_{sm}(0)|\right)|\leq K\exp(-K't^3),
\end{equation}
because then, for $k=1$, the lemma follows from the dominated convergence theorem.
The proof for the case $k=2$ is similar.
\tqed

\noindent
We are now in position to complete the proofs of Theorems~\ref{th:main1} and~\ref{th:main2}.

\bigskip

\noindent
\textbf{Proof of Theorems~\ref{th:main1} and~\ref{th:main2}.}
Define
\[
m_1:=\sum_{i<j}\int_{0}^{1}
\E\left|
\frac{{V}_{ti}(t)}{L_{i}'(t)c_i^{1/3}}
-
\frac{{V}_{tj}(t)}{L_{j}'(t)c_j^{1/3}}
\right|
|f_0'(t)|\,\md t
\]
and
\[
m_2:=\sum_{j=1}^J\int_{0}^{1}
\E\left|
\frac{{V}_{ti}(t)}{L_{i}'(t)c_i^{1/3}}
-
\frac{V_{t0}'(t)}{L_0'(t)}
\right|
|f_0'(t)|\,\md t,
\]
where $V_{tj}(t)$ and $V_{t0}'(t)$ are defined in~\eqref{eq:def Vtj} and~\eqref{eq:def Vt0'}.
From Lemma~\ref{lem:approx tilde Vn} we then obtain for $k=1,2$,
$$
n^{1/6}\left(m_k-\E\int_0^1Y_{nk}(t)\,\md t\right)=o(n^{-1/6}).
$$
Note that for $j=1,2,\ldots,J$,
\[
\frac{|f_0'(t)|}{L_j'(t)c_j^{1/3}}V_{tj}(t)
\stackrel{d}{=}
|4f_{0}'(t)|^{1/3}
\frac{L_j'(t)^{1/3}}{c_j^{1/3}}\zeta_j(0)
=
|4f_{0}'(t)|^{1/3}
Y_{sj}(0),
\]
where $\zeta_j$ and $Y_{sj}$ are defined in~\eqref{eq: def zetaj} and~\eqref{eq:def Ysj},
so that
\[
m_1=
\sum_{i<j}
\int_{0}^{1}
|4f_{0}'(t)|^{1/3}
\E\left|
Y_{si}(0)-Y_{sj}(0)
\right|
\,\md t
\]
Similarly
\[
\frac{|f_0'(t)|}{L_0'(t)}V_{t0}'(t)
\stackrel{d}{=}
|4f_{0}'(t)|^{1/3}
L_0'(t)^{1/3}\widetilde{\zeta}_{t0}'(0),
\]
where $\widetilde{\zeta}_{t0}'$ is defined in~\eqref{eq: def zetat0}, so that
\[
m_{2}
=
\sum_{j=1}^J
\int_{0}^{1}
|4f_{0}'(t)|^{1/3}
\E\left|
L_0'(t)^{1/3}\widehat{\zeta}_{t0}(0)
-
\frac{L_j'(t)^{1/3}}{c_j^{1/3}}\zeta_j(0)
\right|
\,\md t.
\]
Therefore, in order to prove Theorems~\ref{th:main1} and~\ref{th:main2}, it remains to show that, for $k=1,2$,
under $f_1=f_2=\cdots=f_J$,
$$
n^{1/6}\left(\int_0^1Y_{nk}(t)\,\md t-\E\int_0^1Y_{nk}(t)\,\md t\right)
$$
converges in distribution to a centered Gaussian variable with variance $\sigma_k^2$.
This can be done using the method of big blocks and small blocks,
similar to the proof of Theorem~4.1 in~\cite{groeneboom-hooghiemstra-lopuhaa1999}.
The details are omitted.
\tqed

\subsection{Proofs for Section~\ref{sec:calibration}}
\label{subsec:appendix examples}
The following Kiefer-Wolfowitz type result for our general setting is proved in~\cite{durot-lopuhaa2013} and will be
used to prove Lemma~\ref{lem: brevefn}.
\begin{theorem}
\label{theorem:KW}
Assume  \mono, \mart, {\embed} with some $q\geq 3$, and {\regul} with some $\holder>0$.
For $S=E,B,W$, let $\hat F_n^S$ be the least concave majorant of $F_{n}^S$.
Then
$$
\sup_{x\in[a,b]}|\hat F_n^S(x)-F_{n}^{S}(x)|=O_{p}\left(n^{-2/3}(\log n)^{2/3}\right).
$$
\end{theorem}
\noindent\textbf{Proof of  Lemma~\ref{lem: brevefn}.}
For convenience, we denote by $\widetilde f_n(t)$ the estimator defined by~\eqref{eq:tildefmain} for all $t\in[a+h_{n},b-h_{n}]$
and either by~\eqref{eq:tildefKosorok} or by~\eqref{def:boundary kernel estimate} on the boundaries $[a,a+h_{n})$ and $(b-h_{n},b_n]$,
and we denote by $\widehat f_{n}$ the estimator defined in the same manner as $\widetilde f_{n}$ with~$F_{n0}$ replaced by~$\widehat{F}_{n0}$.
For $l=0,1, 2$, we have
\begin{equation}
\begin{split}\label{eq:smoothgren}
&
\sup_{t\in[a+h_{n},b-h_{n}]}\left|
\widehat{f}_n^{(l)}(t)
-
\widetilde{f}_n^{(l)}(t)
\right|\\
&\quad=
\left|
\frac{1}{h_n^{1+l}}\int_\R (\widehat{F}_{n0}(t-uh_n)-F_{n0}(t-uh_n))K^{(l+1)}(u)\,\md u
\right|\\
&\quad\leq
\frac{1}{h_n^{1+l}}\sup_{s\in[a,b]}|\widehat{F}_{n0}(s)-F_{n0}(s)|
\int_{-1}^1 |K^{(l+1)}|(u)|\,\md u\\
&\quad=
h_{n}^{-(l+1)}O_{p}\left(n^{-2/3}(\log n)^{2/3}\right).
\end{split}
\end{equation}
Moreover, in the proof of Lemma~\ref{lem: tildefn}, it is proved that both functions $\phi$ and $\psi$ are bounded,
and that the supremum norm of $\widetilde f_{n}''$ is of order $h_{n}^{-5/2}n^{-1/2}\sqrt{\log (1/h_{n})}$
if~$f_{0}$ is twice continuously differentiable
Hence, one easily derives Lemma~\ref{lem: brevefn} from~\eqref{eq:smoothgren} and Lemma~\ref{lem: tildefn}.
\tqed

\noindent
\textbf{Sketch of proof for Theorems \ref{th:regression}, \ref{th:density}, \ref{th:censoring} and \ref{th:censoring2}. }
Note first that it suffices to prove the results in the particular case $[a,b]=[0,1]$, see Remark~\ref{remark}, so in the sequel, we consider only $[a,b]=[0,1]$.
That assumptions \mono,\dots,\regul~ are fulfilled in all the considered models is proved in \cite{durot2007},
see her Theorems 3, 5, 6.
In the monotone regression model, the embedding is with a Brownian motion and $L_j(t)=t\tau_j^2$;
in the density model, the embedding is with a Brownian Bridge and $L_j=F_j$;
in the random censorship model, the embedding is with a Brownian motion and
$$
L_j(t)= \int_0^t\frac{f_j(x)}{(1-G_j(x))(1-H_j(x))}\,\md x.
$$
Moreover, Assumption~\indep~ is clearly satisfied since all the original observations are mutually independent,
so Theorems~\ref{th:main1} and~\ref{th:main2} apply in all these models and it suffices to prove that (\ref{eq:cvboot}) also holds under $H_0$. Hereafter, we assume $f_1=\dots=f_J=f_0$ and $\widetilde f_n$ satisfies \monostar.

Note that, in order to prove convergence in probability, we can restrict ourselves to an event whose probability tends to one as $n\to\infty$. Thus, thanks to Lemma~\ref{lem: tildefn}, we assume in the sequel that $\widetilde f_n$ is decreasing on $[0,1]$ and satisfies (\ref{eq: tildef'}), (\ref{eq: approxtildef}), (\ref{eq: approxtildef'}) and (\ref{eq: boundtildef'}) for some positive $C_0,C_1$, $s\in(3/4,1]$, and $\varepsilon_n$ such that $n^{-\gamma}/\varepsilon_n$ tends to zero as $n\to\infty$ for any $\gamma>0$. Moreover, recall that a sequence of random variables $X_n$ converges in probability to a random variable $X$ if, and only if, every subsequence has a further subsequence along which $X_n$ converges almost surely to $X$. Thus, in the three considered models, if a sequence of random variables, which is measurable with respect to the original observations, converges in probability, we can assume for simplicity that it converges almost surely (otherwise, argue along subsequences). Then, we aim to prove that almost surely,
$n^{1/6}(S_{nk}^\star-m_k)$
converges in distribution to the Gaussian law with mean zero and variance $\sigma_k^2.$

In each model, we define the bootstrap versions of $f_0$ and $F_{nj}$ as $\widetilde f_n$ and~$F_{nj}^\star$ respectively, while the bootstrap version of $M_{nj}$ is $M_{nj}^\star=F_{nj}^\star-\widetilde F_n$ with
$\widetilde F_n(t)=\int_0^t\widetilde f_n(u)\,\md u$.
From what precedes, a bootstrap version of assumptions \indep~ and \mono~ holds (we mean that conditionally on the original observations, these assumptions hold with $f_0$ and $F_{nj}$ replaced by their bootstrap version). Moreover,
it can be proved following the line of reasonning used in the proof of Theorems 3,5,6 in \cite{durot2007}, that a
bootstrap version of assumptions \mart~ and \embed~hold. The boostrap version $L_j^\star$ of $L_j$ we consider in each considered model is given below.

\begin{enumerate}
\item In the regression model, $L_j^\star(t)=t\widehat\tau_j^2$ with $\widehat\tau_j^2$ the conditional variance of $\epsilon_{ij}^\star$. Let us notice that the conditional moment of order $q$ of $\epsilon_{ij}^\star$ is equal to
$$\E^\star|\epsilon_{ij}^\star|^q=\frac{1}{n_j}\sum_{i=1}^{n_j}|\widehat\epsilon_{ij}-\bar\epsilon_j|^q.$$
Under the assumptions of Theorem \ref{th:regression}, we have $\sup_t|\widetilde f_n(t)-f_0(t)|=o_p(1)$ so it follows from the law of large numbers that $\bar\epsilon_j=o_p(1)$ for all $j=1,\dots,J$. Using again the law of large numbers we obtain that
$$\max_{ij}\E^\star|\epsilon_{ij}^\star|^q= \max_{ij}\E|\epsilon_{ij}|^q+o_p(1).$$
Thus, we can assume here that $\max_{i,j}\E^\star|\epsilon_{ij}^\star|^q\leq C$
for some positive number $C$ that does not depend on $n$, which is the main ingredient to establish the bootstrap version of \embed~ in the regression model. Likewise, we can assume that $\widehat\tau_j^2=\tau_j^2+o(n^{-1/3})$ for all $j=1,\dots,J$.
\item In the density model,  $L_j^\star=\widetilde F_n$ for every $j$. Note that $f_0$ is assumed bounded from above and from below, so thanks to (\ref{eq: approxtildef}) we can assume here that $\widetilde f_n$ is bounded independently of $n$ from above and from below.
\item In the random censorship model, we still can assume that $\widetilde f_n$ is bounded independently of $n$ from above and from below. Besides, we can assume that $\widetilde G_n(1)<1-\varepsilon$ and $\lim_{t\uparrow 1}H_{nj}(t)<1-\varepsilon$ for some positive number $\varepsilon$ that does not depend on $n$, where
$\widetilde G_n=1-\exp(-\widetilde F_n)$
is the distribution function of the $T_{ij}^\star$'s. Now, the natural bootstrap version of $L_j$ to consider is
\begin{equation}\label{eq:Lboot1}
t\mapsto \int_0^t\frac{\widetilde f_n(x)}{(1-\widetilde G_n(x))(1-H_{nj}(x))}\,\md x.
\end{equation}
But the maximal heigth of the jumps of $H_{nj}$ on $[0,1)$ is of the order $O_p(1/n)$ and therefore,
$$\sup_{t\in[0,1)}|H_{nj}(t)-\widetilde H_j(t)|=O_p(1/n),$$
where
$\widetilde H_j$ is the continuous version of $H_{nj}$ (we mean, the polygonal function on $[0,1]$ that coincides with $H_{nj}$ at every discontinuity points of $H_{nj}$ on $[0,1)$ and such that $\widetilde H_j(1)=\lim_{t\uparrow 1}H_{nj}(t)$). One can then check that the bootstrap version of the embedding also holds with the bootstrap version of $L_j$ defined by
$$L_j ^\star(t)=\int_0^t\frac{\widetilde f_n(x)}{(1-\widetilde G_n(x))(1-\widetilde H_j(x))}\,\md x.$$
The advantage of this proposal as compared to (\ref{eq:Lboot1}) is that it has a continuous derivative that is bounded from below and from above.
\end{enumerate}
Finally, with the above definitions of the bootstrap versions of $f_0$ and $L_j$, the following bootstrap version of
\regul\ clearly holds in both regression and density models: there exist an $s\in(3/4,1]$ that does not depend on $n$ such that for all $x,y,j$,
$$
|\widetilde f_n'(x)-\widetilde f_n'(y)|\leq |x-y|^{s}/\epsilon_n
\quad\text{ and }\quad
|{L_j^\star}''(x)-{L_j^\star}''(y)|\leq |x-y|^{s}/\epsilon_n.
$$
In the random censorship model, $\widetilde f_n'$ satisfies the above H\"olderian assumption, but $L_j^\star$ is not twice differentiable so the bootstrap version of the smoothness assumption on $L_j$ is not satisfied. Instead, we have
$$\sup_{t\in[0,1]}|{L_j^\star}'(t)-L_{j}'(t)|=o_p(n^{-1/3}),$$
where $L_{j}'$ is smooth.

Now that we have bootstrap versions of \indep,\dots,\regul, it can be proved, following the line of reasoning used in the proof of
Theorems~\ref{th:main1} and~\ref{th:main2}, that $n^{1/6} (S_{nk}^\star - m_k)$
converges in distribution to the Gaussian law with mean zero and variance $\sigma_k^2.$
This can be done at the price of additional difficulties which are mainly due to the facts that the bootstrap versions $\widetilde f_n$ and $L_j^\star$ of $f$ and $L_j$ depend on $n$, and $L_j^\star$ is less smooth than the original~$L_j$ in the random censorship model.
In particular, in this model we cannot use a bootstrap version of Lemma~6.3 in the supplement,
which is the key lemma that makes the transition possible from Brownian bridge to Brownian motion.
However, here $q$ can be chosen as large as we wish and $B_{nj}$ is a Brownian motion (see assumption \embed), so
in this model one can avoid the use of a bootstrap version for Lemma~6.3 in the supplement
and directly obtain
$$n^{1/3}\int_0^1|\widehat U_{ni}^{E\star}(t)-\widehat U_{nj}^{E\star}(t)|\,\md t=n^{1/3}\int_{f(1)}^{f(0)}|\widehat U_{ni}^{W\star}(a)-\widehat U_{nj}^{W\star}(a)|\,\md a+o_p(n^{-1/6})$$
(where as usual, $\widehat U_{ni}^{E\star}$ and $\widehat U_{ni}^{W\star}$ are defined in the same manner as $\widehat U_{ni}^{E}$ and $\widehat U_{ni}^{W}$ respectively, just replacing the original observations with their bootstrap version in the definition), by using Lemma 5 in \cite{durot2007}. On the other hand, the presence of $\epsilon_n$ in the bootstrap version of \regul\ does not cause any trouble thanks to the assumption that $n^{-\gamma}/\epsilon_n$ tends to zero as $n\to\infty$ for any positive $\gamma$. Nevertheless, the additional difficulties are not essential so, to alleviate the paper, we do not provide a detailed proof for the consistency of the bootstrap.

It should be mentioned that the proof of Theorems~\ref{th:main1} and~\ref{th:main2} could be simplified if,
instead of considering the general framework of Section \ref{sec:main}, we restrict ourselves to one of models considered in
Section~\ref{sec:calibration}; the proof of the bootstrap version simplifies in the same manner.
For example, in the regression model as well as in the random censoring model, the embedding in \embed~ is with a Brownian motion $B_{nj}$ so the cases $S=B$ and $S=W$ coincide in these models; in particular, Lemmas~\ref{lem:conv slopes}, \ref{lem:mixing} and \ref{lem:transition} are pointless and there is no need to prove a bootstrap version of these lemmas.  On the other hand, in the monotone regression model as well as in the monotone density model,
the processes~$W_{t0}$ and $\widetilde W_{t0}$ coincide since the functions $L_j$'s are proportional to each other, see (\ref{eq:Wt0 in Wtj}) and~(\ref{eq:def tildeWtj}).
The proof thus simplifies a bit, and the same hold with its bootstrap version.
 \tqed

\noindent{\bf Acknowledgements.}
The research of C\'{e}cile Durot is partly supported by the
French Agence Nationale de la Recherche [ANR 2011 BS01 010 01 projet Calibration].


\def\cprime{$'$}


\pagebreak

\centerline{\LARGE\bf Testing equality of functions under monotonicity}
\centerline{\LARGE\bf constraints}
\bigskip
\centerline{\LARGE Supplementary Material}
\bigskip
\centerline{\Large C\'ecile Durot, Piet Groeneboom and Hendrik P.~Lopuha\"{a}}
\medskip
\centerline{\it Universit\'e Paris Ouest Nanterre La D\'efense, Nanterre and Delft University of Technology}
\bigskip
\centerline{\date{\today}}

\setcounter{equation}{60}
\setcounter{section}{5}

\section{Supplement}

\subsection{Proofs of Section~\ref{sec:main}}

\noindent
{\bf Proof of Lemma~\ref{lem:embed0}:}
Because $f_0=f_1=\cdots=f_J$ and $\sum_j c_j=1$, according to~\eqref{eq:def Fn0}, we can write
\[
M_{n0}=F_{n0}-F_0=\sum_{j=1}^J c_j(F_{nj}-F_j)=\sum_{j=1}^J c_jM_{nj}.
\]
Since $c_j=n_j/n$, for $j=1,2,\ldots,J$, this means that
\[
M_{n0}(t)-n^{-1/2}B_{n0}\circ L_0(t)
=\frac{1}{n}
\sum_{j=1}^Jn_j
\left(
M_{nj}(t)-n_j^{-1/2}B_{nj}\circ L_j(t)
\right).
\]
We now have that
\begin{equation*}
\begin{split}
&\prob\left\{
n^{1-1/q}
\sup_{t\in[0,1]}
\left|
M_{n0}(t)-n^{-1/2}B_{n0}\circ L_{0}(t)
\right|
>x
\right\}\\
&\leq
\prob\left\{
n^{-1/q}
\sum_{j=1}^J
n_j
\sup_{t\in[0,1]}
\bigg|
M_{nj}(u)-n_j^{-1/2}B_{nj}\circ L_j(u)
\bigg|>
x
\right\}\\
&\leq
\sum_{j=1}^J
\prob\left\{
n_j^{1-1/q}
\sup_{t\in[0,1]}
\bigg|
M_{nj}(u)-n_j^{-1/2}B_{nj}\circ L_j(u)
\bigg|>
\frac{x}Jc_j^{-1/q}
\right\}.
\end{split}
\end{equation*}
Applying \embed~to each summand yields (\ref{eq:embed0}) for all
$0<x\leq nJ(\min c_j)^{1+1/q}$.
Possibly enlarging $C$, we obtain the result for all $x\in(0,n]$.
\tqed

\bigskip

\noindent
To prove Lemma~\ref{lem:switch} for all $i,j=0,1,\ldots,J$ and $S=E,B,W$, we first need to establish
two results on the tail probabilities of $\widehat{U}_{nj}^S$, that are an extension of Lemmas~3 and~4 in~\cite{durot2007}.
For $j=0,1,\ldots,J$, let $g_j$ denote the inverse of $f_j$, defined for $a\in\R$ by
\begin{equation*}
g_j(a)=\sup\{t\in[0,1],\ f_j(t)\geq a\},
\end{equation*}
with the convention that the supremum of an empty set is 0.
\begin{lemma}
\label{lem:tail 01}
Assume \mono, \mart\ and suppose $f_0=f_1=\dots=f_J$.
Then, there exists~$C>0$, such that for all $j=0,1,\ldots,J$, $S=E,B,W$, $x>0$, and $a\notin [f_j(1),f_j(0)]$,
\[
\prob
\left(
|\widehat U_{nj}^S(a)-g_j(a)|> x
\right)
\leq
\frac{C}{nx(f_j\circ g_j(a)-a)^2}.
\]
\end{lemma}
{\bf Proof:}
We follow the line of reasoning used in the proof of Lemma 3 in \cite{durot2007}.
It suffices to prove the result for $x\in(0,1]$, since the considered probability vanishes for $x\geq 1$.
Moreover, for simplicity, we restrict ourselves to the case $a>f_j(0)$, so that $g_j(a)=0$.
The case $a<f_j(1)$ can be treated likewise.
By definition of $\widehat U_{nj}^S$ we have:
\begin{equation}
\label{lem:sup}
\prob\left(|\widehat U_{nj}^S(a)-g_j(a)|>x\right)
\leq
\prob\left(\sup_{|u-g_j(a)|>x}
\left\{F_{nj}^S(u)-au\right\}\geq F_{nj}^S(g_j(a))-ag_j(a)\right).
\end{equation}
Since $f_j$ is decreasing according to~\mono, we have $F_j(u)-F_j(0)\leq uf_j(0)$, for all $u\in[0,1]$.
Therefore,
\[
\prob\left(|\widehat U_{nj}^S(a)-g_j(a)|>x\right)
\leq
\prob\left(\sup_{u>x}
\left\{M_{nj}^S(u)- M_{nj}^S(0)-u(a-f_j(0))\right\}
\geq0\right),
\]
where
\begin{equation}
\label{eq:def M_nj^S}
M_{nj}^S(u)=F_{nj}^S(u)-F_j(u).
\end{equation}
It then follows from Markov's inequality that
\[\begin{split}
\prob\left(|\widehat U_{nj}^S(a)-g_j(a)|>x\right)
&\leq
\sum_{k\geq 1}
\prob\left(\sup_{u\in[x2^{k-1},x2^k]}
\left\{M_{nj}^S(u)- M_{nj}^S(0)\right\}\geq x2^{k-1}(a-f_j(0))\right)\\
&\leq
\sum_{k\geq 1}
\dfrac{\E\left[\sup_{u\in[x2^{k-1},x2^k]}\left(M_{nj}^S(u)- M_{nj}^S(0)\right)^2\right]}{x^22^{2k-2}(a-f_j(0))^2}.\end{split}\]
If we are able to show that there exists a constant $C>0$, such that for all $x\geq 0$,
\begin{equation}
\label{eq:bound E}
\mathds{E}
\left[
\sup_{x/2\leq u\leq x}
\left(M_{nj}^S(u)-M_{nj}^S(0)\right)^2
\right]
\leq
\frac{Cx}{n},
\end{equation}
then we would obtain the required result since in that case,
$$
\prob\left(|\widehat U_{nj}^S(a)-g_j(a)|>x\right)
\leq
\sum_{k\geq 1}\dfrac{Cx2^k}{nx^22^{2k-2}(a-f_j(0))^2}\leq \dfrac{K}{nx(a-f_j(0))^2}.
$$
It thus remains to prove (\ref{eq:bound E}) for $S=E,B,W$.
First consider $j=1,2,\dots,J$.
In the case~$S=E$, inequality~(\ref{eq:bound E}) follows from condition~\mart.
In the case~$S=W$,
since $L_j$ is increasing, the process
$M_{nj}^W(u)=n_j^{-1/2}W_{nj}(L_j(u))$
is a mean zero martingale
and it follows from Doob's inequality, that
for all $x\geq 0$:
\[
\mathbb{E}
\left[
\sup_{0\leq u\leq 1\wedge x}
\left(M_{nj}^W(u)-M_{nj}^W(0)\right)^2
\right]
\leq
4
\mathbb{E}
\left[
\left(M_{nj}^W(1\wedge x)-M_{nj}^W(0)\right)^2
\right]
\leq
\frac{Cx}{n},
\]
because $M_{nj}^W(1\wedge x)-M_{nj}^W(0)$ has a normal distribution with mean zero and variance
$n_j^{-1}(L_j(1\wedge x)-L_j(0))$, which is bounded by $n_j^{-1}(\sup_{t} |L_j'(t)|)x$.
This yields~(\ref{eq:bound E}) for $S=W$.
In the case~$S=B$, we write
\[
M_{nj}^B(t)-M_{nj}^B(u)
=
M_{nj}^W(t)-M_{nj}^W(u)
-n_j^{-1/2}\xi_{nj}(L_j(t)-L_j(u)).
\]
Next, inequality~\eqref{eq:jensen} yields:
\[
\begin{split}
\mathbb{E}
\left[
\sup_{0\leq u\leq 1\wedge x}
\left(M_{nj}^B(u)-M_{nj}^B(0)\right)^2
\right]
&\leq
\frac{2Cx}{n}
+
\frac{2\mathbb{E}[\xi_{nj}^2]}{n_j}
\sup_{0\leq u\leq 1\wedge x}
\left(L_j(u)-L_j(0)\right)^2\\
&\leq
\frac{2Cx}{n}
+
\frac{2}{n_j}
\left(
\sup_{t\in[0,1]} |L_j'(t)|
\right)^2(1\wedge x)^2.
\end{split}
\]
For $x\geq 0$, we have$ (1\wedge x)^2\leq x$,
whence (\ref{eq:bound E}) holds for $S=B$.

Finally, consider the case $j=0$ and $S=E,B,W$. From \eqref{eq:def F_n0^S}, \eqref{eq:def M_nj^S}, and \eqref{eq:jensen},
we find:
\begin{equation}
\mathds{E}
\left[
\sup_{x/2\leq u\leq x}
\left(
M_{n0}^S(u)-M_{n0}^S(0)
\right)^2
\right]\leq J
\sum_{j=1}^J
c_j^2
\mathds{E}
\left[
\sup_{x/2\leq u\leq x}
\left(
M_{nj}^S(u)-M_{nj}^S(0)
\right)^2
\right],
\end{equation}
for all $x\geq 0$. Since (\ref{eq:bound E}) holds for $j=1,2,\dots,J$ and $S=E,B,W$,
it also holds for $j=0$ and~$S=E,B,W$.
\tqed
\begin{lemma}
\label{lem:tail}
Assume \mono, \embed\ and suppose $f_0=f_1=\dots=f_J$.
Then, there exist $c>0$ and $C>0$ such that for every $j=0,1,\ldots,J$,  $x>0$ and $a\in\R$,
\begin{equation}
\label{eq:tailBW}
\prob\left(|\widehat U_{nj}^S(a)-g_j(a)|>x\right)\leq 2\exp(-cnx^3),
\end{equation}
and
\begin{equation*}
\prob\left(|\widehat U_{nj}^S(a)-g_j(a)|>x\right)\leq C(nx^3)^{1-q},
\end{equation*}
for $S=E,B,W$.
\end{lemma}
{\bf Proof:}
We follow the line of reasoning used in the proof of Lemma 4 in \cite{durot2007} and assume $x\in(0,1]$.
Let $\beta$ satisfy $0<\beta<\inf_t|f_j'(t)|/2$, for every $j=1,2,\ldots,J$,
which is possible according to~\mono.
From Taylor's expansion,
$F_j(u)-F_j(g_j(a))\leq(u-g_j(a))a-\beta(u-g_j(a))^2$,
for all $u\in[0,1]$, so (\ref{lem:sup}) yields:
\begin{equation}
\label{eq:tail}
\begin{split}
&
\prob\left(|\widehat U_{nj}^S(a)-g_j(a)|>x\right)\\
&\quad\leq
\prob\left(\sup_{|u-g_j(a)|>x}
\left\{M_{nj}^S(u)-M_{nj}^S(g_j(a))-\beta(u-g_j(a))^2\right\}\geq 0\right).
\end{split}
\end{equation}
To cover all cases $j=0,1,\ldots,J$ simultaneously, write
\begin{equation*}
n_0=n
\quad\text{and define}\quad
\xi_{nj}^S(t)
=
\begin{cases}
\xi_{n0}(t) & ,j=0\\
\xi_{nj}t  & ,j=1,2,\ldots,J.
\end{cases}
\end{equation*}
where $\xi_{n0}$ is defined in~\eqref{eq:def xi0}.
Note that $\xi_{nj}^W(t)\equiv0$, for all $j=0,1,\ldots,n$.
For $S=B,W$, we have $M_{nj}^S(u)=M_{nj}^W(u)-n_j^{-1/2}\xi_{nj}^S(L_j(u))$,
so that in both cases,
\begin{equation}
\label{eq:Pbound U}
\begin{split}
&\prob\left(|\widehat U_{nj}^S(a)-g_j(a)|>x\right)\\
&\quad\leq
\prob\left(\sup_{|u-g_j(a)|>x}
\left\{
M_{nj}^W(u)-M_{nj}^W(g_j(a))-\frac{\beta}{2}(u-g_j(a))^2
\right\}\geq 0\right)\\
&\qquad+
\prob\left(\sup_{|u-g_j(a)|>x}
\left\{
\xi_{nj}^S(L_j(g_j(a)))-\xi_{nj}^S(L_j(u))-\frac{\beta n_j^{1/2}}{2}(u-g_j(a))^2
\right\}\geq 0\right).
\end{split}
\end{equation}
Let $W$ be a standard two-sided Brownian motion on $\R$.
Setting $v=L_j(u)-L_j(g_j(a))$, $k_0\leq \inf_t|L_j'(t)|$ and $k_1\geq\sup_t|L_j'(t)|^2$,
one can derive from the assumptions on $L_j$ and scaling properties of $W$ that
\[
\begin{split}
&\prob\left(\sup_{|u-g_j(a)|>x}
\left\{
M_{nj}^W(u)-M_{nj}^W(g_j(a))-\frac{\beta}{2}(u-g_j(a))^2
\right\}\geq 0\right)\\
&\quad\leq
\prob\left(\sup_{|v|>k_0x}
\left\{n_j^{-1/2}W(v)-\frac{\beta v^2}{2k_1}\right\}\geq0\right)
\leq
\prob\left(\sup_{|v|>k_0x}\left\{\frac{W(v)}{|v|}\right\}\geq n_j^{1/2}\frac{\beta k_0 x}{2k_1}\right).
\end{split}
\]
The process $\{uW(1/u),\ u>0\}$ is distributed like $\{W(u),\ u>0\}$, and the distribution of~$W$ is symmetric about zero.
By Proposition~1.8 on page~52 in~\cite{revuz-yor1991},
we conclude that there exists~$C,K_1>0$ depending only on $L_j$, such that
\begin{equation}
\label{eq:bound tailprob BM}
\begin{split}
&\prob\left(\sup_{|u-g_j(a)|>x}
\left\{M_{nj}^W(u)-M_{nj}^W(g_j(a))-\frac{\beta}{2}(u-g_j(a))^2\right\}\geq 0\right)\\
&\quad\leq
2\prob\left(\sup_{v>k_0x}\left\{\frac{W(v)}{v}\right\}\geq n_j^{1/2}\frac{\beta k_0 x}{2k_1}\right)
\leq 2\exp(-K_1n\beta^2 x^3).
\end{split}
\end{equation}
From (\ref{eq:Pbound U}), there thus exists $K_2>0$, such that
$$
\prob\left(|\widehat U_{nj}^S(a)-g_j(a)|>x\right)\leq 2\exp\left(-K_1n\beta^2x^3\right)+C(j)\prob(|\xi|>K_2\beta n^{1/2}x),
$$
where $\xi$ is a standard Gaussian variable, $C(j)=1$ for $j=1,2,\dots,J$ and $C(0)=J$.
Since $x\in(0,1]$, the latter probability is bounded from above by
$\exp\left(-K_1n\beta^2x^3\right)$,
provided $K_1$ is small enough. Therefore, there exists $C>0$ with
$$\prob\left(|\widehat U_{nj}^S(a)-g_j(a)|>x\right)\leq 3\exp\left(-K_1n\beta^2x^3\right)\leq C(nx^3)^{1-q},$$
which proves the lemma for the cases $S=B,W$.  In the case $S=E$, (\ref{eq:tail}) yields
\[
\begin{split}
\prob\left(|\widehat U_{nj}^E(a)-g_j(a)|>x\right)
&\leq
\prob\left(\sup_{|u-g_j(a)|>x}
\left\{
M_{nj}^B(u)-M_{nj}^B(g_j(a))-\frac{\beta}{2}(u-g_j(a))^2
\right\}\geq 0\right)\\
&\quad+
\prob\left(\sup_{u\in[0,1]}|M_{nj}^E(u)-M_{nj}^B(u)|\geq \frac{\beta x^2}{4}\right),
\end{split}
\]
where we recall that $M_{nj}^B(u)=n_j^{-1/2}B_{nj}(L_j(u))$.
From what precedes, the first probability on the right hand side is bounded from above by $ C(nx^3)^{1-q}$ for some $C>0$.
Note that $x\in(0,1]$ and we can assume $\beta<1$. In the cases $j=1,2,\dots,J$, it follows from assumption \embed\
that there is $K_4>0$ with
$$
\prob\left(\sup_{u\in[0,1]}|M_{nj}^E(u)-M_{nj}^B(u)|\geq \frac{\beta x^2}{4}\right)
\leq
K_4\beta^{-q}x^{-2q}n^{1-q}
\leq
K_4\beta^{-q}(nx^3)^{1-q}.
$$
From Lemma \ref{lem:embed0}, this inequality still holds with $j=0$, which completes the proof.
\tqed

\bigskip

\noindent
{\bf Proof of Lemma~\ref{lem:switch}:}
We have
$$
\int_0^1(\widehat f_{ni}^S(t)-\widehat f_{nj}^S(t))_+\,\md t
=
\int_0^1\int_0^\infty\indicator_{a\leq\widehat f_{ni}^S(t)-\widehat f_{nj}^S(t)}\,\md a\,\md t.
$$
For all $t\in(0,1]$ and $a\in\R$, $\widehat U_{ni}^S(a)\geq t$ if and only if $\widehat f_{ni}^S(t)\geq a$,
so the change of variable $b=a+\widehat f_{nj}^S(t)$ yields
\[
\begin{split}
\int_0^1(\widehat f_{ni}^S(t)-\widehat f_{nj}^S(t))_+\,\md t
&=
\int_\R\int_0^1\indicator_{\widehat f_{nj}^S(t)<b\leq\widehat f_{ni}^S(t)}\,\md t\,\md b\\
&=
\int_\R\int_0^1\indicator_{\widehat U_{nj}^S(b)<t\leq\widehat U_{ni}^S(b)}\,\md t\,\md b\\
&=
\int_\R(\widehat U_{ni}^S(b)-\widehat U_{nj}^S(b))_+\,\md b,
\end{split}
\]
whence
$$
\int_0^1|\widehat f_{nj}^S(t)-\widehat f_{ni}^S(t)|\,\md t=\int_\R|\widehat U_{ni}^S(b)-\widehat U_{nj}^S(b)|\,\md b.
$$
From Fubini's theorem and Lemmas~\ref{lem:tail 01} and~\ref{lem:tail}, there exists $K>0$,
such that for all~$j=0,1,\ldots,J$ and~$S=E,B,W$,
\[
\begin{split}
\E\int_{f_j(0)}^\infty|\widehat U_{nj}^S(b)-g_j(b)|\,\md b
&=
\int_0^\infty\int_{f_j(0)}^\infty\prob(|\widehat U_{nj}^S(b)-g_j(b)|>x)\,\md b\,\md x\\
&\leq
K\int_0^\infty\int_0^\infty1\wedge\dfrac{1}{nx^3}\wedge\dfrac{1}{nxb^2}\,\md b\,\md x\\
&\leq
K\int_0^{n^{-1/3}}\left( \dfrac{1}{\sqrt{nx}}+\int_{1/\sqrt{nx}}^\infty\dfrac{1}{nxb^2}\,\md b\right) \,\md x\\
&\quad+
K\int_{n^{-1/3}}^\infty\left( \dfrac{1}{nx^2}+\int_x^\infty\dfrac{1}{nxb^2}\,\md b\right) \,\md x\\
&=
6Kn^{-2/3}.
\end{split}
\]
The integral $\int_{-\infty}^{f_j(1)}|\widehat U_{nj}^S(b)-g_j(b)|\,\md b$ can be treated likewise,
so we obtain
$$
\int_0^1|\widehat f_{nj}^S(t)-\widehat f_{ni}^S(t)|\,\md t=\int_{f(1)}^{f(0)}|\widehat U_{ni}^S(b)-\widehat U_{nj}^S(b)|\,\md b+O_\prob(n^{-2/3}).
$$
\tqed

\bigskip

\noindent
In order to prove Lemma~\ref{lem:E to B}, the line of reasoning is similar to that in~\cite{groeneboom-hooghiemstra-lopuhaa1999},
except that in the current situation we work in a more general model.
The key result that we need is an extended version of Lemma~3.4 in~\cite{groeneboom-hooghiemstra-lopuhaa1999},
provided by the following lemma.
\begin{lemma}
\label{lem:jump}
Assume \mono, \embed, \regul\ and suppose $f_0=f_1=\dots=f_J$.
For $j=0,1,\ldots,J$ and~$S=E,B,W$ let
\begin{equation}
\label{eq:def Vn^S}
\widehat{V}_{nj}^S(a)=n_j^{1/3}\left(\widehat{U}_{nj}^S(a)-g_j(a)\right),
\end{equation}
and
\[
J_{nj}
=
\left[
f_j(1)+n_j^{-1/3}\log n_j, f_j(0)-n_j^{-1/3}\log n_j
\right],
\]
where $n_0=n$.
Then there exists a constant $\beta_j>0$, independent of $a\in J_{nj}$, such that
for~$S=B,W$ and for all $h\in(0,1)$,
\[
\prob
\left\{
\widehat{V}_{nj}^S\text{ jumps in }\left(a-hn_j^{-1/3},a+hn_j^{-1/3}\right)
\right\}
\leq
\beta_j\delta_{n,h}+o(\delta_{n,h})
\]
as $\delta_{n,h}\downarrow0$, where $\delta_{n,h}=h\vee(n^{-\holder/3}(\log n)^{\holder+1})$,
where $\holder\in(3/4,1]$ is taken from~\regul.
\end{lemma}
\textbf{Proof:}
We follow the line of reasoning used in the proof of Lemma 3.4 in~\cite{groeneboom-hooghiemstra-lopuhaa1999}.
In order to transform the process $t\mapsto W_{nj}\circ L_j(g_j(a)+n_j^{-1/3}t)$ into the process
$t\mapsto W_{nj}(L_j(g_j(a))+n_j^{-1/3}t)$,
we define for $|c|\leq 1$,
\[
\widehat{L}_j^S(a,c)
=
n_j^{1/3}
\left\{
L_j\big(\widehat{U}_{nj}^S(a+n_j^{-1/3}c)\big)-L_j\big(g_j(a)\big)
\right\}.
\]
Then $\widehat{V}_{nj}^S$ has a jump in $(a-hn_j^{-1/3},a+hn_j^{-1/3})$
if and only if $c\mapsto\widehat{L}_j^S(a,c)$ has a jump in~$(-h,h)$.
Then in the case $S=W$ we have
\[
\widehat{L}_j^W(a,c)=
\argmax_{y\in I_{nj}(a)}
\left\{
W_{g_j(a),j}(y)-p_{nj}(c,y)
\right\},
\]
where
\begin{equation}
\label{eq:def Inj Wtj supplement}
\begin{split}
I_{nj}(a)
&=
\left[
n_j^{1/3}(L_j(0)-L_j(g_j(a))),\, n_j^{1/3}(L_j(1)-L_j(g_j(a)))
\right],\\
W_{tj}(y)
&=
n_j^{1/6}
\left\{
W_{nj}\big(L_j(t)+n_j^{-1/3}y)-W_{nj}(L_j(t)\big)
\right\},
\end{split}
\end{equation}
with $W_{nj}$, for $j=1,2,\ldots,J$, being independent Brownian motions from~\eqref{eq:Bnj and Wnj}
and $W_{n0}$ is the Brownian motion defined by (\ref{eq:def Wn0}),
and where
\[
\begin{split}
p_{nj}(c,y)
&=
-n_j^{2/3}\left\{
\left(F_j\circ H_j\right)\left(L_j(g_j(a))+n_j^{-1/3}y\right)
-
F_j\left(g_j(a)\right)
\right\}\\
&\qquad
+
n_j^{2/3}(a+n_j^{-1/3}c)
\left\{
H_j\left(L_j(g_j(a))+n_j^{-1/3}y\right)
-
g_j(a)
\right\},
\end{split}
\]
with $H_j=L_j^{-1}$.
To deal with the case $S=B$ for all $j=0,1,2,\ldots,J$ simultaneously,
write $B_{nj}(t)=W_{nj}(t)-\xi_{nj}^B(t)$, where $\xi_{nj}^B$ is defined in~\eqref{eq:def xi_nj^S},
and write $t_j=g_j(a)$.
Then,
\[
\widehat{L}_j^B(a,c)=
\argmax_{y\in I_{nj}(a)}
\left\{
B_{t_j,j}(y)-p_{nj}(c,y)
\right\},
\]
where
\[
B_{t_j,j}(y)
=W_{t_j,j}(y)-n_j^{1/6}
\left\{
\xi_{nj}^B\big(L_j(t_j)+n_j^{-1/3}y\big)-\xi_{nj}^B\big(L_j(t_j)\big)
\right\}.
\]
Now, to deal simultaneously with the cases $S=B,W$,  define $\widehat{\psi}_{nj}^W(c)=\widehat{L}_j^W(a,c)$ and
\[
\widehat{\psi}_{nj}^B(c)
=
\argmax_{y\in I_{nj}(a)}
\left\{
B_{t_j,j}(y)-p_{nj}\big(c-n_j^{-1/6}L_j'(t_j){\xi_{nj}^B}'(L_j(t_j)),y\big)
\right\},
\]
so that $\widehat{L}_j^S(a,c)=\widehat{\psi}_{nj}^S(c+n_j^{-1/6}L_j'(t_j){\xi^S_{nj}}'(L_j(t_j))$, for $S=B,W$.
It follows that
\[
\widehat{\psi}_{nj}^S(c)
=
\argmax_{y\in I_{nj}(a)}
\left\{
W_{t_j,j}(y)
-
q_{nj}^S(c,y)
\right\},
\]
with
\[
\begin{split}
q_{nj}^S(c,y)
&=
n_j^{1/6}
\left\{
\xi_{nj}^S\big(L_j(t_j)+n_j^{-1/3}y\big)-\xi_{nj}^S\big(L_j(t_j)\big)
\right\}\\
&\qquad+
p_{nj}\big(c-n_j^{-1/6}L_j'(t_j){\xi_{nj}^S}'(L_j(t_j)),y\big).
\end{split}
\]
Suppose that the process $\widehat{V}_{nj}^S$ jumps in the interval $(a-hn_j^{-1/3},a+hn_j^{-1/3})$.
Then $\widehat{\psi}_{nj}^S$ has a jump at some
\[
c^*\in
\left(
-h+n_j^{-1/6}L_j'(t_j){\xi_{nj}^S}'(L_j(t_j)),h+n_j^{-1/6}L_j'(t_j){\xi_{nj}^S}'(L_j(t_j))
\right).
\]
This means that if we drop the function $y\mapsto q_{nj}^S(c^*,y)+\beta$ for varying $\beta\in\R$ onto
the process~$W_{t_j,j}$, it first touches $W_{t_j,j}$ simultaneously in two points,
say $(y_i^S,w_i^S)$, for $i=1,2$, where
$w_i^S=W_{t_j,j}(y_i^S)=q_{nj}^S(c^*,y_i^S)+\beta^S$,
for some $\beta^S\in\R$.
Define the event $A_{nj}^S=\{|\widehat{L}_j^S(a,c)|\leq \log n_j\text{, for all }|c|\leq 1\}$.
By Lemma~\ref{lem:tail} it follows that
$\prob((A_{nj}^S)^c)=o(\delta_{n,h})$.
Furthermore, define
$A_{n}'=
\{|\xi_{nj}|\leq n_j^{1/6},\text{ for all }j=1,2,\ldots,J\}$.
Then, $\prob((A_{n}')^c)=o(\delta_{n,h})$.
Hence,
\[
\begin{split}
&
\prob\left(\widehat{V}_{nj}^S\text{ jumps in }(a-hn_j^{1/3},a+hn_j^{1/3})\right)\\
&\leq
\prob\left(\widehat{V}_{nj}^S\text{ jumps in }(a-hn_j^{1/3},a+hn_j^{1/3}),A_{nj}^S\cap A_n'
\right)+o(\delta_{n,h}),
\end{split}
\]
and we can restrict ourselves to the event $A_{nj}^S\cap A_{n}'$. On this event, we have $|y_1^S-y_2^S|\leq 2\log n_j$.

Next, we show that for each $y_i^S$, $i=1,2$, we can construct
a parabola that lies above $q_{nj}^S(c^*,y)+\beta^S$ for all $|y|\leq\log n_j$ and that touches $W_{t_j,j}(y)$
at $(y_i^S,w_i^S)$. To this end first consider the derivatives of $p_{nj}(c,y)$:
\[
\begin{split}
\frac{\,\md p_{nj}(c,y)}{\,\md y}
&=
-n_j^{1/3}\big(F_j\circ H_j\big)'\big(L_j(t_j)+n_j^{-1/3}y\big)
+
n_j^{1/3}(a+n_j^{-1/3}c)
H_j'\big(L_j(t_j)+n_j^{-1/3}y\big),\\
\frac{\,\md^2 p_{nj}(c,y)}{\,\md y^2}
&=
-\big(F_j\circ H_j\big)''\big(L_j(t_j)+n_j^{-1/3}y\big)
+
(a+n_j^{-1/3}c)H_j''\big(L_j(t_j)+n_j^{-1/3}y\big).
\end{split}
\]
With assumption \regul, we find that
\[\begin{split}
\frac{\,\md^2 p_{nj}(c,y)}{\,\md y^2}&=
\frac{\,\md^2 p_{nj}(c,0)}{\,\md y^2}+(1+cn^{-1/3})O((n^{-1/3}y)^\holder)\\
&=\frac{|f_j'(t_j)|}{L_j'(t_j)^2}+(1+cn^{-1/3})O((n^{-1/3}y)^\holder)+O(cn^{-1/3}).\end{split}
\]
The second derivative of $q_{nj}^S(c,y)$ is
\[\frac{\,\md^2 q_{nj}^S(c,y)}{\,\md y^2}=
n_j^{-1/2}{\xi_{nj}^S}''\big(L_j(t_j)+n_j^{-1/3}y\big)
+
\frac{\,\md^2 p_{nj}\big(c-n_j^{-1/6}L_j'(t_j){\xi_{nj}^S}'(L_j(t_j)),y\big)}{\,\md y^2}.
\]
Using that
$p_{nj}(c+d,y)
=
p_{nj}(c,y)
+dn_j^{1/3}
\{H_j(L_j(t_j)+n_j^{-1/3}y)
-
t_j\}$,
and the fact that on the event $A_{nj}^S\cap A_{n}'$, one has
$n_j^{-1/2}|{\xi_{nj}^S}''\big(L_j(t_j)+n_j^{-1/3}y\big)|
\leq
K_1n^{-1/3}$,
for some $K_1>0$ only depending on the $c_j$'s and $L_j$'s, we conclude that on $A_{nj}^S\cap A_{n}'$,
$$\frac{\,\md^2 q_{nj}^S(c,y)}{\,\md y^2}=\frac{|f_j'(t_j)|}{L_j'(t_j)^2}\left(1+O(n^{-1/3}\log n)^\holder\right),$$
where the approximation is uniform in $|y|\leq\log n_j$ and $|c|\leq 1.$
Therefore, there exists a $K_2>0$ that only depends on $f_j$, the $c_j$'s and $L_j$'s,
such that on $A_{nj}^S\cap A_{n}'$,
\[
0<\frac{\,\md^2 q_{nj}^S(c,y)}{\,\md y^2}
\leq
\frac{|f_j'(t_j)|}{L_j'(t_j)^2}
\left\{
1+K_2(n^{-1/3}\log n)^\holder
\right\},
\]
for all $|y|\leq\log n_j$ and $|c|\leq 1.$ Choose $M>K_2$ and define the parabola
\begin{equation*}
\pi_{nj}(c,y)
=
\frac{c}{L_j'(t_j)}y+\alpha_{nj}y^2,
\end{equation*}
where
\begin{equation}
\label{eq:def alpha_nj}
\alpha_{nj}=\frac{|f_j'(t_j)|}{2L_j'(t_j)^2}
\left\{
1+M(n^{-1/3}\log n)^\holder
\right\}.
\end{equation}
Then it follows that for all $|y|\leq\log n_j$, $|c|\leq 1$, and $b\in\R$:
\[
\frac{\,\md^2 \pi_{nj}(b,y)}{\,\md y^2}
>
\frac{\,\md^2 q_{nj}^S(c,y)}{\,\md y^2}.
\]
Now, for each $i=1,2$, choose $b_i$ such that
\[
\frac{b_i}{L_j'(t_j)}+2\alpha_{nj}y_i=\frac{\,\md q^S_{nj}(c^*,y_i)}{\,\md y},
\]
so that the functions $\pi_{nj}(b_i,y)$ and $q^S_{nj}(c^*,y)$ have the same tangent at $y_i$.
If we also take $\gamma_i=q^S_{ nj}(c^*,y_i)-\pi_{nj}(b_i,y_i)$,
then it follows that the parabola $\pi_{nj}(b_i,y)+\gamma_i$ lies above
$q^S_{nj}(c^*,y)$ and touches $q_{nj}^S(c^*,y)$ at $y_i$ in such a way that they have the same tangent.
This implies that if we drop $\pi_{nj}(b_i,y)+\gamma$,
for varying  $\gamma\in\R$ onto the process $W_{t_j,j}$,
it first touches $W_{t_j,j}$ at $y_i$.
Hence, if we define
\[
\widehat{V}_{nj}^{\pi}(c)
=
\argmax_{y\in I_{nj}(a)}
\left\{
W_{t_j,j}(y)-\pi_{nj}(c,y)
\right\},
\]
then from the above construction, it follows that the process $\widehat{V}_{nj}^{\pi}$ has a jump in the interval
$[b_1,b_2]$ of maximal size $|y_1-y_2|\leq 2\log n_j$.
Because

\begin{equation*}
\begin{split}
\frac{\,\md q_{nj}^S(c,y)}{\,\md y}&=
n_j^{-1/6}{\xi_{nj}^S}'\big(L_j(t_j)+n_j^{-1/3}y\big)
+
\frac{\,\md p_{nj}\big(c-n_j^{-1/6}L_j'(t_j){\xi_{nj}^S}'(L_j(t_j)),y\big)}{\,\md y}\\
&=
\frac{c}{L_j'(t_j)}
+
\frac{|f_j'(t_j)|}{L_j'(t_j)^2}y
+
\,O(n^{-1/3}y)
+
O((n^{-1/3}y)^\holder),
\end{split}
\end{equation*}
and
\[
\frac{\,\md \pi_{nj}(b,y)}{\,\md y}
=
\frac{b}{L_j'(t_j)}+2\alpha_{nj}y
=
\frac{b}{L_j'(t_j)}
+
\frac{|f_j'(t_j)|}{L_j'(t_j)^2}y
+y\,O((n^{-1/3}\log n)^\holder),
\]
it follows from $\,\md q^S_{nj}(c^*,y_i)/\,\md y=\,\md \pi_{nj}(b_i,y_i)/\,\md y$, that
there exists a constant $K_3>0$, such that $|b_i-c^*|\leq K_3|y_i+1|(n^{-1/3}\log n)^\holder$,
for $i=1,2$.
Therefore,
\[
|b_i-c^*|\leq K_4n^{-\holder/3}(\log n)^{\holder+1}.
\]
Because
\[
c^*\in (-h+n_j^{-1/6}L_j'(t_j){\xi_{nj}^S}'(L_j(t_j)),h+n_j^{-1/6}L_j'(t_j){\xi_{nj}^S}'(L_j(t_j))),
\]
this means that
\[
\begin{split}
[b_1,b_2]
\subset
\bigg(
-h&+n_j^{-1/6}L_j'(t_j){\xi_{nj}^S}'(L_j(t_j))-K_4n_j^{-\holder/3}(\log n_j)^{\holder+1},\\
&h+n_j^{-1/6}L_j'(t_j){\xi_{nj}^S}'(L_j(t_j))+K_4n_j^{-\holder/3}(\log n_j)^{\holder+1}
\bigg).
\end{split}
\]
On the event $A_{nj}^S\cap A_{n}'$, it follows that
there exists a constant $K_5>0$ such that $[b_1,b_2]$ is contained in
\[
\mathcal{I}_{nj}'
=
\left(
n_j^{-1/6}L_j'(t_j){\xi_{nj}^S}'(L_j(t_j))-K_5\delta_{n,h},
n_j^{-1/6}L_j'(t_j){\xi_{nj}^S}'(L_j(t_j))+K_5\delta_{n,h}
\right),
\]
where $\delta_{n,h}=h\vee(n^{-\holder/3}(\log n)^{\holder+1})$,
so that the process $\widehat{V}_{nj}^\pi$ jumps in the interval $\mathcal{I}_{nj}'$. Hence
\[\begin{split}
\prob\left(\widehat{V}_{nj}^S\text{ jumps in }(a-hn_j^{1/3},a+hn_j^{1/3})\right)
\leq
\prob\left(\widehat{V}_{nj}^\pi\text{ jumps in }\mathcal{I}_{nj}'
\right)+o(\delta_{n,h}).
  \end{split}\]
In the case $S=W$, we have
\[\begin{split}
\prob\left\{
\widehat{V}_{nj}^\pi\text{ jumps in }\mathcal{I}_{nj}'
\right\}=
\prob\left(\widehat{V}_{nj}^\pi\text{ jumps in }\left(K_5\delta_{n,h},K_5\delta_{n,h}\right)
\right).
  \end{split}\]
In the case $S=B$, because $\xi_{nj}^S$ is independent of the process $W_{t_j,j}$,
\[
\begin{split}
\prob
\left\{
\widehat{V}_{nj}^\pi\text{ jumps in }\mathcal{I}_{nj}'
\right\}
=
\int_{-\infty}^{\infty}
\prob
\left\{
\widehat{V}_{nj}^\pi\text{ jumps in }
\left(a-K_5\delta_{n,h},a+K_5\delta_{n,h}\right)
\right\}
h_{nj}(a)\,\md a,
\end{split}\]
where $h_{nj}$ denotes the density of $n_j^{-1/6}L_j'(t_j){\xi_{nj}^S}'(L_j(t_j))$.
Moreover, in both cases $S=W,B$, the process $W_{t_j,j}$ is distributed like Brownian motion $W$, so
the variable $\widehat{V}_{nj}^\pi(c)$ is distributed as
\[
\argmax_{y\in I_{nj}(a)}
\left\{
W(y)-\frac{c}{L_j'(t_j)}y-\alpha_{nj}y^2
\right\},
\]
where $\alpha_{nj}$ and $I_{nj}(a)$ are defined in~\eqref{eq:def alpha_nj} and~\eqref{eq:def Inj Wtj supplement}.
Since $a\in J_ {nj}$ and the $L_j$'s have a bounded first derivative, there exists $K_6$ such that this random variable only differs from
\[
\widetilde{V}_{nj}(c)
=
\argmax_{y\in \R}
\left\{
W(y)-\alpha_{nj}\left(y+\frac{c}{2\alpha_{nj}L_j'(t_j)}\right)^2
\right\},
\]
if $|\widetilde{V}_{nj}(c)|>K_6\log n_j$.
Hence
\[
\begin{split}
&
\prob\left\{
\widehat{V}^\pi_{nj}\text{ jumps in }\left(a-K_5\delta_{n,h},a+K_5\delta_{n,h}\right)\right\}\\
&\leq
\prob\left\{
\widetilde{V}_{nj}\text{ jumps in }\left(a-K_5\delta_{n,h},a+K_5\delta_{n,h}\right)
\right\}+
\prob\left\{
\sup_{c\in \mathcal{I}_{nj}}|\widetilde{V}_{nj}(c)|>\log n_j
\right\}.
\end{split}
\]
According to Lemma 3.3(iii) in \cite{groeneboom-hooghiemstra-lopuhaa1999} the first probability is of the order
$\delta_{n,h}$, uniformly in $a$.
Finally, from the monotonicity of $\widetilde{V}_{nj}(c)$
together with property (3.16) and Lemma~3.3(ii) in~\cite{groeneboom-hooghiemstra-lopuhaa1999},
it follows that the second probability is of smaller order.
This proves the lemma.\tqed

\bigskip

Lemma~\ref{lem:jump} enables us the make the transition from the empirical inverse process
to the Brownian bridge process and establish Lemma~\ref{lem:E to B}.
In the sequel, for $x,t\in\R$ and a function $G:\R\to\R$, we use the notation
$G(t,x]=G(x)-G(t)$.

\bigskip

\noindent
\textbf{Proof of Lemma~\ref{lem:E to B}:}
The proof is along the lines of the proof of Corollary 3.1 in \cite{groeneboom-hooghiemstra-lopuhaa1999}.
For $j=0,1,\ldots,J$, let
\[
K_{nj}=
\left\{
n_j^{2/3}\sup_{t\in[0,1]}
\left|M_{nj}(t)-n_j^{-1/2}B_{nj}(L_j(t))\right|\leq \delta_{nj}
\right\},
\]
where $\delta_{nj}=n_j^{-1/3+1/q}\log n$ and $n_0=n$.
By condition~\embed\ and Lemma~\ref{lem:embed0} it follows that $\prob(K_{nj})\to 1$, as $n_j\to\infty$.
Hence, we can restrict ourselves to the event $K_{nj}$.
Also fix $a\in(f_j(1),f_j(0))$ and let
\[
A_{nj}=\left\{|\widehat{V}_{nj}^E(a)|\leq n_j^\gamma,|\widehat{V}_{nj}^B(a)|\leq \log n_j\right\},
\qquad
\text{where }
\frac{1}{6(q-1)}<\gamma<\frac{\holder}{3}-\frac{1}{6},
\]
(such a choice of $\gamma$ is allowed since $\holder>3/4$ and $q>6$) and write $A_{nj}'=K_{nj}\cap A_{nj}$.
Then by Lemma \ref{lem:tail} we have
\[
\prob(K_{nj}\cap A_{nj}^c)
\leq
Kn_j^{-3\gamma(q-1)}.
\]
Hence, since $|\widehat{V}_{nj}^E(a)-\widehat{V}_{nj}^B(a)|\leq 2n_j^{1/3}$,
we have for $a\in(f_j(1),f_j(0))$,
\[
\mathbb{E}|\widehat{V}_{nj}^E(a)-\widehat{V}_{nj}^B(a)|\indicator_{K_{nj}}
\leq
\mathbb{E}|\widehat{V}_{nj}^E(a)-\widehat{V}_{nj}^B(a)|\indicator_{A_{nj}'}
+
2n_j^{1/3}\cdot Kn_j^{-3\gamma(q-1)}.
\]
The second term on the right hand side is $o(n_j^{-1/6})$.
To bound the first term, write
\[
\begin{split}
\mathbb{E}|\widehat{V}_{nj}^E(a)-\widehat{V}_{nj}^B(a)|\indicator_{A_{nj}'}
&=
\int_0^{\epsilon_{nj}}
\prob\left(|\widehat{V}_{nj}^E(a)-\widehat{V}_{nj}^B(a)|>x,A_{nj}'\right)\,\md x\\
&\qquad+
\int_{\epsilon_{nj}}^{2n_j^\gamma}
\prob\left(|\widehat{V}_{nj}^E(a)-\widehat{V}_{nj}^B(a)|>x,A_{nj}'\right)\,\md x\\
&\leq
\epsilon_{nj}+\int_{\epsilon_{nj}}^{2n_j^\gamma}
\prob\left(|\widehat{V}_{nj}^E(a)-\widehat{V}_{nj}^B(a)|>x,A_{nj}'\right)\,\md x,
\end{split}
\]
where $\epsilon_{nj}=n_j^{-1/3+\gamma}(\log n_j)^2$.
By definition of $\widehat{V}_{nj}^S(a)$ for $S=E,B$ and $j=0,1,\ldots,J$, we have:
\begin{equation*}
\widehat{V}_{nj}^S(a)
=
\argmax_{g_j(a)+n_j^{-1/3}u\in[0,1]}
Z_{nj,g_j(a)}^S(u)
\end{equation*}
where
\begin{equation*}
\begin{split}
Z_{nj,t}^S(u)
&=
n_j^{2/3}
\left\{
M_{nj}^S\big(t,t+n_j^{-1/3}u\big]
+
F_j\big(t,t+n_j^{-1/3}u\big]-f_j(t)n_j^{-1/3}u
\right\}.
\end{split}
\end{equation*}
Since $n_j^{2/3}|M_{nj}^E(t)-M_{nj}^B(t)|\leq\delta_{nj}$ on the event $A_{nj}'$, we can only have
$|\widehat{V}_{nj}^E(a)-\widehat{V}_{nj}^B(a)|>x$, if
\[
|Z_{nj,g_j(a)}^B\left(\widehat{V}_{nj}^B(a)\right)-Z_{nj,g_j(a)}^B(u)|\leq 2\delta_{nj},
\]
for some $u\in[-n_j^{1/3}g_j(a),n_j^{1/3}(1-g_j(a))]$,
such that $|u-\widehat{V}_{nj}^B(a)|>x$.
From here on the argument is identical to the one in the proof of Corollary 3.1 in  \cite{groeneboom-hooghiemstra-lopuhaa1999},
and it follows that
\[
\begin{split}
\mathbb{E}|\widehat{V}_{nj}^E(a)-\widehat{V}_{nj}^B(a)|\indicator_{A_{nj}'}
&\leq
\epsilon_{nj}+\int_{\epsilon_{nj}}^{2n_j^\gamma}
\prob\left(|\widehat{V}_{nj}^E(a)-\widehat{V}_{nj}^B(a)|>x,A_{nj}'\right)\,\md x\\
&\leq
\epsilon_{nj}+\int_{\epsilon_{nj}}^{2n_j^\gamma}
\prob\left(\widehat{V}_{nj}^B\text{ has a jump in }\left[a-\frac{2\delta_{nj}}{n_j^{1/3}x},a+\frac{2\delta_{nj}}{n_j^{1/3}x}\right]\right)\,\md x.
\end{split}
\]
According to Lemma~\ref{lem:jump} the right hand side is bounded by
\[
\epsilon_{nj}+\beta_j(1+o(1))
\int_{\epsilon_{nj}}^{2n_j^\gamma}
\left(\frac{2n_j^{-1/3+1/q}\log n}{x}
\vee
\left(n^{-\holder/3}(\log n)^{\holder+1}\right)\right)\,\md x
=
o\left(n_j^{-1/6}\right),
\]
uniformly in $a\in(f_j(1),f_j(0))$, since $q>6$ and $\gamma<\holder/3-1/6$.
The lemma then follows from the Markov inequality.
\tqed

\bigskip

\noindent
\textbf{Proof of Lemma~\ref{lem:mixing}.}
Let $t\in(0,1)$ arbitrary and take $0<s_1\leq s_2\leq\cdots\leq s_k=t<t+d=u_1\leq u_2\leq\cdots\leq u_l<1$.
Consider events
\[
\begin{split}
E_1
&=
\bigcap_{j=0}^J
\left\{\hf_{nj}^W(s_1)\in B_{1j},\ldots,\hf_{nj}^W(s_k)\in B_{kj}\right\},\\
E_2
&=
\bigcap_{j=0}^J
\left\{\hf_{nj}^W(u_1)\in C_{1j},\ldots,\hf_{nj}^W(u_l)\in C_{lj}\right\},
\end{split}
\]
for Borel sets $B_{1j},\ldots,B_{kj}$ and $C_{1j},\ldots,C_{lj}$ of $\R$, for $j=0,1,\ldots,J$.
Note that cylinder sets of the form $E_1$ and $E_2$ generate the
$\sigma$-algebras $A\in\sigma\{\hf_{nj}^W(s):j=0,1,\ldots,J,\,0<s\leq t\}$  and
$B\in\sigma\{\hf_{nj}^W(u):j=0,1,\ldots,J,\,t+d\leq u<1\}$, respectively.
Let $\CM_I Z$ denote the least concave majorant of $Z$ on the interval $I$.
Define the event $S=\bigcap_{j=0}^J S_j$, where
\[\begin{split}
S_j
=
&
\Big\{\left[\CM_{[0,1]}F_{nj}^W\right](u)=\left[\CM_{[0,t+d/2]}F_{nj}^W\right](u), \text{ for  } u\in [0,t] \\
&\qquad
\text{ and }
\left[\CM_{[0,1]}F_{nj}^W\right](u)=\left[\CM_{[t+d/2,1]}F_{nj}^W\right](u),
\text{ for } u\in [t+d,1]\Big\},
\end{split}\]
where $F_{nj}^W$ is defined in~\eqref{eq:def F_nj^S} and~\eqref{eq:def F_n0^S}.
Let $E_1'=E_1\cap S$ and $E_2'=E_2\cap S$.
Then $E_1'$ only depends on the processes $F_{n0}^W,F_{n1}^W,\ldots,F_{nJ}^W$ before time $t+d/2$
and $E_2'$ only depends on the processes $F_{n0}^W,F_{n1}^W,\ldots,F_{nJ}^W$  after time $t+d/2$.
Hence, since each process $W_{nj}$, for $j=0,1,\ldots,J$, is distributed like Brownian motion,
it has independent increments, so that the events $E_1'$ and $E_2'$ are independent.
This yields
\[
|\prob(E_1\cap E_2)-\prob(E_1)\prob(E_2)|\leq 3\prob\left(S^c\right)
\leq
3\sum_{j=0}^J\prob\left(S_j^c\right).
\]
Note that the concave majorants of $F_{nj}^W$ on $[0,1]$ and on $[0,t]$ coincide on $[0,t]$
as soon as they coincide at the boundary points.
Hence, for each $j=0,1,\ldots,J$, we have that on the event $S_j^c$ both concave majorants differ at $t$, which implies that they differ
on the interval $[t-dn_j^{-1/3},t+dn_j^{-1/3}]$.
According to Lemma~1.3 in \cite{kulikov-lopuhaa2006SPL}, this happens with probability bounded by $8\exp(-C_jnd^3)$,
where the constant $C_j>0$ only depends on $f_j$ and $c_j$.
This proves the lemma.
\tqed

\bigskip

\noindent
{\bf Proof of Lemma~\ref{lem:transitionV}.}
For $j=0,1,\dots,J$, we have
\[
n_j^{1/3}
\left(
L_j(\widehat{U}_{nj}^W(a))-L_j(g_j(a))
\right)
=
\argmax_{u\in I_{nj}(a)}
\left\{
W_{g_j(a),j}(u)-
\frac{|f_j'(t)|}{2(L_j'(t))^2}u^2+R_{nj}(a,u)
\right\}
\]
where $W_{tj}$ and $I_{nj}(a)$ are defined in~\eqref{eq:def Inj Wtj supplement} and
\[
\begin{split}
R_{nj}(a,u)
&=
n_j^{2/3}\left(F_{j}\circ L_j^{-1}-aL_j^{-1}\right)\left(L_j(g_j(a))+n_j^{-1/3}u\right)\\
&\qquad-
n_j^{2/3}
\left(
F_j(g_j(a))-ag_j(a)
\right)+
\frac{|f_j'(t)|}{2(L_j'(t))^2}u^2.
\end{split}
\]
Let
$$\widetilde U^W_{nj}(a)=\argmax_{u\in [-\log n,\log n]}
\left\{
W_{g_j(a),j}(u)-
\frac{|f_j'(t)|}{2(L_j'(t))^2}u^2+R_{nj}(a,u)
\right\}.$$
Since $[-\log n,\log n]\subset I_{nj}(a),$ we can only have
$n_j^{1/3}\left(L_j(\widehat{U}_{nj}^W(a))-L_j(g_j(a))\right)\neq \widetilde U^W_{nj}(a)$,
if
$$
n_j^{1/3}
\left|
L_j(\widehat{U}_{nj}^W(a))-L_j(g_j(a))
\right|\geq \log n.
$$
Furthermore, $|\widetilde U^W_{nj}(a)|\leq\log n$ and
\[
\left|
L_j(\widehat{U}_{nj}^W(a))-L_j(g_j(a))
\right|
\leq \sup_{t\in[0,1]}|L_j'(t)||\widehat{U}_{nj}^W(a)-g_j(a)|
\leq \sup_{t\in[0,1]}|L_j'(t)|,
\]
whence
\[\begin{split}
&\E\left|n_j^{1/3}
\left(
L_j(\widehat{U}_{nj}^W(a))-L_j(g_j(a))
\right)-\widetilde U^W_{nj}(a)\right|\\
&\quad\leq
\left(n_j^{1/3}\sup_{t\in[0,1]}|L_j'(t)|+\log n\right)\prob\left(n_j^{1/3}\sup_{t\in[0,1]}|L_j'(t)|
\left|\widehat{U}_{nj}^W(a)-g_j(a)
\right|\geq \log n\right).
  \end{split}\]
From Lemma \ref{lem:tail}, the right hand term is of the order $o(n^{-1/6}/\log n)$. This means that it remains to show that
\begin{equation}\label{eq:transitionU}
\E\left|\widetilde U^W_{nj}(a)
-\widetilde V_{nj}(g_j(a))\right|\leq C\frac{n^{-1/6}}{\log n}.
\end{equation}
It follows from the assumptions on $L_j$ and $f_j$ that there exists $K>0$ depending only on $f_j$ and $L_j$ such that for all $a\in(f_j(1),f_j(0))$,
$$\sup_{|u|\leq\log n}|R_{nj}(a,u)|\leq Kn^{-\holder/3}(\log n)^3.$$
By assumption, $s>3/4$, so Lemma 5 from \cite{durot2007} yields (\ref{eq:transitionU}).
\tqed

\bigskip

\noindent
{\bf Proof of Lemma~\ref{lem:11}.}
By assumption, $L_j'$ is bounded from below and $|L_j''|$ is bounded from above for every $j=0,1,\dots,J$.
Moreover, integrating~(\ref{eq:tailBW}) with $g_j=g$ the inverse of $f_0$,
one obtains that there exists $K>0$ depending on $f$ and $L_j$ only such that
$$
\E|\widehat{U}_{nj}^W(a)-g(a)|^2\leq Kn^{-2/3}
$$
for every $a\in\R$.
Therefore, together with a Taylor expansion, \eqref{eq:transition complete} yields
\[
\begin{split}
&
n^{1/3}\int_0^1
|\widehat{f}_{ni}^E(t)-\widehat{f}_{nj}^E(t)|\,\md t\\
&=
n^{1/3}\int_{f(1)}^{f(0)}
\left|
\frac{L_i(\widehat{U}_{ni}^W(a))-L_i(g(a))}{L_i'(g(a))}
-
\frac{L_j(\widehat{U}_{nj}^W(a))-L_j(g(a))}{L_j'(g(a))}
\right|\,\md a
+
o_p(n^{-1/6})\\
&=
n^{1/3}\int_{f(1)+\delta_n}^{f(0)-\delta_n}
\left|
\frac{L_i(\widehat{U}_{ni}^W(a))-L_i(g(a))}{L_i'(g(a))}
-
\frac{L_j(\widehat{U}_{nj}^W(a))-L_j(g(a))}{L_j'(g(a))}
\right|\,\md a
+
o_p(n^{-1/6}),
\end{split}
\]
where $\delta_n=n^{-1/3}(\log n)^2$.
From Lemma \ref{lem:transitionV} together with the fact that $|\widetilde V_{nj}(t)|\leq\log n$ for every $j$ and $t$, it then follows that
\[n^{1/3}\int_0^1
|\widehat{f}_{ni}^E(t)-\widehat{f}_{nj}^E(t)|\,\md t\\
=
\int_{f(1)}^{f(0)}
\left|
\frac{\widetilde{V}_{ni}(g(a))}{L_i'(g(a))c_i^{1/3}}
-
\frac{\widetilde{V}_{nj}(g(a))}{L_j'(g(a))c_j^{1/3}}
\right|\,\md a
+
o_p(n^{-1/6}).
\]
The change of variable $t=g(a)$ finally yields the lemma.
\tqed

\bigskip

\noindent
\textbf{Proof of Lemma~\ref{lem:approx tilde Vn}:}
Similar to~\eqref{eq:bound tailprob BM},
one can prove that there exists $K_1>0$, such that for all $s,t\in(0,1)$, $j=0,1,\dots,J$ and $x>0$,
\begin{equation}
\label{eq:tailV}
\prob(|V_{tj}(s)|>x)\leq 2\exp(-K_1x^3),
\end{equation}
and integrating (\ref{eq:tailV}) yields
\begin{equation}
\label{eq:espV}
\E|V_{tj}(s)|^{\gamma}\leq K_2,
\end{equation}
for all $s,t\in[0,1]$, $j=0,1,\ldots,J$, and $\gamma>0$, where $K_2>0$ depends only on $\gamma$, the~$L_j$'s and~$f$.
To prove~\eqref{eq:tildeV VV'}, first consider the cases $j=1,2,\dots,J$.
By H\"older's inequality we have
$$
\E\left|\widetilde{V}_{nj}(t)-V_{tj}(t)\right|^r
\leq
\sqrt{\E\left|\widetilde{V}_{nj}(t)-V_{tj}(t)\right|^{2r}\prob\left(\widetilde{V}_{nj}(t)\neq V_{tj}(t)\right)}.
$$
But, according to~\eqref{eq:def tilde Vnj} and~\eqref{eq:def Vtj},
the argmax~$\widetilde{V}_{nj}(t)$ can differ from $V_{tj}(t)$ only if $|V_{tj}(t)|>\log n$, so from~(\ref{eq:tailV}),
we obtain
$$
\prob \left(\widetilde{V}_{nj}(t)
\neq
V_{tj}(t)\right)\leq 2\exp(-K_1(\log n)^3).
$$
Moreover, $|\widetilde{V}_{nj}(t)-V_{tj}(t)|\leq 2|V_{tj}(t)|$.
Therefore, (\ref{eq:espV}) with $\gamma=2r$ yields
\begin{equation}
\label{eq:tildeV V exp}
\E\left|\widetilde{V}_{nj}(t)-V_{tj}(t)\right|^r\leq 2^{r+1/2}K_2^{1/2}\exp(-K_1(\log n)^3/2),
\end{equation}
uniformly in $t\in(0,1)$ and $j=1,2,\dots,J$.
This proves~\eqref{eq:tildeV VV'} for the cases $j=1,2,\ldots,J$.
The case $j=0$ in~\eqref{eq:tildeV VV'} is proven completely similar.

We proceed with proving~\eqref{eq:tildeV V logn}.
From the convexity of $x\mapsto x^r$ we deduce that for all $s,t\in[0,1]$ and $j=0,1,\ldots,J$,
\begin{equation}
\label{eq:bound r-moment}
\E
\left|\widetilde{V}_{nj}(t)-V_{tj}(s)\right|^r\leq 2^{r-1}\E\left|\widetilde{V}_{nj}(t)-V_{tj}(t)\right|^r+
2^{r-1}\E\left|{V}_{tj}(t)-V_{tj}(s)\right|^r.
\end{equation}
To deal with the second term on the right-hand side of~\eqref{eq:bound r-moment}, define the event
\[
A_n=\left\{|{V}_{tj}(t)|\leq\log n,|{V}_{tj}(s)|\leq\log n\right\},
\]
and note that from the H\"older inequality, it follows that for all $s,t\in[0,1]$,
\[
\begin{split}
\E\left|{V}_{tj}(t)-V_{tj}(s)\right|^r
\leq
\E\left[\left|
V_{tj}(t)-V_{tj}(s)\right|^r\indicator_{A_n}
\right]
+
\sqrt{\E\left|{V}_{tj}(t)-V_{tj}(s)\right|^{2r} \prob\left(A_n^c\right)}.
\end{split}\]
Therefore, together with~\eqref{eq:tailV} and~\eqref{eq:espV}, we find
\[
\E\left|{V}_{tj}(t)-V_{tj}(s)\right|^r
\leq
\E\left[\left|
V_{tj}(t)-V_{tj}(s)\right|^r\indicator_{A_n}
\right]
+
K_3\exp\left(-\frac{K_1}{2}(\log n)^3\right).
\]
Since $|s-t|\leq An^{-1/3}\log n$, from~\eqref{eq:def Vtj} and~\regul\ it follows that the difference between the drifts
of the processes corresponding to $V_{tj}(t)$ and $V_{tj}(s)$
satisfies condition~(ii) of Lemma~5 in~\cite{durot2007}.
Since $r<2\holder$, from Lemma~5 in~\cite{durot2007} we conclude that
$$
\E\left[\left|
V_{tj}(t)-V_{tj}(s)\right|^r\indicator_{A_n}
\right]
=
\E\left(\left|{V}_{tj}(t)-V_{tj}(s)\right|^r\indicator_{|{V}_{tj}(t)|\leq\log n}\indicator_{|{V}_{tj}(s)|\leq\log n}\right)
\leq
K_4\left(\frac{n^{-1/6}}{\log n}\right)^r.
$$
Hence, for all $j=0,1,\ldots,J$, and $s,t\in[0,1]$, such that $|s-t|\leq An^{-1/3}\log n$,
\begin{equation}
\label{eq:Er}
\E\left|\widetilde{V}_{nj}(t)-V_{tj}(s)\right|^r
\leq
2^{r-1}\E\left|\widetilde{V}_{nj}(t)-V_{tj}(t)\right|^r
+
K_5\left(\frac{n^{-1/6}}{\log n}\right)^r.
\end{equation}
For $j=1,2,\ldots,J$, the first term on the right hand side of~\eqref{eq:Er}
can be bounded together with~\eqref{eq:tildeV V exp}.
In particular we obtain
\begin{equation}\label{eq:tildeV V}
\E^{1/r}\left|\widetilde{V}_{nj}(t)-V_{tj}(t)\right|^r=o(1/\log n),
\end{equation}
uniformly in $t\in(0,1)$ and $j=1,2,\dots,J$.
Next, consider the case $j=0$.
By H\"older's inequality, for every $\delta>0$  we have
\begin{equation}
\label{eq:j=0}
\E|\widetilde{V}_{n0}(t)-V_{t0}(t)|^r
\leq
\delta
+
\sqrt{\E|\widetilde{V}_{n0}(t)-V_{t0}(t)|^{2r}\prob(|\widetilde{V}_{n0}(t)-V_{t0}(t)|>\delta)},
\end{equation}
where $\widetilde{V}_{n0}(t)$ is defined in~\eqref{eq:def tilde Vnj}, and can be written as
$$
\widetilde{V}_{n0}(t)=\argmax_{|u|\leq\log n}\left\{\widetilde W_{t0}(u)-\frac{|f_0'(s)|}{2(L_0'(s))^2}u^2+R_n(s,t,u)\right\},
$$
with
$$
\sup_{|u|\leq\log n}|R_n(s,t,u)|
\leq
\left|\frac{|f_0'(s)|}{2(L_0'(s))^2}-\frac{|f_0'(t)|}{2(L_0'(t))^2}\right|(\log n)^2
+
\sup_{|u|\leq\log n}\left|\widetilde W_{t0}(u)- W_{t0}(u)\right|.
$$
In order to bound the probability in~\eqref{eq:j=0}, we use
Proposition 1 in~\cite{durot2002} (see also the comments just above this proposition),
i.e.,
\begin{equation}
\label{eq: Proposition 1}
\begin{split}
\prob(|\widetilde{V}_{n0}(t)-V_{t0}(t)|>\delta)
&\leq
\prob\left\{2\sup_{|u|\leq\log n}|R_n(s,t,u)|>x\delta^{3/2}\right\}
+
Cx\log n\\
&+
\prob(|V_{t0}(t)|>\log n),
\end{split}
\end{equation}
for $\delta=n^{-\epsilon}$, for some $\epsilon<1/9$, and
$x=\delta^{-3/2}n^{-1/6}(\log n)^2$.
From~\eqref{eq:Wt0 in Wtj} and~(\ref{eq:def tildeWtj}),
we obtain that for all $|s-t|\leq An^{-1/3}\log n$,
$$
\sup_{|u|\leq\log n}|R_n(s,t,u)|
\leq
K_6n^{-\holder/3}(\log n)^{3}+J\max_{1\leq j\leq J}\sup_{|u|\leq\gamma_n,|u-v|\leq \delta_n}|W_{tj}(u)-W_{tj}(v)|,
$$
where $\gamma_n=K_7\log n$ and $\delta_n=K_7n^{-1/3}(\log n)^2$.
This means that for $y>K_6n^{-\holder/3}(\log n)^{3}$,
\[
\prob
\left(
\sup_{|u|\leq\log n}|R_n(s,t,u)|>y
\right)
\leq
J
\prob
\left(
J\sup_{|u|\leq\gamma_n,|u-v|\leq \delta_n}|W_{tj}(u)-W_{tj}(v)|>\frac{y}{2}
\right)
\]
Let $M=\lfloor K_7 n^{1/3}(\log n)^{-1}\rfloor$
and for $k\in\{-M,\ldots,M\}$, let $t_k=kn^{-1/3}(\log n)^2$.
Then, if $|u|\leq\gamma_n,|u-v|\leq \delta_n$,
it follows that $|u-t_k|\leq n^{-1/3}(\log n)^2$ and $|v-t_k|\leq 2n^{-1/3}(\log n)^2$,
for some $k\in\{-M,\ldots,M\}$.
This means that, together with using
\[
|W_{tj}(u)-W_{tj}(v)|
\leq
|W_{tj}(u)-W_{tj}(t_k)|+|W_{tj}(t_k)-W_{tj}(v)|,
\]
the probability on the right hand side can be bounded by
\[
2J
\sum_{k=-M}^M
\prob
\left(
\sup_{|v-t_k|\leq 2n^{-1/3}(\log n)^2}
|W(t_k)-W(v)|
>
\frac{y}{4J}
\right),
\]
where $W$ is a standard Brownian motion.
We conclude that
\begin{equation}
\label{eq:sup BM}
\begin{split}
\prob
\left(
\sup_{|u|\leq\log n}|R_n(s,t,u)|>y
\right)
&\leq
2J(2M+1)
\prob
\left(
\sup_{|v|\leq 2n^{-1/3}(\log n)^2}
|W(v)|
>
\frac{y}{4J}
\right)\\
&\leq
K_8n^{1/3}(\log n)^{-1}
\exp
\left(
-\frac{K_9 y^2n^{1/3}}{(\log n)^2}
\right).
\end{split}
\end{equation}
Possibly enlarging the constant $K_8$, this inequality still holds for all $y>0$ since for a large enough $K_8$,
the right hand term becomes greater than one for $y\in(0,2K_6n^{-\holder/3}(\log n)^3)$.
Therefore, from~\eqref{eq: Proposition 1},
\[
\begin{split}
\prob(|\widetilde{V}_{n0}(t)-V_{t0}(t)|>\delta)
\leq
2\exp\left(
-K_8(\log n)^2
\right)
&+
Cn^{3\epsilon/2-1/6}(\log n)^3\\
&+
2\exp(-K_1(\log n)^3).
\end{split}
\]
Moreover, from the definitions of $\widetilde V_{n0}$ and $V_{t0}$ in~\eqref{eq:def tilde Vnj}
and~\eqref{eq:def Vt0}, together with the fact that~$W_{t0}$ and~$\widetilde{W}_{t0}$ are both distributed
as standard Brownian motion,
it is easy to see that $\E|\widetilde V_{n0}(t)|^\gamma\leq \E|V_{t0}(t)|^\gamma$, for every~$\gamma>0$,
so it follows from~(\ref{eq:espV}) that for every $\gamma>0$,
$$
\E|\widetilde{V}_{n0}(t)|^\gamma
\leq
\E|{V}_{t0}(t)|^\gamma=O(1),
$$
uniformly in $t\in(0,1)$.
Hence, we deduce from (\ref{eq:j=0}) that (\ref{eq:tildeV V}) also holds with $j=0$.
From~(\ref{eq:Er}) we then obtain
$$
\E^{1/r}\left|\widetilde{V}_{nj}(t)-V_{tj}(s)\right|^r=o(1/\log n)
$$
uniformly in $j=0,1,\dots,J$ and $s,t$ such that $|s-t|\leq An^{-1/3}\log n$.
\tqed

\bigskip

{\bf Proof of Lemma~\ref{lem:limit variance}}
For $k=1,2$, we have
$$
v_{nk}=2\int_0^1\int_s^1
\mathrm{cov}(Y_ {nk}(t),Y_ {nk}(s))\,\md t\,\md s.
$$
Note that by definition of $\widetilde V_{nj}(t)$ in~\eqref{eq:def tilde Vnj},
the random variable $Y_{nk}(t)$, defined in~\eqref{eq:def Ynk}, depends only on the increments of $W_{nj}$ over a neigbourhood
of $L_j(t)$ with radius of the order~$O(n^{-1/3}\log n)$,
for $j=0,1,\dots,J$.
But for every $s,t\in[0,1]$ and $j=0,1,\dots,J$, we have
$$
|L_j(t)-L_j(s)|\geq|t-s|\inf_{u\in[0,1]}|L_j'(u)|,
$$
where the infimum is positive according to~\embed.
Setting $a_n=An^{-1/3}\log n$, for some large enough $A>0$,
we find that $Y_{nk}(t)$ is independent of $Y_{nk}(s)$ for every $|t-s|\geq a_n$.
This means that in for $k=1$,
\begin{equation}
\label{eq:v1}
\begin{split}
v_{n1}
&=
2\int_0^1\int_s^{1\wedge(s+a_n)}\mathrm{cov}(Y_ {n1}(t),Y_ {n1}(s))\,\md t\,\md s\\
&=
2\sum_{i<j}\sum_{l<m}
\int_0^1\int_s^{1\wedge(s+a_n)}
|f_0'(s)||f_0'(t)|
\widetilde{C}_{n,iljm}(s,t)\,\md t\,\md s,
\end{split}
\end{equation}
where
$$
\widetilde{C}_{n,iljm}(s,t)
=
\mathrm{cov}\left(
\left|
\frac{\widetilde{V}_{ni}(s)}{L_i'(s)c_i^{1/3}}
-
\frac{\widetilde{V}_{nj}(s)}{L_j'(s)c_j^{1/3}}
\right|,
\left|
\frac{\widetilde{V}_{nl}(t)}{L_l'(t)c_l^{1/3}}
-
\frac{\widetilde{V}_{nm}(t)}{L_m'(t)c_m^{1/3}}
\right|
\right),
$$
with $\widetilde{V}_{nj}(t)$ defined in~\eqref{eq:def tilde Vnj}.
Next, we approximate both $\widetilde{V}_{nj}(s)$ and $\widetilde{V}_{nj}(t)$ with the variable~$V_{tj}(s)$,
defined in~\eqref{eq:def Vtj} and~\eqref{eq:def Vt0}.
According to~\eqref{eq:tildeV V logn}, we have
$\E^{1/r}|\widetilde{V}_{nj}(t)-V_{tj}(s)|^r=o(1/\log n)$,
uniformly in $j=0,1,\dots,J$ and $s,t$ such that $|s-t|\leq a_n$.
It thus follows from~(\ref{eq:v1}), together with H\"older's inequality, that
\[
v_{n1}
=
2\sum_{i<j}\sum_{l<m}
\int_0^1\int_s^{1\wedge(s+a_n)}
|f_0'(s)|^2C_{iljm}(s,t)\,\md t\,\md s+o(n^{-1/3}),
\]
where
$$
C_{iljm}(s,t)
=
\mathrm{cov}\left(\left|
\frac{{V}_{ti}(s)}{L_{i}'(s)c_i^{1/3}}
-
\frac{{V}_{tj}(s)}{L_{j}'(s)c_j^{1/3}}
\right|,\left|
\frac{{V}_{sl}(s)}{L_{l}'(s)c_{l}^{1/3}}
-
\frac{{V}_{sm}(s)}{L_{m}'(s)c_{m}^{1/3}}
\right|\right).
$$
Let $d_j(s)=|f_0'(s)|/(2L_j'(s)^2)$.
From~\eqref{eq:def Vtj} and~\eqref{eq:def Inj Wtj}, we have for $j=1,2,\ldots,J$,
\[
\begin{split}
&
d_j(s)^{2/3}V_{tj}(s)\\
&=
\argmax_{u\in\R}
\left\{
W_{tj}\left(d_j(s)^{-2/3}u\right)
-
d_j(s)^{-1/3}u^2
\right\}\\
&=
\argmax_{u\in\R}
\left\{
W_{sj}\left(d_j(s)^{-2/3}
\left[
u+n_j^{1/3}\left(L_j(t)-L_j(s)\right)d_j(s)^{2/3}
\right]
\right)
-
d_j(s)^{-1/3}u^2
\right\}\\
&\stackrel{d}{=}
\argmax_{u\in\R}
\left\{
W_{j}\left(
u+n_j^{1/3}\left(L_j(t)-L_j(s)\right)d_j(s)^{2/3}
\right)
-
u^2
\right\}\\
&=
\argmax_{u\in\R}
\left\{
W_{j}\left(
u+n^{1/3}(t-s)|f_0'(s)/2|^{2/3}\frac{c_j^{1/3}}{L_j'(s)^{1/3}}
\right)
-
u^2
+
R_n'(s,t,u)
\right\}
\end{split}
\]
where $W_j$ is a standard Brownian motion and where for every $|t-s|\leq a_n$,
$$
\sup_{|u|\leq\log n}|R_n'(s,t,u)|
\leq
\sup_{|u|\leq K\log n, |u-v|\leq K'n^{-1/3}(\log n)^2}|W_{j}(u)-W_{j}(v)|.
$$
For $j=1,2,\ldots,J$, let $\zeta_j$ be defined by~\eqref{eq: def zetaj}.
By partitioning the interval $[-K\log n,K\log n]$ into subintervals of length $n^{-1/3}\log n$

We can bound tail probabilities of $\sup |R_n'(s,t,u)|$ similar to~\eqref{eq:sup BM}.
Then using the same reasoning as in~\eqref{eq:j=0} and~\eqref{eq: Proposition 1},
from which we obtained~\eqref{eq:tildeV V} for the case $j=0$,
we now conclude that for $j=1,2,\ldots,J$,
\[
\mathbb{E}
\left|
\left(\frac{|f_0'(s)|}{2(L_{j}'(s))^2}\right)^{2/3}V_{tj}(s)
-
\zeta_j\left(n^{1/3}(t-s)|f_0'(s)/2|^{2/3}\frac{c_j^{1/3}}{L_j'(s)^{1/3}}\right)
\right|
=
o(1/\log n),
\]
or equivalently
$$
\mathbb{E}
\left|
\frac{V_{tj}(s)}{L_{j}'(s)c_j^{1/3}}
-
\left(\frac{4L_{j}'(s)}{c_j|f_0'(s)|^2}\right)^{1/3}
\zeta_j\left(n^{1/3}(t-s)|f_0'(s)/2|^{2/3}\frac{c_j^{1/3}}{L_j'(s)^{1/3}}\right)
\right|
=
o(1/\log n).
$$
Change of variable $t'=n^{1/3}(t-s)|f'(s)/2|^{2/3}$, then gives
\[
\begin{split}
n^{1/3}v_{n1}
=
8\sum_{i<j}\sum_{l<m}
\int_0^1\int_0^{a_n'}
\mathrm{cov}&\left(|Y_{si}(t')-Y_{sj}(t')|,|Y_{sl}(0)-Y_{sm}(0)|\right)\,\md t'\,\md s
+
o(1),
\end{split}
\]
where $a_n'=A\log n|f_0'(s)/2|^{2/3}$ and where for $j=1,2,\ldots,J$,
\[
Y_{sj}(t)
=
\frac{L_{j}'(s)^{1/3}}{c_j^{1/3}}
\zeta_{j}\left(\frac{c_j^{1/3}t}{L_{j}'(s)^{1/3}}\right),
\]
We finish the proof for the case $k=1$, by showing that
there exist absolute constants $K$ and~$K'$ such that
\begin{equation}
\label{eq: cov bound}
|\mathrm{cov}\left(|Y_{si}(t)-Y_{sj}(t)|, |Y_{sl}(0)-Y_{sm}(0)|\right)|\leq K\exp(-K't^3),
\end{equation}
because then, for $k=1$, the lemma follows from the dominated convergence theorem.
To prove~\eqref{eq: cov bound}, for $d>0$, let
\begin{equation*}
\zeta_{jd}(c)
=
\argmax_{|u|\leq d}
\left\{
W_j(u+c)-u^2
\right\},
\end{equation*}
define~$Y_{sj}^d(t)$ similar to $Y_{sj}(t)$, with $\zeta_j$ replaced by $\zeta_{jd}$.
If we take
\[
d=(t/4) \min_{1\leq j\leq J}c_j^{1/3}\left(\sup_{s\in[0,1]} L_j'(s)\right)^{-1/3},
\]
then $Y_{si}^d(t)-Y_{sj}^d(t)$ and $Y_{sl}^d(0)-Y_{sm}^d(0)$ are independent for all
$i,j,l,m=1,2,\ldots,J$ and $s\in[0,1]$,
since $\zeta_{jd}(c)$ only depends on $W_j(u)$, for $u\in[c-d,c+d]$.
Furthermore, $\zeta_{jd}(c)\ne \zeta_{j}(c)$ if and only if $|\zeta_{j}(c)|>d$,
and similar to~\eqref{eq:tailV} and~\eqref{eq:espV}, we find that
$\prob(|\zeta_j(c)|>d)\leq 2\exp\left(-K_1d^3\right)$ and
$\mathbb{E}|\zeta_j(c)|^\gamma|\leq K_2$,
for all $j=1,2,\ldots,J$, $\gamma>0$ and $c\in\R$.
Therefore, there exist $C_1,C_2>0$ such that
\begin{equation}
\label{eq:tail Y}
\prob
\left(
|Y_{si}(t)|>d
\right)
\leq
2\exp(-C_1d^3)
\end{equation}
and
\begin{equation}
\label{eq:bound Y}
\mathbb{E}|Y_{si}^d(t)|^2
\leq
\mathbb{E}|Y_{si}(t)|^2
\leq
\min_{1\leq j\leq J}
c_j^{-2/3}
\sup_{s\in[0,1]}L_j'(s)^{2/3}
\sup_{c\in\R}
\mathbb{E}|\zeta_j(c)|^2
\leq C_2.
\end{equation}
By H\"older's inequality, for any real valued random variables $X,X',Y$ and $Y'$, we have
\[
\left|\mathrm{cov}(X,Y)-\mathrm{cov}(X',Y')\right|
\leq
\mathbb{E}^{1/2}|X|^2\mathbb{E}^{1/2}|Y-Y'|^2
+
\mathbb{E}^{1/2}|Y'|^2\mathbb{E}^{1/2}|X-X'|^2,
\]
and similar to~\eqref{eq:tildeV V exp},
\[
E\left|Y_{si}^d(t)-Y_{sj}^d(t)\right|^2
\leq
C_3
\exp
\left(
-C_4t^3
\right),
\]
uniformly in $s\in[0,1]$ and $i,j=1,2,\ldots,J$.
Together with~\eqref{eq:tail Y} and~\eqref{eq:bound Y}, this proves~\eqref{eq: cov bound}.

The proof for the case $k=2$ is similar.
As before, from~\eqref{eq:tildeV V logn},
\[
v_{n2}
=
2\sum_{i=1}^J\sum_{j=1}^J
\int_0^1\int_s^{1\wedge(s+a_n)}
|f_0'(s)|^2C_{ij}(s,t)\,\md t\,\md s+o(n^{-1/3}),
\]
where
$$
C_{ij}(s,t)
=
\mathrm{cov}\left(\left|
\frac{{V}_{ti}(s)}{L_{i}'(s)c_i^{1/3}}
-
\frac{{V}_{t0}(s)}{L_{0}'(s)}
\right|,\left|
\frac{{V}_{sj}(s)}{L_{j}'(s)c_{j}^{1/3}}
-
\frac{{V}_{s0}(s)}{L_{0}'(s)}
\right|\right),
$$
and
\[
\mathbb{E}
\left|
\frac{V_{t0}(s)}{L_{0}'(s)}
-
\left(\frac{4L_{0}'(s)}{|f_0'(s)|^2}\right)^{1/3}
\widetilde{\zeta}_{t0}\left(n^{1/3}(t-s)|f_0'(s)/2|^{2/3}\frac{1}{L_0'(s)^{1/3}}\right)
\right|
=
o(1/\log n).
\]
where $\widetilde{\zeta}_{t0}$ is defined in~\eqref{eq: def zetat0}.
After change of variables $t'=n^{1/3}(t-s)|f_0'(s)/2|^{2/3}$, we obtain
\[
\begin{split}
n^{1/3}v_{n2}
=
2\sum_{i=1}^J\sum_{j=1}^J
\int_0^1\int_0^{a_n'}
\mathrm{cov}&\left(|Y_{si}(t')-Y_{s0}(t')|,|Y_{sj}(0)-Y_{s0}(0)|\right)\,\md t'\,\md s
+
o(1),
\end{split}
\]
where $a_n'=A\log n|f_0'(s)/2|^{2/3}$ and
\[
Y_{s0}(t)
=
L_{0}'(s)^{1/3}\widetilde{\zeta}_{t0}\left(\frac{t}{L_{0}'(s)^{1/3}}\right),
\]
Similar to the reasoning above, \eqref{eq: cov bound} can be shown for the case $j=m=0$,
so that the lemma follows for the case $k=2$, by application of the dominated convergence theorem.
\tqed

\subsection{Proof of Lemma~\ref{lem: tildefn}}
Note first that it suffices to prove the lemma for $[a,b]=[0,1]$,
since one can go from this specific case to the general case by considering
$(b-a)f_{0}(a+x(b-a))$ for $x\in[0,1]$ and the corresponding estimator $(b-a)\widetilde f_{n}(a+x(b-a))$.
For the sake of brevity we write $h$ instead of $h_n$.
For every $t\in[h,1-h]$, define
\begin{equation*}
\widetilde F_n(t)=\frac{1}{h}
\int_{t-h}^{t+h}
F_{n0}(x)K\left(\frac{t-x}{h}\right)\,\md x,
\end{equation*}
so that $\widetilde f_n=\widetilde F_n'$ on $[h,1-h]$, where we consider a right-derivative at $h$ and a left-derivative at $1-h$.
With $l=0,1,2$, the $l$th derivative of $\widetilde f_n$ on $[h,1-h]$ is given by
$$
\widetilde f_n^{(l)}(t)=\frac{1}{h^{2+l}}
\int_{t-h}^{t+h}
F_{n0}(x)K^{(l+1)}\left(\frac{t-x}{h}\right)\,\md x.
$$
It follows from Lemma \ref{lem:embed0} that
$$
\sup_{t\in[0,1]}|F_{n0}(t)-F_{0}(t)-\frac{1}{\sqrt n}B_{n0}\circ L_0(t)|=O_p(n^{-1+1/q}),
$$
where $F_{0}=F_1=\dots=F_J$ and $L_0$ and $B_{n0}$ are defined by~(\ref{eq:defL0}) and~(\ref{eq:defBn0}), respectively.
Combining this with the fact that $K^{(l+1)}$ is bounded and supported on $[-1,1]$, yields that on $[h,1-h]$,
$$
\widetilde f_n^{(l)}(t)
=
\frac{1}{h^{2+l}}
\int_{t-h}^{t+h}
\left(F_{0}(x)+\frac{1}{\sqrt n}B_{n0}\circ L_0(x)\right)
K^{(l+1)}\left(\frac{t-x}{h}\right)\,\md x+O_p\left(\frac{n^{-1+1/q}}{h^{l+1}}\right),
$$
where the $O_p$-term is uniform in $t\in[h,1-h]$.
Now, since $h\sim n^{-\hn}$, with the change of variable $u=(t-x)/h$ we obtain
\begin{equation}
\label{eq: tildefl}
\begin{split}
\widetilde f_n^{(l)}(t)
=
\frac{1}{h^{1+l}}
\int_{-1}^{1}
\bigg(
F_{0}(t-uh)
&+
\frac{1}{\sqrt n}B_{n0}\circ L_0(t-uh)\bigg) K^{(l+1)}(u)\,\md u\\
&+
O_p\left(n^{-1+1/q+(l+1)\hn}\right).
\end{split}
\end{equation}
We have $K^{(l)}(1)=K^{(l)}(-1)=0$, since $K$ is supported on $[-1,1]$ and is three times continuously differentiable on $\R$.
Therefore, we have for $t\in[h,1-h]$,
\[
\begin{split}
&
\left|\frac{1}{h^{1+l}\sqrt n}
\int_{-1}^{1}
B_{n0}\circ L_0(t-uh)K^{(l+1)}\left(u\right)\,\md u\right|\\
&\quad=
\left|\frac{1}{h^{1+l}\sqrt n}\int_{-1}^1 \left(B_{n0}\circ L_0(t-uh)-B_{n0}\circ L_0(t)\right)K^{(l+1)}\left(u\right)\,\md u\right|\\
&\quad\leq
\frac{1}{h^{1+l}\sqrt n}\sup_{x,y\in[L_0(0),L_0(1)],\, |x-y|\leq Ch}|B_{n0}(x)-B_{n0}(y)|
\int_\R
\left|K^{(l+1)}\left(u\right)\right|\,\md u,
\end{split}
\]
where $C=\sup_{t\in[0,1]}L_0'(t)$.
Although $B_{n0}$ is neither Brownian motion nor Brownian bridge, it is a linear combination of $B_{nj}$, $j=1,2,\ldots,J$.
For this reason, using standard properties of Brownian motion and Brownian bridge together with the fact that
$h\sim n^{-\hn}$, we conclude that
$$
\frac{1}{h^{1+l}\sqrt n}
\int_{-1}^{1}
B_{n0}\circ L_0(t-uh)K^{(l+1)}\left(u\right)\,\md u
=
O_p\left(n^{\hn(l+1/2)-1/2}\sqrt{\log n}\right),
$$
where the $O_p$ is uniform in $t\in[h,1-h]$.
Using that $\hn\leq1/5$, this yields
$$
\frac{1}{h^{1+l}\sqrt n}
\int_{-1}^{1}
B_{n0}\circ L_0(t-uh)K^{(l+1)}\left(u\right)\,\md u
=
O_p\left(n^{-2/5+l/5}\sqrt{\log n}\right).
$$
Since, by assumption $q>6$ and $\hn\leq 1/5$, together with (\ref{eq: tildefl}), this yields
\[
\widetilde f_n^{(l)}(t)
=
\frac{1}{h^{1+l}}
\int_{-1}^{1}
F_{0}(t-uh) K^{(l+1)}\left(u\right)\,\md u+O_p\left(n^{-2/5+l/5}\sqrt{\log n}\right),
\]
where the $O_p$-term is uniform in $t\in[h,1-h]$.
Uniformly in $t\in[h,1-h]$ and $u\in[-1,1]$ we have
$$
F_{0}(t-uh)=F_{0}(t)-uhf_0(t)+\frac{1}{2}u^2h^2f_0'(t)+O(h^3).
$$
Since, for $l=0,1,2$, we have $\int K^{(l+1)}(u)\,\md u=0$, we conclude that uniformly on $[h,1-h]$,
\[
\begin{split}
\widetilde f_n^{(l)}(t)
=
-\frac{f_0(t)}{h^{l}}\int_{-1}^{1} u K^{(l+1)}\left(u\right)\,\md u
&+
\frac{1}{2h^{l-1}}f_0'(t)\int_{-1}^{1} u^2K^{(l+1)}(u)\,\md u\\
&+
O(h^{2-l})+O_p\left(n^{-2/5+l/5}\sqrt{\log n}\right).
\end{split}\]
From our assumptions on $K$ together with integration by parts, it follows that
\begin{equation}\label{eq:K}
\int_{-1}^{1} uK^{(3)}(u)\,\md u=\int_{-1}^{1} u^2K^{(3)}(u)\,\md u=\int_{-1}^{1} uK^{(2)}(u)\,\md u=\int_{-1}^{1} u^2K^{(1)}(u)\,\md u=0.
\end{equation}
Therefore, in the case $l=2$ we obtain that $\widetilde f_n^{(2)}(t)=O_p(\sqrt {\log n})$ uniformly on $[h,1-h]$.
In the cases $l=0,1$, using~(\ref{eq:K}) together with the facts that
$\int uK^{(1)}(u)\,\md u=-1$ and $\int u^2K^{(2)}(u)\,\md u=2$, we obtain
\begin{equation}
\label{eq:order difference f and f'}
\widetilde f_n^{(l)}(t)=f_0^{(l)}(t)+O(h^{2-l})+O_p\left(n^{-2/5+l/5}\sqrt{\log n}\right),
\end{equation}
uniformly on $[h,1-h]$.

Now, first consider the boundary correction~\eqref{eq:tildefKosorok} on $[0,h]$ and $[1-h,1]$.
Since $\widetilde f_n^{(2)}$ is constant on $[0,h]$ and $[1-h,1]$, we get
\begin{equation}
\label{eq:order fn''}
\sup_{t\in[0,1]}|\widetilde f_n^{(2)}(t)|=O_p(\sqrt {\log n}).
\end{equation}
In particular, if $A_n'=\{\sup_{t\in[0,1]}|\widetilde f_n^{(2)}(t)|\leq \log n\}$, then $\prob(A'_n)\to1$.
Moreover, since $\widetilde f_n$ is twice differentiable, (\ref{eq: tildef'}) holds on $A_n'$
with $\holder=1$ and $\varepsilon_n=1/\log n$.
Using the definition of $\widetilde f_n$ on $[0,h]$ and $[1-h,1]$,
together with a Taylor expansion yields that on $[0,h)\cup(1-h,1],$
\[
\begin{split}
\widetilde f_{n}(t)-f_{0}(t)
&=
\widetilde f_{n}(u_{n})+\widetilde f_{n}'(u_{n})(t-u_{n})-f_{0}(u_{n})-f_{0}'(u_{n})(t-u_{n})+O(h^2)\\
&=
\widetilde f_{n}(u_{n})-f_{0}(u_{n})+\left(\widetilde f_{n}'(u_{n})-f_{0}'(u_{n})\right)(t-u_{n})+O(h^2).
\end{split}
\]
Together with~\eqref{eq:order difference f and f'}, we conclude that
$$
\sup_{t\in[0,1]}|\widetilde f_{n}(t)-f_{0}(t)|=O_p(n^{-2\hn})+O_p\left(n^{-2/5}\sqrt{\log n}\right).
$$
Similarly, on $[0,h)\cup(1-h,1]$ we have
\begin{eqnarray*}
\widetilde f_{n}'(t)-f_{0}'(t)=\widetilde f_{n}'(u_{n})-f_{0}'(u_{n})+O(h),\end{eqnarray*}
from which we conclude that
$$\sup_{t\in[0,1]}|\widetilde f_{n}'(t)-f_{0}'(t)|=O_p(n^{-\hn})+O_p\left(n^{-1/5}\sqrt{\log n}\right).$$

Therefore, if
$$
A_n''
=
\left\{
\sup_{t\in[0,1]}|\widetilde f_n(t)-f(t)|\leq n^{-1/3}/\log n\text{ and } \sup_{t\in[0,1]}|\widetilde f_n'(t)-f'(t)|\leq n^{-1/6}/\log n
\right\},
$$
then $\prob(A_n'')\to 1$,
and both~(\ref{eq: approxtildef}) and~(\ref{eq: approxtildef'}) hold on the event~$A_n''$.
Setting $A_n=A_n'\cap A_n''$ completes the proof of the lemma in the case of
boundary correction~\eqref{eq:tildefKosorok}.

Next, consider boundary correction~\eqref{def:boundary kernel estimate}.
First of all, note that for $s\in[-1,1]$,
\begin{equation}
\label{eq:coefficients booundary kernel}
\bcoefA(s)
=
\frac{\mathbb{K}^{(2)}(s)}{\mathbb{K}^{(0)}(s)\mathbb{K}^{(2)}(s)-\mathbb{K}^{(1)}(s)^2}
\quad\text{and}\quad
\bcoefB(s)
=
\frac{-\mathbb{K}^{(1)}(s)}{\mathbb{K}^{(0)}(s)\mathbb{K}^{(2)}(s)-\mathbb{K}^{(1)}(s)^2},
\end{equation}
where for $j=0,1,2$ and $s\in[-1,1]$,
\begin{equation}
\label{def:int K}
\mathbb{K}^{(j)}(s)=\int_{-1}^s u^jK(u)\,\md u.
\end{equation}
Note that the denominator in~\eqref{eq:coefficients booundary kernel} is bounded away from zero.
To see this, first note that from the Cauchy-Schwarz inequality, it follows that for every $s\in [-1,1]$,
$$
|\mathbb{K}^{(1)}(s)| \leq \int_{-1}^s |uK(u)| \md u \leq\left(\mathbb{K}^{(0)}(s)\right)^{1/2} \left(\mathbb{K}^{(2)}(s)\right)^{1/2},
$$
with equality if and only if $|uK^{1/2}(u)|=a|K^{1/2}(u)|$ Lesbesgue almost everywhere for some $a\geq0$.
Since for all $u$ in the support of $K$, we can have $|uK^{1/2}(u)|=a|K^{1/2}(u)|$ only for $u$ equal to $-a$ and $a$, this means that
$\mathbb{K}^{(0)}(s)\mathbb{K}^{(2)}(s)-\mathbb{K}^{(1)}(s)^2>0$,
for all $s\in [-1,1]$,
and hence, as being a continuous function of $s$ on $[-1,1]$, it has a strictly positive minimum.
This means that $\bcoefA(s)$ and $\bcoefB(s)$ are bounded and twice differentiable with bounded derivatives.
For $t\in[0,h]$, write
\[
\tilde f_n(t)-f_0(t)
=
\tilde f_n(t)-f_{0,B}(t)
+
f_{0,B}(t)-f_0(t),
\]
where
\begin{equation*}
f_{0,B}(t)
=
\int
\frac1{h}\bK\left(\frac{t-x}{h}\right)
\,\md F_{0}(x).
\end{equation*}
Then, since $t-h\leq 0$,
\[
\begin{split}
\tilde f_n(t)-f_{0,B}(t)
&=
\int_0^{t+h}
\frac{1}{h}
\bK\left(\frac{t-x}{h}\right)
\,\md (F_{n0}-F_0)(x)\\
&=-
\frac{1}{h}
\bK\left(\frac{t}{h}\right)
(F_{n0}(0)-F_0(0))\\
&\qquad+
\int_0^{t+h}
\frac{1}{h^2}
\bK'\left(\frac{t-x}{h}\right)
(F_{n0}(x)-F_0(x))
\,\md x,
\end{split}
\]
where $\bK'(u)=\md\bK(u)/\md u$.
After change of variables and using that both $K_{B,t}$ and~$K'_{B,t}$ are bounded,
thanks to Lemma~\ref{lem:embed0} and the fact that $h\leq 1-h$ for $n$ sufficiently large,
we get for $t\in[0,h]$,
\begin{equation}\label{eq:variance}
\begin{split}
\tilde f_n(t)-f_{0,B}(t)
&=-
\frac{1}{\sqrt nh}
\bK\left(\frac{t}{h}\right)
B_{n0}\circ L_{0}(0)\\
&\quad+
\int_{-1}^{t/h}
\frac{1}{\sqrt nh}
\bK'(u)
B_{n0}\circ L_{0}(t-hu)
\,\md u
+
O_p(h^{-1}n^{-1+1/q}).
\end{split}
\end {equation}
Using the fact that
\begin{equation}\label{eq:accB}
\sup_{t\in[0,h],\, |x-t|\leq h}\left|B_{n0}\circ L_{0}(t)-B_{n0}\circ L_{0}(x)\right|=O_{p}\left(\sqrt{h}\right),
\end{equation}
we obtain
\[
\begin{split}
&\int_{-1}^{t/h}
\frac{1}{\sqrt nh}
\bK'\left(u\right)
B_{n0}\circ L_{0}(t-hu)
\,\md u\\
&\qquad=
\int_{-1}^{t/h}
\frac{1}{\sqrt nh}
\bK'\left(u\right)B_{n0}\circ L_{0}(t)
\,\md u+O_{p}\left(\frac{1}{\sqrt {nh}}\right)\\
&\qquad=
\frac{1}{\sqrt nh}
\bK\left(\frac{t}{h}\right)B_{n0}\circ L_{0}(t)+O_{p}\left(\frac{1}{\sqrt {nh}}\right).
\end{split}
\]
Using again \eqref{eq:accB}, we then conclude from \eqref{eq:variance} that
\begin{equation}
\label{eq:order term 1}
\tilde f_n(t)-f_{0,B}(t)
=O_{p}\left(\frac{1}{\sqrt {nh}}\right),
\end{equation}
since $q>6$ and $h$ is of greater order than $n^{-2/3}$.
Next, consider $f_{0,B}(t)-f_0(t)$.
Note that for $t<h$:
\begin{equation}
\label{bias_expansion}
\begin{split}
f_{0,B}(t)
&=
\int_{0}^{t+h}
\frac{1}{h}
\bK\left(\frac{t-x}{h}\right)f_0(x)\,\md x\\
&=
\int_{-1}^{t/h}
\bK(u)f_0(t-hu)\,\md u
=
f_0(t)+O(h^2),
\end{split}
\end{equation}
because
\[
\int_{-1}^{t/h}
\bK(u)\,\md u=1
\quad\text{and}\quad
\int_{-1}^{t/h}
u\bK(u)\,\md u=0,
\]
by definition~\eqref{eq:def coef boundary kernel} of the coefficients in $\bK$.
Since, $\gamma\in(1/6,1/5]$, together with~\eqref{eq:order term 1}, this yields
\begin{equation*}
\sup_{t\in[0,h]}
|\tilde f_n(t)-f_0(t)|
=
O_{p}\left(\frac{1}{\sqrt {nh}}\right)+O(h^2)
=
o_p(n^{-1/3}).
\end{equation*}
The proof for $t\in[1-h,1]$ is completely similar, using the symmetry of $K$.
Together with~\eqref{eq:order difference f and f'}, this proves~\eqref{eq: approxtildef}
for boundary correction~\eqref{def:boundary kernel estimate}.

Similarly, for $t\in[0,h]$, write
\[
\tilde f_n'(t)-f_0'(t)
=
\tilde f_n'(t)-f_{0,B}'(t)
+
f_{0,B}'(t)-f_0'(t),
\]
where
\begin{equation}
\label{def:boundary kernel f0'}
f_{0,B}'(t)
=
\frac{\md}{\md t}
\int
\frac1{h}
\bK\left(\frac{t-x}{h}\right)
\,\md F_{0}(x).
\end{equation}
Note that
\[
\frac{\md}{\md t}
\bK\left(\frac{t-x}{h}\right)
=
\frac1h D_{B,t}\left(\frac{t-x}{h}\right),
\]
where $D_{B,t}:\R\mapsto\R$ is the bounded function defined by
\[
D_{B,t}(u)
=
\bcoefA'\left(\frac{t}{h}\right)K(u)
+
\bcoefB'\left(\frac{t}{h}\right)uK(u)
+
\bK'(u).
\]
Then, as before
\[
\begin{split}
\tilde f_n'(t)-f_{0,B}'(t)
&=-
\frac{1}{h^2}
D_{B,t}\left(\frac{t}{h}\right)
(F_{n0}(0)-F_0(0))\\
&\qquad+
\int_0^{t+h}
\frac{1}{h^3}
D_{B,t}'\left(\frac{t-x}{h}\right)
(F_{n0}(x)-F_0(x))
\,\md x,
\end{split}
\]
where $D_{B,t}'(u)=\md D_{B,t}(u)/\md u$.
Similar to~\eqref{eq:variance} we find
\begin{equation*}
\begin{split}
\tilde f_n'(t)-f_{0,B}'(t)
&=-
\frac{1}{\sqrt nh^2}
D_{B,t}\left(\frac{t}{h}\right)
B_{n0}\circ L_{0}(0)\\
&\quad+
\int_{-1}^{t/h}
\frac{1}{\sqrt nh^2}
D_{B,t}'(u)
B_{n0}\circ L_{0}(t-hu)
\,\md u
+
O_p(h^{-2}n^{-1+1/q}),
\end{split}
\end {equation*}
where the second integral on the right hand side is equal to
\[
\frac{1}{\sqrt nh^2}
D_{B,t}\left(\frac{t}{h}\right)B_{n0}\circ L_{0}(t)+O_{p}\left(\frac{1}{h\sqrt{nh}}\right).
\]
Together with~\eqref{eq:accB}, we conclude that
\begin{equation}
\label{eq:order term 1'}
\tilde f_n'(t)-f_{0,B}'(t)
=
O_{p}\left(\frac{1}{h\sqrt{nh}}\right).
\end{equation}
Next, consider $f_{0,B}'(t)-f_0'(t)$.
For $t\leq h$, define
\[
\bIK(s)=\int_{-1}^s\bK(u)\,\md u
=
\bcoefA\left(\frac{t}{h}\right)\IK^{(0)}(s)
+
\bcoefB\left(\frac{t}{h}\right)\IK^{(1)}(s)
\]
where $\IK^{(j)}$ is defined in~\eqref{def:int K}.
Then from integration by parts, from~\eqref{def:boundary kernel f0'} we get
\begin{equation}
\label{eq:decomposition f0B'}
\begin{split}
f_{0,B}'(t)
&=
\frac{\md}{\md t}
\left\{
f_0(0)
+
\int_0^{t+h}
\bIK\left(\frac{t-x}{h}\right)f_0'(x)\,\md x
\right\}\\
&=
\frac1h\int_0^{t+h} \bK\left(\frac{t-x}{h}\right)f_0'(x)\,\md x\\
&\qquad+
\frac1h
\int_0^{t+h}
\left\{
\bcoefA'\left(\frac{t}{h}\right)\IK^{(0)}\left(\frac{t-x}h\right)
+
\bcoefB'\left(\frac{t}{h}\right)\IK^{(1)}\left(\frac{t-x}h\right)
\right\}
f_0'(x)\,\md x.
\end{split}
\end{equation}
Similar to~\eqref{bias_expansion}, for the first term on the right hand side, we find
\begin{equation}
\label{bias_expansion'}
\frac1h\int_0^{t+h} \bK\left(\frac{t-x}{h}\right)f_0'(x)\,\md x\\
=
\int_{-1}^{t/h}
\bK(u)f_0'(t-hu)\,\md u
=
f_0'(t)+O(h),
\end{equation}
whereas the second term is equal to
\[
\begin{split}
&\int_{-1}^{t/h}
\left\{
\bcoefA'\left(\frac{t}{h}\right)\IK^{(0)}(u)
+
\bcoefB'\left(\frac{t}{h}\right)\IK^{(1)}(u)
\right\}
f_0'(t-hu)\,\md u\\
&=
f_0'(t)
\int_{-1}^{t/h}
\left\{
\bcoefA'\left(\frac{t}{h}\right)\IK^{(0)}(u)
+
\bcoefB'\left(\frac{t}{h}\right)\IK^{(1)}(u)
\right\}\,\md u
+
O(h).
\end{split}
\]
Note that for $s\in[-1,1]$,
\begin{equation}
\label{eq:derivative bK zero}
\int_{-1}^{s}
\left\{
\bcoefA'(s)\IK^{(0)}(u)
+
\bcoefB'(s)\IK^{(1)}(u)
\right\}\,\md u
=0.
\end{equation}
This can be seen as follows.
From~\eqref{eq:def coef boundary kernel} it follows that
\[
\begin{split}
&\int_{-1}^s
\left\{
\bcoefA(s)\IK^{(0)}(u)+\bcoefB(s)\IK^{(1)}(u)
\right\}
\,\md u\\
&\qquad=
s
\left\{
\bcoefA(s)\IK^{(0)}(s)+\bcoefB(s)\IK^{(1)}(s)
\right\}
-
\int_{-1}^s
\left\{
\bcoefA(s)uK(u)+\bcoefB(s)u^2K(u)
\right\}
\,\md u=s.
\end{split}
\]
Differentiating this equality with respect to $s$, gives
\[
\int_{-1}^s
\left\{
\bcoefA'(s)\IK^{(0)}(u)+\bcoefB'(s)\IK^{(1)}(u)
\right\}
\,\md u
+
\bcoefA(s)\IK^{(0)}(s)+\bcoefB(s)\IK^{(1)}(s)
=
1,
\]
which proves~\eqref{eq:derivative bK zero}, due to~\eqref{eq:def coef boundary kernel}.
Hence, from~\eqref{eq:decomposition f0B'}, \eqref{bias_expansion'} and~\eqref{eq:order term 1'},
and using the fact that $1/6<\gamma\leq 1/5$,
we conclude that
\begin{equation*}
\sup_{t\in[0,h]}
|\tilde f_n'(t)-f_0'(t)|
=
O_{p}\left(\frac{1}{h\sqrt{nh}}\right)
+
O(h)
=
o_p(n^{-1/6}).
\end{equation*}
The proof for $t\in[1-h,1]$ is completely similar.
Together with~\eqref{eq:order difference f and f'}, this proves~\eqref{eq: approxtildef'}
for boundary correction~\eqref{def:boundary kernel estimate}.

Finally, consider the second derivative
\[
\tilde f_n''(t)-f_0''(t)
=
\tilde f_n''(t)-f_{0,B}''(t)
+
f_{0,B}''(t)-f_0''(t),
\]
where
\[
f_{0,B}''(t)
=
\frac{\md^2}{\md t^2}
\int
\frac1{h}
\bK\left(\frac{t-x}{h}\right)
\,\md F_{0}(x)
=
\frac1{h^2}
\int
\frac{\md}{\md t}
D_{B,t}\left(\frac{t-x}{h}\right)
\,\md F_{0}(x).
\]
Similar to~\eqref{eq:order term 1} and~\eqref{eq:order term 1'},
we obtain
\begin{equation}
\label{eq:order term 1''}
\tilde f_n''(t)-f_{0,B}''(t)
=
O_{p}\left(\frac{1}{h^2\sqrt{nh}}\right).
\end{equation}
By differentiating the right hand side of~\eqref{eq:decomposition f0B'}, we find
\begin{equation}
\label{eq:decomposition f0B''}
\begin{split}
f_{0,B}''(t)
=&
\frac{\md}{\md t}
\int_0^{t+h}\frac1h \bK\left(\frac{t-x}{h}\right)f_0'(x)\,\md x\\
&+
\frac{\md}{\md t}
\int_0^{t+h}
\frac1h\left\{
\bcoefA'\left(\frac{t}{h}\right)\IK^{(0)}\left(\frac{t-x}h\right)
+
\bcoefB'\left(\frac{t}{h}\right)\IK^{(1)}\left(\frac{t-x}h\right)
\right\}
f_0'(x)\,\md x.
\end{split}
\end{equation}
For the first term on the right hand side of~\eqref{eq:decomposition f0B''} we get
similar to~\eqref{eq:decomposition f0B'}:
\[
\begin{split}
&
\frac{\md}{\md t}
\int_0^{t+h}\frac1h \bK\left(\frac{t-x}{h}\right)f_0'(x)\,\md x\\
&=
\frac{\md}{\md t}
\left\{
f_0'(0)
+
\int_0^{t+h}
\bIK\left(\frac{t-x}{h}\right)f_0''(x)\,\md x
\right\}\\
&=
\frac1h\int_0^{t+h} \bK\left(\frac{t-x}{h}\right)f_0''(x)\,\md x\\
&\qquad+
\frac1h
\int_0^{t+h}
\left\{
\bcoefA'\left(\frac{t}{h}\right)\IK^{(0)}\left(\frac{t-x}h\right)
+
\bcoefB'\left(\frac{t}{h}\right)\IK^{(1)}\left(\frac{t-x}h\right)
\right\}
f_0''(x)\,\md x,
\end{split}
\]
which is bounded on $[0,h]$, because $f_0''$ is bounded.
For the second term on the right hand side of~\eqref{eq:decomposition f0B''},
first note that from integration by parts
\[
\begin{split}
&
\int_0^{t+h}
\frac1h\left\{
\bcoefA'\left(\frac{t}{h}\right)\IK^{(0)}\left(\frac{t-x}h\right)
+
\bcoefB'\left(\frac{t}{h}\right)\IK^{(1)}\left(\frac{t-x}h\right)
\right\}
f_0'(x)\,\md x\\
&=
\left\{
\bcoefA'\left(\frac{t}{h}\right)\IK^{(0,\mathrm{int})}\left(\frac{t}h\right)
+
\bcoefB'\left(\frac{t}{h}\right)\IK^{(1,\mathrm{int})}\left(\frac{t}h\right)
\right\}
f_0'(0)\\
&\qquad+
\int_0^{t+h}
\left\{
\bcoefA'\left(\frac{t}{h}\right)\IK^{(0,\mathrm{int})}\left(\frac{t-x}h\right)
+
\bcoefB'\left(\frac{t}{h}\right)\IK^{(1,\mathrm{int})}\left(\frac{t-x}h\right)
\right\}
f_0''(x)\,\md x\\
&=
\int_0^{t+h}
\left\{
\bcoefA'\left(\frac{t}{h}\right)\IK^{(0,\mathrm{int})}\left(\frac{t-x}h\right)
+
\bcoefB'\left(\frac{t}{h}\right)\IK^{(1,\mathrm{int})}\left(\frac{t-x}h\right)
\right\}
f_0''(x)\,\md x,
\end{split}
\]
due to~\eqref{eq:derivative bK zero}, where for $j=0,1$,
\[
\IK^{(j,\mathrm{int})}(s)
=
\int_{-1}^s\IK^{(j)}(u)\,\md u.
\]
This means that
\[
\begin{split}
&\frac{\md}{\md t}
\int_0^{t+h}
\frac1h\left\{
\bcoefA'\left(\frac{t}{h}\right)\IK^{(0)}\left(\frac{t-x}h\right)
+
\bcoefB'\left(\frac{t}{h}\right)\IK^{(1)}\left(\frac{t-x}h\right)
\right\}
f_0'(x)\,\md x\\
&=
\frac1{h}
\int_0^{t+h}
\left\{
\bcoefA'\left(\frac{t}{h}\right)\IK^{(0)}\left(\frac{t-x}h\right)
+
\bcoefB'\left(\frac{t}{h}\right)\IK^{(1)}\left(\frac{t-x}h\right)
\right\}
f_0''(x)\,\md x\\
&\qquad+
\frac1{h}
\int_0^{t+h}
\left\{
\bcoefA''\left(\frac{t}{h}\right)\IK^{(0,\mathrm{int})}\left(\frac{t-x}h\right)
+
\bcoefB''\left(\frac{t}{h}\right)\IK^{(1,\mathrm{int})}\left(\frac{t-x}h\right)
\right\}
f_0''(x)\,\md x,
\end{split}
\]
which is bounded uniformly for $t\in[0,h]$.
Together with~\eqref{eq:order term 1''}, using that $1/6<\gamma\leq 1/5$,
it follows that
\begin{equation*}
\sup_{t\in[0,h]}
|\tilde f_n''(t)-f_0''(t)|
=
O_{p}\left(\frac{1}{h^2\sqrt{nh}}\right)
+
O(1)
=
O_p(1).
\end{equation*}
The proof for $t\in[1-h,1]$ is completely similar.
Together with~\eqref{eq:order difference f and f'}, this establishes~\eqref{eq:order fn''}
and proves~\eqref{eq: tildef'} for boundary correction~\eqref{def:boundary kernel estimate}.
\tqed

\end{document}